\documentclass{amsart}
\usepackage{latexsym, amsmath, amsfonts, amscd, amssymb}

\newtheorem{theorem}{Theorem}

\newcommand{\PA}{\par \hskip0.4cm}
\newtheorem{proposition}{Proposition}
\newtheorem{notat}{\bf Notation}
{\end{notat}}

\newtheorem{clai}{\bf Claim}
\newenvironment{claim}{\begin{clai}- \rm }%
{\end{clai}}
\newtheorem{exemp}{\bf Example}
{\end{exemp}}
\newtheorem{lemma}{Lemma}

\begin{document}
\begin{abstract}
The goal of this work is to find the asymptotics of the hitting
probability of a distant point for the voter model on the integer
lattice started from a single $1$ at the origin. In dimensions $d=2$
or $3$, we obtain the precise asymptotic behavior of this
probability. We use the scaling limit of the voter model started
from a single $1$ at the origin in terms of super-Brownian motion
under its excursion measure. This invariance principle was stated by
Bramson, Cox and Le Gall, as a consequence of a theorem of Cox,
Durrett and Perkins. Less precise estimates are derived in dimension
$d \geq 4$.
\end{abstract}
\begin{center}
 {\Large \bf{HITTING PROBABILITY OF A DISTANT POINT FOR THE VOTER MODEL STARTED WITH A SINGLE ONE}}
\vskip0.8cm
{\large Mathieu {\sc Merle}}
\vspace{0.2cm}
\end{center}
{ \small \paragraph{\sc{Abstract}}
The goal of this work is to find the asymptotics of the hitting probability
of a distant point for the voter model on the integer lattice started from a single $1$ at the origin.
In dimensions $d=2$ or $3$, we obtain the precise asymptotic
behaviour of this probability. We use the scaling limit of the voter model
started from a single
$1$ at the origin in terms of super-Brownian motion under its excursion
measure. This invariance principle
was stated by Bramson, Cox and Le Gall, as a consequence of a theorem
of Cox, Durrett and Perkins. Less precise estimates are derived in dimension $d \geq 4$.}
 \vspace{0.4cm} \normalfont
\section{Introduction, Notation and Statement of result}
\PA The voter model is one of the most classical
 interacting particle systems. This model is of great interest because it
exhibits a range of interesting phenomena and also because
it is dual to a system of coalescing random walks.
The voter model was first introduced in \cite{CS}, \cite{HL}, and some of
its basic properties were investigated by Liggett \cite{Li}, Sawyer \cite{Sa},
Arratia \cite{A2}, Bramson and Griffeath \cite{BG}. \PA
More recently, Cox, Durrett and
Perkins \cite{CDP} showed an important invariance principle,
establishing that, after a suitable renormalization, voter models in dimension
$d \geq 2$ converge to super-Brownian motion. Super-Brownian motion
is a continuous measure-valued
process which arises as the weak limit of branching particle systems
(see Watanabe \cite{W}). It was discussed by Dawson
\cite{D1}, and studied extensively in the nineties (see in particular
 \cite{D2}, \cite{Pe}, \cite{LG}).
In the recent years, it was shown that super-Brownian motion also appears in
scaling limits of a wide range of lattice systems such as lattice trees,
contact processes or oriented percolation.
The main idea of this work is to exploit known properties of super-Brownian
motion to get asymptotic results for the voter model.

\PA Let us now describe the voter model and state our main result.
 Let $d \geq 2$. At each site of the integer lattice
$\mathbb{Z}^{d}$ there is a voter holding an opinion. We will study here
a two-type model, where
there are only two possible opinions, say $0$ or $1$.
At rate $1$ exponential times, the voter at $x \in \mathbb{Z}^{d}$
chooses a neighbor $y$ according to a given jump kernel
$p$ and adopts the opinion of $y$.
The voting times and neighbor selections are supposed independent.
The jump kernel $p : \mathbb{Z}^{d}
\times \mathbb{Z}^{d} \to [0,1]$ will be supposed symmetric,
translation invariant, irreducible, centered, isotropic, and
having exponential moments :
\begin{itemize}
\item $p(x,y) = p(0,y-x), \qquad p(x,y)= p(y,x), \qquad p(0,0)=0$,
\item $\sum_{y\in \mathbb{Z}^{d}} y p(0,y) =0$,
\item $\sum_{y\in \mathbb{Z}^{d}}
p(0,y)y^{i}y^{j}= \sigma^{2}\delta_{ij}$ for some $0<\sigma^{2}<\infty$,
\item there exists a constant $C>0$ such that $\sum_{y \in \mathbb{Z}^{d}}
p(0,y) \exp(C|y|) < \infty$.
 \end{itemize}
If $t \geq 0$, we denote by $\xi_{t}$ the set of sites where voters
hold opinion $1$ at time $t$; $(\xi_{t})_{t \geq 0}$ is the two-type
voter model. If $A \subset \mathbb{Z}^{d}$, we write $P_{A}$ for the
probability measure under which $\xi_{0}=A$.
 Throughout this paper, we will consider the particular case
when $\xi_{0}=\{0\}$. In this case,
$(\xi_{t}^{0})_{t \geq 0}$ will denote the two-type voter model
started from a single opinion $1$ at the origin, and
for simplicity, we will write $P$ for $P_{\{0\}}$.  \PA
 It is often convenient to work with the associated measure-valued
processes
$$X_{t} : = \sum_{y \in \xi_{t}} \delta_{y}, \quad
X_{t}^{0} : = \sum_{y \in \xi_{t}^{0}} \delta_{y}.$$
 For $\alpha >0$ we define the conditional probability
 $$ P_{\alpha}^{*}(.) : = P(.|\xi_{\alpha}^{0} \neq \emptyset). $$
 \vskip0.1cm \PA
  We are interested in estimating
 the probability that a voter located at a distance
of order $c$ from the origin ever holds opinion $1$.
If $x \in \mathbb{R}^{d}$, we denote by $[x]_{c}$ the point in
$c^{-1}\mathbb{Z}^{d}$ closest to $x$. If there is more than one such point, we
choose the point closest to the origin.
 Our goal is to find the
asymptotic order as $c \to \infty$ of
$$P( \exists t \geq 0 : c[x]_{c} \in \xi_{t}^{0}).$$
We introduce the notation $T_{c[x]_{c}} = \inf \{ t \geq 0 : c[x]_{c}
\in \xi_{t}^{0} \}$ so that the previous quantity can also be written
$P(T_{c[x]_{c}}<\infty)$. Set $\beta_{2}=2\pi$, and for
$d \geq 3$, let $\beta_{d}$ be the probability that a rate $1$ continuous time
random walk with jump kernel $p$ started from the origin
 never returns to it.
\begin{theorem} \label{main}
Let $x \in \mathbb{R}^{d}\setminus 0 $ be fixed.
Let us define
$$ \phi_{d}(c) =
\begin{cases} \frac{c^{2}}{2\ln(c)} & \mbox{ if } d =2, \\
              c^{2} & \mbox{ if } d=3, \\
 c^{d-2} & \mbox{ if } d \geq 5.  \end{cases}
$$
Then, if $d=2$ or $d=3$,
\begin{eqnarray*}
\lim_{c \to \infty} \phi_{d}(c) P(T_{c[x]_{c}} < \infty ) =
\frac{2\sigma^{2}}{\beta_{d}}\left(2-
\frac{d}{2}\right) |x|^{-2}.
\end{eqnarray*}
If $d \geq 5$, there exist positive constants $a_{d}, b_{d}$
depending on $x$ such that
\begin{eqnarray*}
a_{d} \leq \liminf_{c \to \infty} \phi_{d}(c) P(T_{c[x]_{c}} < \infty )  \leq
\limsup_{c \to \infty} \phi_{d}(c) P(T_{c[x]_{c}} < \infty ) \leq b_{d}.
\end{eqnarray*}
\end{theorem}
\PA In dimension $4$ we obtain less precise results. We will
prove the existence of a positive constant $a_{4}$ and we
conjecture the existence of a positive $b_{4}$ such that a
statement similar to the one in $d \geq 5$ holds for $d=4$ with
the function $\phi_{4}(c): = c^{2}\ln(c)$. The upper bound in dimension $4$
seems more difficult than the corresponding results in other dimensions.
 Adapting the proof of the upper bound
for $2 \leq d \leq 3$ to the case $d=4$ only gives
$\limsup_{c \to \infty} c^{2} P(T_{c[x]} < \infty) = 0$.
\PA Theorem \ref{main} immediately extends to the multitype voter model
$\overline{\xi}_{t}$, which is described as follows. We assume that the
initial opinions are all distinct.
The dynamics of the multitype voter model are
the same as those of the two-type voter model. In this multitype setting,
Theorem \ref{main} gives the asymptotics of the probability that the voter
at $x$ ever adopts the initial opinion of $y$, as $|x-y|$ tends to infinity.
\vskip0.2cm
 \PA In dimensions $2$ and $3$, we will let $T>0$ and argue under the
measure $P_{c^{2}T}^{*}$. Motivated by the results of \cite{CDP},
 Bramson, Cox and Le
Gall \cite{BCLG} proved
that for $T>0$, the voter model $\xi^{0}$ under $P_{c^{2}T}^{*}$
 converges as $c \to \infty$ modulo a suitable rescaling
to a nondegenerate limit that can be expressed in terms of the excursion
measure $\mathbb{N}_{0}$
of super-Brownian motion (see Theorem \ref{inv} below). This invariance
principle of \cite{BCLG} will be our main tool in the proof of Theorem
\ref{main} for small dimensions. We will also need properties of super-Brownian
motion under its excursion measure $\mathbb{N}_{0}$.
The Brownian snake
 approach of Le Gall \cite{LG} gives a good understanding of the measure
$\mathbb{N}_{0}$, and will be used to prove an intermediate result.

\PA As mentioned earlier, the voter model and coalescing random walks are
dual processes. In a
coalescing random walk system, particles are assumed to execute rate $1$
random walks with jump kernel $p$. Particles move independently until they
meet, then coalesce and move together afterwards. The duality
property also serves as a major tool for our results.
\vskip0.2cm
\PA In Section 2.1, we introduce super-Brownian motion and its excursion
measure $\mathbb{N}_{0}$.
 Scaling limits of the voter model
(invariance principles) are discussed in Section 2.2.
The duality property
is explained in Section 2.3, and preliminary
results on rate $1$ random walks and system of coalescing random
walks are discussed in Section 2.4
 and 2.5.\PA
 We establish
the asymptotic upper bounds on $P(T_{c[x]}<\infty)$ in Section 3. This requires
 interesting intermediate results. Lemma \ref{speed} expresses
that the probability for the voter model under $P_{\alpha}^{*}$ to escape
$B(0,A)$ before time $2 \alpha$
 decays exponentially with $A$. Lemma \ref{enough} informally
expresses that for any fixed $\epsilon >0$,
then,
$\cup_{t \geq \epsilon \alpha} \xi_{t}^{0}$
does not contain any ``isolated'' point, with arbitrarily high probability
under $P_{\alpha}^{*}$, when $\alpha$ is taken large enough. \PA
We prove the asymptotic lower bounds in Section 4. Sections 4.1
is devoted to the case $d \geq 4$, and Sections 4.2 and 4.3 to the case
$d=2$ or $3$.
Finally, we prove the results of Sections 2.4 and 2.5 in Section 5.

\section{Further notation and preliminary results}
Let $f$ and $g$ be two functions from $\mathbb{R}$ into $(0,\infty)$.
 We will write $f(x) = o(g(x))$ as $x \to \infty$, respectively
$f(x) \sim g(x)$ as $x \to \infty$ whenever $\lim_{x \to \infty}
f(x)(g(x))^{-1}$ is equal to $0$, respectively $1$. \PA For $x \in
\mathbb{R}^{d}, r>0$ we denote by $B(x,r)$ the open ball in
$\mathbb{R}^{d}$ centered at $x$ with radius $r$, and $B(x,r)^{c}$
its complement. \PA For real numbers $x \leq y$, the set $\{n \in
\mathbb{Z} : x \leq n \leq y \}$ of integers between $x$ and $y$
will be denoted by $\big[ \! |x,y| \!\big]$; also, the integer part
of $x$ : $\max\{ n \in \mathbb{Z} : n \leq x \}$ will be denoted by
$\left\lfloor x \right\rfloor$, while $\left\lfloor x \right\rfloor
+1 = \min\{n \in \mathbb{Z} : n > x \}$
 will be denoted by $\left\lceil x \right\rceil$. \PA
\subsection{Super-Brownian motion}
Let $M_{F}(\mathbb{R}^{d})$ be the space of all
finite measures on $\mathbb{R}^{d}$,
equipped with the topology of weak convergence. For
$\mu \in M_{F}(\mathbb{R}^{d})$, $f$ a function on $\mathbb{R}^{d}$, the
notation $\mu(f)$ will stand for $\int_{\mathbb{R}^{d}} f(x)\mu(dx)$ whenever
this integral is well-defined. We let  $\mathcal{C}(\mathbb{R}_{+},
M_{F}(\mathbb{R}^{d}))$ be the space
 of continuous paths from $\mathbb{R}_{+}$ into
$M_{F}(\mathbb{R}^{d})$, and we let
$D(\mathbb{R}_{+}, M_{F}(\mathbb{R}^{d}))$ be the Skorohod space of
cadlag functions from $\mathbb{R}^{d}$ into $M_{F}(\mathbb{R}^{d})$.
We denote by $(Y_{t}, t \geq 0)$ the canonical process on either
$\mathcal{C}(\mathbb{R}_{+}, M_{F}(\mathbb{R}^{d}))$ or
$D(\mathbb{R}_{+}, M_{F}(\mathbb{R}^{d}))$.
\PA The law of
 super-Brownian motion
with branching rate $\gamma$ and
diffusion coefficient $\sigma^{2}$, starting from
$\mu \in M_{F}(\mathbb{R}^{d})$, is the probability measure
$\mathbf{Q}_{\mu}^{\gamma,\sigma^{2}}$ on
$\mathcal{C}(\mathbb{R}_{+}, M_{F}(\mathbb{R}^{d}))$ that solves the
following well-posed martingale problem (see \cite{Pe}, Theorem II.5.1) :
\begin{enumerate}
\item[(MP)]  For any $\phi\in C_{b}^{\infty}(\mathbb{R}^{d})$,
\begin{equation*}
Y_{t}(\phi)\ = \  \mu(\phi) + M_{t}(\phi) + \frac{1}{2} \int_{0}^{t}
Y_{s}(\sigma^{2}\Delta \phi) ds,
\end{equation*}
where $M_{t}(\phi)$ is a $\mathbf{Q}_{\mu}^{\gamma,\sigma^{2}}$ continuous
square integrable martingale such
that $M_{0}(\phi)=0$ and the quadratic variation of $M(\phi)$ is
   $$<M(\phi)>_{t}=\int_{0}^{t} Y_{s}(\gamma \phi^{2}) ds. $$
\end{enumerate}
One can show
(see for example Section II.7 of \cite{Pe}) that there exists a family
$\{ R_{t}^{\gamma,\sigma^{2}}(y,.), y \in \mathbb{R}^{d}, t>0 \}$
of finite measures on
$M_{F}(\mathbb{R}^{d})$, called the canonical measures of super-Brownian motion,
which assign zero mass to the $0$ measure and are such that the following holds.
The law of $Y_{t}$ under $\mathbf{Q}_{\mu}^{\gamma,\sigma^{2}}$ is the same as
the law of $\sum_{i \in I} Y_{t}^{i}$, where
$\sum_{i \in I} \delta_{Y_{t}^{i}}$ is a Poisson measure with intensity
$\int_{M_{F}(\mathbb{R}^{d})} R_{t}^{\gamma,\sigma^{2}}(y,.)\mu(dy)$.
It follows that for any Borel subset $\mathcal{Y}$ of $M_{F}(\mathbb{R}^{d})$
 with $0 \notin \mathcal{Y}$,
\begin{equation} \label{convR}
 \lim_{\epsilon \to 0} \epsilon^{-1}
\mathbf{Q}_{\epsilon \delta_{y}}^{\gamma,\sigma^{2}}(Y_{t}\in \mathcal{Y})
=R_{t}^{\gamma,\sigma^{2}}(y,\mathcal{Y}).
\end{equation}
It is also well-known (see \cite{Pe}, Theorem II.7.2)
 that for any $y \in \mathbb{R}^{d}$,
\begin{equation} \label{totalmass}
R_{t}^{\gamma,\sigma^{2}}(y,M_{F}(\mathbb{R}^{d}))= \frac{2}{\gamma t}.
\end{equation}

From \cite{Pe}, Theorem II.7.3 (see also formula (3.10) in \cite{BCLG}),
 for each $y \in \mathbb{R}^{d}$ there is a
$\sigma$-finite measure $\mathbb{N}_{y}$ on $\mathcal{C}(\mathbb{R}_{+},
M_{F}(\mathbb{R}^{d}))$ called the excursion measure of
super-Brownian motion with branching rate $\gamma$ and diffusion coefficient
$\sigma^{2}$ such that the following holds. For any $\alpha >0$ fixed, then
for any bounded continuous function $F$ on
$\mathcal{C}(\mathbb{R}_{+}, M_{F}(\mathbb{R}^{d}))$, such that $F(\omega)=0$
for any $\omega$ with $\omega(t)=0$ for all $t \geq \alpha$,
\begin{equation} \label{convNy}
\lim_{\epsilon \to 0} \epsilon^{-1}
\mathbf{Q}_{\epsilon \delta_{y}}^{\gamma,\sigma^{2}}(F((Y_{t},t \geq 0))) =
\mathbb{N}_{y}(F).
\end{equation}
The convergence (\ref{convR}) is a particular case of (\ref{convNy}).
Thus,
for any Borel subset $\mathcal{Y}$ of  $M_{F}(\mathbb{R}^{d})$ with
$0 \notin \mathcal{Y}$, we have
$$ \mathbb{N}_{y}(Y_{t} \in \mathcal{Y}) =
R_{t}^{\gamma, \sigma^{2}}(y,\mathcal{Y}).$$
Also, for any $T >0, y \in \mathbb{R}^{d}$, we get from (\ref{totalmass})
\begin{equation} \label{surviv}
 \mathbb{N}_{y}(Y_{T} \neq 0) = \frac{2}{\gamma T},
\end{equation}
and we can define the probability measure
 $\mathbb{N}_{y}^{(T)} : = \mathbb{N}_{y}(. | Y_{T} \neq 0)$. \PA
A better understanding of the measures $\mathbb{N}_{y}$ is given by
 the Brownian snake approach of Le Gall
(see \cite{LG}, and Section 4.4 below). The Brownian snake approach
corresponds to $\gamma =4$, but scaling properties of super-Brownian motion
can then be used to deal with a general value of $\gamma$. \PA
 Finally, we will use the following result about hitting
probabilities of a single point. Let
$\mathcal{R}_{t}$ denote the topological
 support of the measure $Y_{t}$,
and $\mathcal{R} = \bigcup_{t > 0} \mathcal{R}_{t}$.
It follows from \cite{LG}, Section 6.1 that
\begin{equation} \label{N0x}
\mathbb{N}_{0}(x \in \mathcal{R}) =
\frac{4\sigma^{2}}{\gamma} \left( 2-\frac{d}{2} \right)^{+} |x|^{-2}.
\end{equation}
In particular, in the case $d \geq 4$,
$\mathbb{N}_{0}(x \in \mathcal{R})=0 $,
which explains why our results are less precise.
 Also, as (\ref{N0x})
suggests, the case of dimension $4$ is critical, and thus
harder.  \PA
  Since $\mathbb{N}_{0}(Y_{T}=0 \big| x \in \mathcal{R}) \to 0$ as $T$ goes
to $0$, we deduce from (\ref{surviv}) and (\ref{N0x}) that
\begin{equation} \label{N0xT}
\mathbb{N}_{0}^{(T)}(x \in \mathcal{R}) \underset{T \to 0}{\sim}
2 \sigma^{2} T \left( 2-\frac{d}{2} \right)^{+} |x|^{-2}.
\end{equation}

\subsection{Extinction probability, invariance principle}
 Set $p_{t}:=P(\xi_{t}^{0} \neq \emptyset)$. The asymptotic rate at which
$p_{t}$ converges to $0$ was found in \cite{BG}. As $t \to \infty$,
\begin{eqnarray} \label{pt}
 p_{t} \sim \begin{cases} \log(t)/(\beta_{2}t)
 \mbox{ if } d=2 \\ 1/(\beta_{d}t) \mbox{ if } d \geq 3, \end{cases}
\end{eqnarray}
where $\beta_{d}, d\geq 2$ was defined before Theorem \ref{main}.
Hence, for any $d \geq 2$ there exist a positive $\kappa_{0}$
depending only on $d$ such that for any $1/4 < t' \leq t$,
\begin{equation} \label{asy}
\frac{p_{t'}}{p_{t}} \leq \kappa_{0} \frac{t}{t'}
\end{equation}
\vskip0.1cm
 If $|C|$ denote the cardinality
of a finite set $C$, Bramson and Griffeath (\cite{BG}) established
that the law of $p_{t}|\xi_{t}^{0}|$ under $P_{t}^{*}$ converges as
$t \to \infty$ to
 an exponential distribution with parameter $1$.
\PA
 Bramson and Griffeath \cite{BG} also conjectured that
$\xi_{t}^{0}$ would obey a certain asymptotic shape theorem. Such a result
 was derived in $2001$ by Bramson, Cox and Le Gall
\cite{BCLG} using the invariance principle relating the voter
model and super-Brownian motion, which was proved by Cox, Durrett
and Perkins in \cite{CDP}.
We rescale the voter model as follows. For $N>0$, the lattice is
 now $S_{N} := \mathbb{Z}^{d}/\sqrt{N}$.
Individuals change opinion at rate $N$ instead of $1$,
and the jump kernel becomes $p_{N} : S_{N}
\times S_{N} \to \mathbb{R}_{+}$ such that $p_{N}(x,y) =
p(\sqrt{N}x,\sqrt{N}y)$. We denote by $(\xi_{t}^{N,0})_{t \geq 0}$
the corresponding process ($\xi_{t}^{N,0}$ represents the set of
sites having opinion $1$ at time $t$). If we let
$$m_{N} := \frac{N}{\log(N)} \mbox{ if } d=2, \qquad m_{N} := N \mbox{ if }
d \geq 3,$$ we can define an associated measure-valued processes :
 $$ X_{t}^{N,0} :=
\frac{1}{m_{N}} \sum_{y \in \xi_{t}^{N,0}} \delta_{y}. $$ \PA
Similarly, when at time $0$, opinion $1$ is started from a given set
$\xi_{0}$, we may define for $N>0$ a rescaled voter model
$\xi_{t}^{N}$ and the corresponding measure valued process
$X_{t}^{N}$. Theorem 1.2 of \cite{CDP} states that whenever
$X_{0}^{N}$ converges to a non-degenerate measure $X_{0} \in
M_{F}(\mathbb{R}^{d})$,
 then $(X_{t}^{N})_{t \geq 0}$
 converges to a super-Brownian motion on $\mathbb{R}^{d}$ with branching
rate $2 \beta_{d}$ and diffusion coefficient $\sigma^{2}$, started from
$X_{0}$.  \PA
Theorem 2 below states the convergence in law of the process
$(X_{t}^{N,0})_{t \geq 0}$ under the conditional distribution
$P(.\big|X_{\alpha}^{N,0} \neq 0)$
towards super-Brownian motion under $\mathbb{N}_{0}^{(\alpha)}$.
 This result, which is taken from \cite{BCLG} (Theorem 4) will
be a key ingredient of the proof of Theorem
\ref{main} in dimensions $2$ and $3$.
\begin{theorem} \label{inv}
Assume $d \geq 2$, and let $\mathbb{N}_{0}$ be the excursion measure of
super-Brownian motion on $\mathbb{R}^{d}$ with branching rate $2\beta_{d}$
and
diffusion coefficient $\sigma^{2}$. Let $\alpha >0$, and let $F$ be a bounded
continuous function on $D(\mathbb{R}_{+},M_{F}(\mathbb{R}^{d}))$.
Then
\begin{equation} \label{convtoN0}
\lim_{N \to \infty} E \left[ F\left( (X_{t}^{N,0})_{t \geq 0}\right)
\big| X_{\alpha}^{N,0} \neq 0 \right] = \mathbb{N}_{0}^{(\alpha)}[F].
\end{equation}
\end{theorem}

\PA Let us now turn to the well-known relation between the voter
model and coalescing random walks.
\subsection{Dual process to the voter model}
 Let us introduce further notation in order to
describe the dual process to the voter model.
The times at which the voter at $x$ adopts the opinion of the voter at
$y$ are the jump times of a standard Poisson process with rate $p(x,y)$. We
 denote by  $\Lambda(x,y)$ this set of times.
 Then,
$\{ \Lambda(x,y), x,y \in \mathbb{Z}^{d} \} $ forms a family of
independent Poisson point processes on $[0,\infty)$. \PA We now
describe the useful graphical representation of the voter model.
Horizontal axis represents $\mathbb{Z}^{d}$, vertical axis
represents time. For $x,y \in \mathbb{Z}^{d}$ we draw a horizontal
arrow from $y$ to $x$ at each time $s \in \Lambda(x,y)$. \PA
  For $s<t$ we say there is a path up from $(y,s)$ to $(x,t)$ or equivalently
a path down from $(x,t)$ to $(y,s)$ and we will
write $$ (y,s) \nearrow (x,t) \Leftrightarrow (x,t) \searrow (y,s)$$
if there exist times
$s=s_{0}< s_{1} <... < s_{n} \leq s_{n+1}=t$ and sites
$y=x_{0}, x_{1},...,x_{n}=x$ such that
\begin{itemize}
\item for $1\leq i \leq n $ there is an arrow pointing from $x_{i-1}$ towards
$x_{i}$ at time $s_{i}$,
\item for $0 \leq i \leq n$, there is no arrow pointing towards $x_{i}$ in
the time interval $(s_{i}, s_{i+1})$.
\end{itemize}
Clearly for every $x \in \mathbb{Z}^{d}$ and every choice of $0 \leq s \leq t$,
 there is a unique $y \in \mathbb{Z}^{d}$ such that $(y,s) \nearrow (x,t)$.
In such a case, the opinion of $(x,t)$ is "descended" from that at
$(y,s)$. We will say that $x$ at time $t$ is a "descendant" of $y$
at time $s$, or equivalently that $y$ at time $s$ is an "ancestor"
of $x$ at time $t$.
 \PA We are now in a position to describe the dual process to the voter model.
For $t>0$ and $x \in \mathbb{Z}^{d}$ we define
$(Z_{s}^{x,t})_{0 \leq s \leq t}$ by setting $Z_{0}^{x,t}=x$ and for
$0 < s \leq t $, $Z_{s}^{x,t}=y$ if and only if
$(x,t) \searrow (y,t-s)$. Clearly, $(Z_{s}^{x,t})_{0 \leq s \leq t}$ is a
rate $1$ random walk with jump kernel $p$ starting from $x$. Moreover, for
$x \in \mathbb{Z}^{d}$, $y \in \mathbb{Z}^{d}$, the two walks
$(Z_{s}^{x,t})_{0 \leq s \leq t}, (Z_{s}^{y,t})_{0 \leq s \leq t}$ start
respectively from $x$ and $y$, move
independently until they meet, and move together afterwards. That is,
 $(Z_{s}^{x,t})_{0 \leq s \leq t, x \in \mathbb{Z}^{d}}$
forms a coalescing random walk system with jump kernel
$p$. Furthermore \begin{equation} \label{dual}
 \xi_{t}^{0} = \{ y
\in \mathbb{Z}^{d} : Z_{t}^{y,t} =0 \}. \end{equation}
   For $t \geq 0$, we denote by $\hat{\xi}_{s}^{y,t}$ the set of
 descendants at time $t+s$ of $y$ at time $t$, that is
 $$  \hat{\xi}_{s}^{y,t} : =
\{ z \in \mathbb{Z}^{d} : (y,t) \nearrow (z,t+s) \}
 = \{ z \in \mathbb{Z}^{d}  :
Z_{s}^{z,t+s} = y \} $$
Notice that
$(\hat{\xi}_{s}^{y,t})_{s \geq 0}$ has the same law as
$(\xi_{s}^{0}+y)_{s \geq 0}$.
For $u \leq t$ we will denote by $\Omega_{u}^{t}$ the set of points having
opinion $1$ at time $u$ and having descendants  at time $t$, that is
$$ \Omega_{u}^{t} : =
\{y \in \xi_{u}^{0} : \hat{\xi}_{t-u}^{y,u} \neq \emptyset \} =
 \xi_{u}^{0} \cap \{Z_{t-u}^{z,t}, z \in
\mathbb{Z}^{d} \}. $$

\vskip0.2cm \hskip0.4cm The coalescing random walk perspective,
combined with the Bramson and Griffeath results and Theorem
\ref{inv}, gives us a heuristic explanation of our main result
Theorem \ref{main}. If $c[x]_{c}$ has opinion $1$ at time $t$, then
$Z_{t}^{c[x]_{c},t}=0$ so that from well-known properties of random
walks, $t$ should be of order $c^{2}$. The probability for the voter
model to survive a time of order $c^{2}$ is of order $p_{c^{2}}$,
and conditionally on that event, the rescaled voter model converges
to super-Brownian motion under its excursion measure. Informally,
formula (\ref{N0x}) is then exactly what we need to conclude in the
case $2 \leq d \leq 3$. Also, not rigourously, one should expect
that for $d \geq 4$, the probability of hitting $c[x]_{c}$ should be
of order $p_{c^{2}} \times \mathbb{N}_{0}(Y \mbox{ hits } B(x,1/c))
\approx A_{d}(x) \times \phi_{d}(c)^{-1}$,
 where $A_{d}(x)$ is a constant depending only on $d$ and $|x|$ (see \cite{DIP}). \PA
In the following paragraph, we present a few well-known properties
of random walks, then some estimates for coalescing random walks.
These will prove useful when using the duality property in the
course of the proof of Theorem \ref{main}.

\subsection{Random walks with jump kernel $p$}
We denote by $(Z_{t}, t\geq 0)$ a continuous-time random walk on
$\mathbb{Z}^{d}$ with jump kernel $p$ and exponential holding times with
parameter $1$. For $x \in \mathbb{Z}^{d}$, $Z$ starts from $x$ under the
probability measure $\mathbb{P}_{x}$.
For $x,y \in \mathbb{Z}^{d}, t\geq 0$ we let
$$q_{t}(x,y) = q_{t}(y-x) := \mathbb{P}_{x}(Z_{t}=y)$$  be
the transition kernel of our random walk.
For $x,y \in \mathbb{R}^{d}$ and $t>0$ let
$$ p_{t}(x,y) = p_{t}(x-y) := (2 \pi \sigma^{2}t)^{-d/2}
\exp \left( -\frac{|x-y|^{2}}{2\sigma^{2}t} \right) $$ be the transition
density of $d$-dimensional Brownian motion.
We denote by $P_{t}$ the associated semigroup.
For $d \geq 3$ and
$x \in \mathbb{R}^{d} \setminus 0$, we also denote
by $G(x)$ the Green function associated with $p$ :
$$ G(x) = \int_{0}^{\infty} p_{s}(x)ds = c_{d} |x|^{2-d}.$$
\PA The asymptotic behaviour of $q_{t}(y)$ as $t \to \infty$
is given by standard local limit theorems (see \cite{Sp}, and \cite{La}
for an equivalent statement for discrete random walks).
\begin{theorem} \label{LLT}
If $q$ and $p$ are defined as above,
$$\lim_{t \to \infty} \sup_{y \in \mathbb{Z}^{d}} \left| t^{d/2}q_{t}(y)
-p_{1}(y t^{-1/2}) \right| =0.$$
 \end{theorem}
\PA We will also need an upper bound on the transition kernel $q$ that is valid
for any $t \geq 1/2$ :
\begin{lemma} \label{boundq}
There exist two positive constants $\kappa_{1},\kappa_{2}$ such that
for every $t \geq 1/2$, $y \in \mathbb{Z}^{d}$,
$$ q_{t}(y) \leq \frac{\kappa_{1}}{t^{d/2}} \exp\left( -\frac{\kappa_{2}|y|}{\sqrt{t}}
\right). $$
\end{lemma}
For the reader's convenience, we provide a short proof of Lemma \ref{boundq}
in Section 5.
For $t>0$ and $y \in \mathbb{R}^{d}$ let
us define
$$f_{t}(y)  : = \frac{\kappa_{1}}{t^{d/2}} \exp\left( -\frac{\kappa_{2}|y|}{\sqrt{t}}
\right).$$
We also set for $t>0$ and $y \in \mathbb{R}^{d}$
\begin{eqnarray*}
&& \tilde{\kappa}_{1} := 2 \kappa_{1}; \ \tilde{\kappa}_{2} :=
\kappa_{2}/4; \qquad \tilde{f}_{t}(y) :=
\frac{\tilde{\kappa}_{1}}{t^{d/2}} \exp\left(
-\frac{\tilde{\kappa}_{2}|y|}{\sqrt{t}}
 \right); \\
&& \hat{\kappa}_{1} := 4 \kappa_{1}; \ \hat{\kappa}_{2} :=
\kappa_{2}/8; \qquad \hat{f}_{t}(y) :=
\frac{\hat{\kappa}_{1}}{t^{d/2}} \exp\left(
-\frac{\hat{\kappa}_{2}|y|}{\sqrt{t}}
 \right),
\end{eqnarray*}
so that $f_{t}(x) \leq \tilde{f}_{t}(x) \leq \hat{f}_{t}(x)$. We
need to control integrals of these functions. Note that, for $x \neq
0$, the supremum of the function $t \to \hat{f}_{t}(x)$ is reached
at $t_{0} = \frac{|x|^{2}\hat{\kappa}_{2}^{2}}{d^{2}}$. Let us
introduce for $r>0$
\begin{eqnarray*}
&& \psi_{2}(r)= 2 \ln(r \vee e), \\
&& \psi_{d}(r)=r^{d-2} \mbox{ if } d \geq 3.
\end{eqnarray*}
We then observe that for $T>0$, there
exists a constant $L_{0}$ depending only on $d$ and $T$ such that
for any $x \in \mathbb{R}^{d} \setminus 0$,
\begin{equation} \label{integf}
\int_{0}^{T} f_{t}(x)dt \leq \int_{0}^{T} \tilde{f}_{t}(x) dt \leq \int_{0}^{T}
\hat{f}_{t}(x)dt \leq L_{0} \psi_{d}\left( |x|^{-1}\right).
\end{equation}
Furthermore, whenever $|x| \geq \frac{d \sqrt{T}}{\hat{\kappa}_{2}}$
we have
\begin{equation} \label{integf2}
\int_{0}^{T} \hat{f}_{t}(x)dt \leq \hat{\kappa}_{1} T^{1-d/2} \exp(
-\frac{\hat{\kappa}_{2}|x|}{\sqrt{T}} ).
\end{equation}
Finally, when $d=2$, the integral $\int_{0}^{T}f_{t}(x) dt$ diverges
 when $T \to \infty$, but, when $d \geq 3$,
there exist a constant $L_{1}$ depending only on $d$ such that
\begin{equation} \label{integfoverR}
\int_{0}^{\infty} \hat{f}_{t}(x)dt \leq L_{1} \psi_{d}\left(
|x|^{-1}\right) = L_{1}|x|^{2-d}.
\end{equation}
\PA We also need an exponential bound on the probability for a random walk with jump kernel
$p$ to escape $B(0,A\sqrt{t})$ before time $t$. As a consequence of Lemma \ref{boundq} and Doob's maximal inequality applied to
a suitable exponential martingale of the random walk,
 there exist positive constants $\kappa_{3}, \kappa_{4}$ such that
for any $t \geq 1/2$, for any $A>0$,
\begin{equation} \label{exp}
\mathbb{P}_{0}(\sup_{s \in [0,t]} |Z_{s}| \geq A \sqrt{t}) \leq
\kappa_{3} \exp(-\kappa_{4}A).
\end{equation}
We may and will assume that the constant $\kappa_{2}$ in Lemma
\ref{boundq} is such that $\kappa_{4} \geq 4\kappa_{2}$.  \PA We
then deduce easy consequences of Theorem \ref{LLT} and Lemma
\ref{boundq}. From Theorem \ref{LLT},
 we obtain, for $x \neq 0$ and $s>0$,
\begin{equation} \label{LLT2}
c^{d} q_{c^{2}s}(c[x]_{c}) = s^{-d/2} (c^{2}s)^{d/2} q_{c^{2}s}(c[x]_{c})
\underset{c \to \infty}{\longrightarrow} s^{-d/2}
 p_{1}\left(\frac{x}{\sqrt{s}}\right) = p_{s}(x).
\end{equation}
On the other hand, using (\ref{exp}), we get
$$c^{d}\int_{0}^{c^{-2}} q_{c^{2}s}(c[x]_{c}) ds \leq \kappa_{3}c^{d-2}
\exp\left( -\kappa_{4}c[x]_{c} \right) \underset{c \to
\infty}{\longrightarrow} 0, $$ whereas, from Lemma \ref{boundq}, for
any $s \geq c^{-2}$ we have
\begin{equation} \label{boundqbyf}
c^{d}q_{c^{2}s}(c[x]_{c}) \leq f_{s}([x]_{c}).
\end{equation}
We can use (\ref{LLT2}) and dominated convergence to deduce that
 for $x \neq 0$ and $T>0$ we have
\begin{eqnarray} \label{EUT1}
c^{d} \int_{0}^{T} q_{c^{2}s}(c[x]_{c}) ds
\underset{c \to \infty}{\longrightarrow} \int_{0}^{T} p_{s}(x) ds.
\end{eqnarray}
By a similar argument, we obtain, for any
$T>0$, $y \in \mathbb{Z}^{d}$,
$$ \int_{0}^{T} q_{c^{2}s}(y) ds \leq c^{-2}\kappa_{3}\exp(-\kappa_{4}|y|) +
c^{-d} \int_{c^{-2}}^{T} f_{s}(y/c) ds.  $$
 Using (\ref{integf}),
(\ref{integf2}),
 it is then easy to establish that there exist constants $ L_{2} , L_{2}'$,
depending only on $T$ and $d$, such that for any $c \geq 1$, we have
\begin{eqnarray} \label{boundqc2s}
\int_{0}^{T} q_{c^{2}s}(y) ds \leq \begin{cases}
 L_{2} c^{-d} \psi_{d}(c) & \mbox{ if } y=0, \\
 L_{2} c^{-d} \psi_{d}\left( c|y|^{-1}\right) & \mbox{ if }
y \in \mathbb{Z}^{d} \setminus 0, \\
L_{2} c^{-d} \exp \left(- L_{2}' \frac{|y|}{c} \right) & \mbox{ if }
y \in \mathbb{Z}^{d}, |y| > c.
\end{cases}
\end{eqnarray}

We now discuss some preliminary results on coalescing random walks.
\subsection{Preliminary results on coalescing random walks}
 Consider two independent copies
$Z^{1},Z^{2}$ of the random walk $Z$ with transition kernel $q$,
starting respectively at points $y_{1},y_{2} \in \mathbb{Z}^{d}$
under the probability measure $\mathbb{P}_{y_{1},y_{2}}$. The time
at which $Z^{1}$ and $Z^{2}$ first meet is the stopping time $T_{1}
= \inf \{t \geq 0 : Z_{t}^{1}=Z_{t}^{2} \}$. We will need the
following result. The first bound below holds in the case $d \geq
3$, for which we recall that $\psi_{d}(r)=r^{d-2}$. The second bound
holds in the case $d=2$, for which we recall $\psi_{2}(r)=2\ln(r
\vee e)$.
\begin{lemma} \label{boundcoal}
Let $d\geq 2$ and $T>0$. There exists a positive constant $L_{4}$
depending only on $T$ and $d$ such that for any $x \in
\mathbb{R}^{d}\setminus 0$, for any $c \geq 1 \vee |x|^{-2}$ and for
any  $y \in \mathbb{Z}^{d} \setminus 0$,
\begin{eqnarray*}
\begin{cases}
\mbox{ if } d \geq 3, \quad
c^{d} \psi_{d}(|y|) \!\!\!&\!\!\!
 \int_{0}^{T} dt \mathbb{P}_{0,y} \left[ T_{1}
\leq c^{2}t, Z_{c^{2}t}^{1} = c [x]_{c} \right]  \leq
L_{4} \ \psi_{d}(|x|^{-1}),  \\
\mbox{ if } d = 2, \quad \ \ c^{2}
\frac{\psi_{2}(|y|)}{\psi_{2}(c/|y|)} \!\!\! & \!\!\! \int_{0}^{T}
dt \mathbb{P}_{0,y} \left[ T_{1} \leq c^{2}t, Z_{c^{2}t}^{1} = c
[x]_{c} \right]  \leq L_{4} \psi_{2}(|x|^{-1}). \end{cases}
\end{eqnarray*}
\end{lemma}
We postpone the proof of this result to Section 5.

\section{Upper bound}
In the case $d \geq 5$, the upper bound of Theorem \ref{main} follows
from the next proposition.
\begin{proposition} \label{uppd5}
Let $d \geq 5, x \in \mathbb{R}^{d}$. For $c$ large enough
$$ P(T_{c[x]_{c}} < \infty) \leq 2e c^{2-d} G(x). $$
\end{proposition}
 In the case $d \leq 3$, we will argue under $P_{c^{2}T}^{*}$ and use
 Theorem \ref{inv} to establish the following sharp asymptotic upper bound.
This bound also holds when $d \geq 4$ but is not sharp in that case.
\begin{proposition} \label{uppd3}
Let $d\geq 2$, $T>0, x \in \mathbb{R}^{d}\setminus \{ 0\}$,
 \begin{equation} \label{majsmall}
\limsup_{c \to \infty} P_{c^{2}T}^{*} (T_{c[x]_{c}}< \infty)
 \leq \mathbb{N}_{0}^{(T)} \ (x \in \mathcal{R}).
\end{equation}
\end{proposition}
In the cases $d=2$ or $d=3$, we will see in Section 3.3
that Proposition \ref{uppd3} implies the asymptotic upper bound in Theorem
\ref{main}. Notice that the right-hand side of (\ref{majsmall}) is $0$ if
$d \geq 4$.
 We begin with the proof of Proposition \ref{uppd5}, which only requires very
 simple arguments.
\subsection{The case $d \geq 5$}
Fix $x \in \mathbb{R}^{d}\setminus \{ 0\}.$
 Proving Proposition \ref{uppd5} reduces to establishing the following two
results :
\begin{equation} \label{upd51}
 E \left[ \int_{0}^{\infty} ds \mathbf{1}_{ \{ c[x]_{c} \in
\xi_{s}^{0} \} }\right] \underset{c \to \infty}{\sim} c^{2-d}
\int_{0}^{\infty} ds p_{s}(x),
\end{equation}
\begin{equation} \label{upd52}
P \left( T_{c[x]_{c}} < \infty \right) \leq e E \left[
\int_{0}^{\infty} ds \mathbf{1}_{ \{ c[x]_{c} \in \xi_{s}^{0} \}
}\right].
\end{equation}
Let us fix $T>0$, and observe that
\begin{eqnarray} \label{occtillT}
E \left[ \int_{0}^{c^{2}T} ds \mathbf{1}_{\{c[x]_{c} \in \xi_{s}^{0}\} }
\right]
\!\!&=&\!\! c^{2} E\left[\int_{0}^{T} ds
\mathbf{1}_{\{c[x]_{c} \in \xi_{c^{2}s}^{0}\} }
 \right] \nonumber \\
\!\!&=&\!\!
 c^{2} \int_{0}^{T} ds P \left( Z_{c^{2}s}^{c[x]_{c},c^{2}s}=0 \right) =
  c^{2} \int_{0}^{T} ds q_{c^{2}s}(c[x]_{c})  \nonumber \\
\!\!& \underset{c \to \infty}{\sim}&\!\! c^{2-d} \int_{0}^{T} p_{s}(x) ds,
\end{eqnarray}
where the asymptotics at the last line come from (\ref{EUT1}).
 Furthermore, using (\ref{boundqbyf}), we have similarly
\begin{eqnarray*}
E \left[ \int_{c^{2}T}^{\infty}
 \mathbf{1}_{\{c[x]_{c} \in \xi_{s}^{0}\} } ds \right]
&=& c^{2} \int_{T}^{\infty}   q_{c^{2}s}(c[x]_{c}) ds \\
& \leq & c^{2-d} \int_{T}^{\infty} f_{s}([x]_{c}) ds,
\end{eqnarray*}
and since $d \geq 3$, $\int_{T}^{\infty} f_{s}([x]_{c})ds $ goes to
$0$ as $ c \to \infty$. Thus from (\ref{occtillT}), we obtain
(\ref{upd51}). Let us now prove (\ref{upd52}). \PA
 When $T_{c[x]_{c}} < \infty$, denote by $N$ the numbers of arrows pointing
towards $c[x]_{c}$ in the time interval $(T_{c[x]_{c}}, T_{c[x]_{c}}+1]$.
Under $P\left(. \big| T_{c[x]_{c}} < \infty \right)$, $N$ is a Poisson variable
with parameter $1$. It follows  that
$$ P(T_{c[x]_{c}}< \infty) = e P(T_{c[x]_{c}}< \infty, N=0).$$
Furthermore, on the event $\{ N=0\}$ we have $c[x]_{c} \in \xi_{s}^{0}$ for
every \\ $s \in [T_{c[x]_{c}}, T_{c[x]_{c}}+1]$. Hence,
$$ E \left[
\int_{0}^{\infty} ds \mathbf{1}_{ \{ c[x]_{c} \in \xi_{s}^{0} \}
}\right] \geq P(T_{c[x]_{c}}< \infty, N=0). $$
This completes the proof of (\ref{upd52}), and of Proposition \ref{uppd5}. $\qquad \Box$

\subsection{Proof of Proposition \ref{uppd3}.}
Let $d \geq 2$ and fix $T>0$, $x \in \mathbb{R}^{d} \setminus 0$, and
$\eta \in (0,|x|/2)$.
Recall the notation $m_{N}$ from Section 2.2. We have for any $\delta >0$,
 $\varepsilon >0$ :
\begin{eqnarray} \label{bddocctime}
P_{c^{2}T}^{*}(T_{c[x]_{c}} < \infty)
 \!\!&\!\! \leq \!\!&\!\! P_{c^{2}T}^{*}
\left[  \int_{0}^{\infty} ds
 \mathbf{1}_{ \{ X_{s}^{0}(\overline{B}(c x ,\eta c))
\geq \delta m_{c^{2}}\} } \geq  \varepsilon c^{2} \right] \nonumber \\
 && \qquad + P_{c^{2}T}^{*}
   \left[  \int_{0}^{\infty} ds
 \mathbf{1}_{ \{ X_{s}^{0}(\overline{B}(c x,\eta c))
\geq \delta m_{c^{2}}\} } <  \varepsilon c^{2},
T_{c[x]_{c}} <\infty \right].
\end{eqnarray}
Intuitively, when $c$ tends to infinity,
the second term of the sum above should
remain small when $\varepsilon$ and $\delta$ are small enough, while the
first term, using the invariance principle, should be bounded by a
corresponding rescaled
quantity under $\mathbb{N}_{0}^{(T)}$. Let us be more precise.
Using rescaling, the first term of the sum
in the right-hand side of (\ref{bddocctime}) is equal to
 $$ P_{c^{2}T}^{*}
\left[  \int_{0}^{\infty} ds
 \mathbf{1}_{ \{ X_{s}^{c^{2},0}(\overline{B}( x ,\eta ))
\geq \delta \} } \geq  \varepsilon  \right].$$
It is easy to see that for any $A>0$, the set
$$\left\{ \omega \in D(\mathbb{R}_{+}, M_{F}(\mathbb{R}^{d})) :
\int_{0}^{A}  ds
 \mathbf{1}_{ \{ \omega_{s}(\overline{B}( x ,\eta ))
\geq \delta \} } \geq  \varepsilon  \right\}  $$
is closed for the Skorohod $J_{1}$ topology. Then, Theorem \ref{inv}
 implies that
\begin{eqnarray*}
  \limsup_{c \to \infty} P_{c^{2}T}^{*}
\left[  \int_{0}^{A} ds
 \mathbf{1}_{ \{ X_{s}^{c^{2},0}(\overline{B}( x ,\eta ))
\geq \delta \} } \geq  \varepsilon  \right]
 \!\!\!&\!\!\!&\!\!\! \leq \mathbb{N}_{0}^{(T)} \left[ \int_{0}^{A} ds
 \mathbf{1}_{ \{ Y_{s}(\overline{B}( x ,\eta ))
\geq \delta \} } \geq  \varepsilon  \right] \\ \!\!\!&\!\!\!&\!\!\!
 \leq \mathbb{N}_{0}^{(T)}
\left[ Y \mbox{ hits } \overline{B}(x,\eta) \right].
\end{eqnarray*}
Furthermore, we have, for $A \geq T$,
\begin{eqnarray*}
 P_{c^{2}T}^{*}
\left[  \int_{A}^{\infty} ds
 \mathbf{1}_{ \{ X_{s}^{c^{2},0}(\overline{B}( x ,\eta ))
\geq \delta \} } \geq  \varepsilon \right] \leq
P_{c^{2}T}^{*}(X_{A}^{c^{2},0} \neq 0) =\frac{p_{c^{2}A}}{p_{c^{2}T}},
\end{eqnarray*}
which goes to $0$ as $A \to \infty$. Hence, we obtain for every
$\delta>0, \epsilon >0$,
\begin{equation} \label{appinv}
  \limsup_{c \to \infty} P_{c^{2}T}^{*}
\left[  \int_{0}^{\infty} ds
 \mathbf{1}_{ \{ X_{s}^{c^{2},0}(\overline{B}( x ,\eta ))
\geq \delta \} } \geq  \varepsilon  \right]
 \leq \mathbb{N}_{0}^{(T)}
\left[ Y \mbox{ hits } \overline{B}(x,\eta) \right]. \end{equation}
\PA To control the second term of the sum in the right-hand side of
(\ref{bddocctime}), we will use the following argument. When $c$ is large
and point
$c[x]_{c}$ is hit by opinion $1$, then with arbitrarily
high probability, a sufficient
 number (of order $m_{c^{2}}$) of its
neighbors (at distance less than $\eta c$) should also be hit by
opinion $1$ during a certain time interval (with length of order $c^{2}$). \PA
We will prove a
somewhat more general result, which will be valid uniformly over all points in
$\xi_{t}^{0}$, with the restriction that $t$ should be at least of order
$c^{2}$.

\begin{lemma} \label{enough}
Let $T>0, \rho > 0$, $\eta>0$ be fixed. We can find $\varepsilon_{0}> 0$
so that for any $\varepsilon \in (0, \varepsilon_{0}]$, there exists
$\delta >0$ such that for $c$ sufficiently large,
\begin{equation} \label{enoughpoints}
P_{c^{2}T}^{*}\left( \exists t \geq 4 \varepsilon c^{2} \
\exists x \in \xi_{t}^{0}
 : \inf_{s \in [t-3\varepsilon c^{2},t-2\varepsilon c^{2}]} |\xi_{s}^{0}\cap
\overline{B}(x,\eta c) |< \delta m_{c^{2}} \right) \leq \rho.
\end{equation}
\end{lemma}
We will also need a useful exponential bound on the probability
for the voter model to escape a ball of radius $A\sqrt{\alpha}$ before time
$2\alpha$ :
\begin{lemma} \label{speed}
 There exists constants $K_{1}>0, K_{2} >0$ such that for any
$\alpha > 1$, for any $A>0$,
\begin{equation} \label{deviation}
P_{\alpha}^{*} \left(\sup_{t \leq 2 \alpha} \sup_{x\in \xi_{t}^{0}} |x| >
A \sqrt{\alpha} \right) \leq K_{1}\exp(-K_{2} A).
\end{equation}
\end{lemma}
Let us postpone the proofs of Lemma \ref{enough} and Lemma \ref{speed},
and finish the proof of Proposition
2. Recall $x,T,\eta \in  (0,|x|/2) $ have been fixed.
Notice that, when $c$ is large enough,
$\overline{B}(c[x]_{c},\eta c/2) \subset \overline{B}(cx,\eta c)$. Thus,
\begin{eqnarray*}
&& P_{c^{2}T}^{*}
   \left[  \int_{0}^{\infty} ds
 \mathbf{1}_{ \{ X_{s}^{0}(\overline{B}(c x,\eta c))
\geq \delta m_{c^{2}}\} } <  \varepsilon c^{2}, 4\varepsilon c^{2} \leq
T_{c[x]_{c}} <\infty \right] \nonumber \\
&& \leq P_{c^{2}T}^{*}
   \left[  \int_{0}^{\infty} ds
 \mathbf{1}_{ \{ X_{s}^{0}(\overline{B}(c [x]_{c},\eta c/2))
\geq \delta m_{c^{2}}\} } <  \varepsilon c^{2}, 4\varepsilon c^{2} \leq
T_{c[x]_{c}} <\infty \right].
\end{eqnarray*}
Hence, using Lemma \ref{enough}, for any $\rho>0$,
we can choose $\varepsilon_{0}>0$
such that for any $\varepsilon \in (0,\varepsilon_{0}]$,
there exists $\delta >0$ such that for $c$ large enough,
\begin{eqnarray} \label{3epsc2}
 P_{c^{2}T}^{*}
   \left[  \int_{0}^{\infty} ds
 \mathbf{1}_{ \{ X_{s}^{0}\left(\overline{B}(c [x]_{c},\eta c/2)\right)
\geq \delta m_{c^{2}}\} } <  \varepsilon c^{2}, 4\varepsilon c^{2} \leq
T_{c[x]_{c}} <\infty \right] \leq \rho.
\end{eqnarray}
 Furthermore, provided
$2\varepsilon \leq T$, we have
$$ P_{c^{2}T}^{*}
   \left[ T_{c[x]_{c}} < 4 \varepsilon c^{2} \right] \leq
\frac{p_{2\varepsilon c^{2}}}{p_{c^{2}T}}
P_{2 \varepsilon c^{2}}^{*}
\left[ T_{c[x]_{c}} < 4 \varepsilon c^{2} \right].$$
If $c$ is sufficiently large, we can thus use (\ref{asy})
and the fact that $c|[x]_{c}| \geq c|x|/\sqrt{2}$,
then Lemma \ref{speed} with $\alpha = 2 \varepsilon c^{2}$ and
$A = \frac{|x|}{2\sqrt{\varepsilon}}$ to get
\begin{eqnarray*}
P_{c^{2}T}^{*}
   \left[ T_{c[x]_{c}} < 4 \varepsilon c^{2} \right] &\leq &
\kappa_{0} \frac{T}{2\varepsilon} P_{2\varepsilon c^{2}}^{*}\left(
\sup_{t \leq 4 \varepsilon c^{2}}
\sup_{y \in \xi_{t}^{0}} |y| > \frac{c|x|}{\sqrt{2}} \right)  \\
&\leq & \kappa_{0}K_{1} \frac{T}{2\varepsilon} \exp \left(
-\frac{K_{2}|x|}{2 \sqrt{\varepsilon}} \right).
\end{eqnarray*}
Combining (\ref{bddocctime}), (\ref{appinv}), (\ref{3epsc2}) and the
last inequality now yields
\begin{equation*}
\limsup_{c \to \infty} P_{c^{2}T}^{*}(T_{c[x]_{c}} < \infty) \leq
\mathbb{N}_{0}^{(T)} \left[ Y \mbox{ hits } \overline{B}(x,\eta)
\right] + \rho + \kappa_{0}K_{1} \frac{T}{2\varepsilon} \exp \left(
-\frac{K_{2}|x|}{2\sqrt{\varepsilon}} \right),
\end{equation*}
for any $\rho>0$ and
$\varepsilon \in (0,\varepsilon_{0}(\eta,\rho)]$.
By letting $\varepsilon$ and then $\rho$ go to $0$, we get
\begin{equation} \label{bornsu}
\limsup_{c \to \infty}
P_{c^{2}T}^{*}(T_{c[x]_{c}} < \infty) \leq \mathbb{N}_{0}^{(T)}
\left[ Y \mbox{ hits } \overline{B}(x,\eta) \right].
\end{equation}
Our reasonning is valid for any $\eta \in (0,|x|/2)$.
Thus, letting $\eta$ go to $0$ in (\ref{bornsu})
finishes the proof of Proposition \ref{uppd3}.
$ \qquad \Box$ \PA
 It remains to prove Lemma \ref{speed} and Lemma \ref{enough}. We start with the proof
of Lemma \ref{speed}, since it will appear to
be a key tool in the proof of Lemma \ref{enough}.

\subsubsection{Proof of Lemma \ref{speed}}
  Let us first outline the proof and summarize the intermediate results. We need to discretize the time scale.
Introduce the integer
$$N := \min \{ n \in \mathbb{N} : \alpha 2^{-n} <1 \} =
\left\lfloor \frac{\ln(\alpha)}{\ln(2)} \right\rfloor,$$
and the time intervals
$$ B_{n} := [(n-1)2^{-N-1}\alpha, n 2^{-N-1}\alpha], \ \ n \in \big[\!| 1,2^{N+2}|\!\big].$$
Let us introduce the set of points having, for some odd $n \in \big[\!| 1,2^{N+2}|\!\big],$ opinion $1$ at a time
belonging to $B_{n}$, and descendants at time $(n+1)2^{-N-1}\alpha$ :
 $$ \Xi_{N+1} := \bigcup_{\stackrel{n=1}{n \ \mathrm{ odd }}}^{2^{N+1}} \bigcup_{u \in B_{n}}
 \Omega_{u}^{(n+1)2^{-N-1}\alpha}.$$
 Informally,
our interest in this set $\Xi_{N+1}$ comes from the fact
that if $x \in \bigcup_{t \leq 2\alpha} \xi_{t}^{0}$, a ``close'' ancestor of $x$ belongs to $\Xi_{N+1}$, and hence,
$\Xi_{N+1}$ should not be too far from $\bigcup_{t\leq 2\alpha} \xi_{t}^{0}$.
More precisely, for $t<1$, set $u_{t}=0$, and
for $t \in [1, 2\alpha]$, let us choose
$$u_{t} \in \left[t- \frac{\alpha}{2^{N}}, t-\frac{\alpha}{2^{N+1}}\right] \bigcap
\bigcup_{\stackrel{n=1}{n \ \mathrm{ odd}}}^{2^{N+1}} B_{n}.$$
We have $t-u_{t} \leq \alpha 2^{-N} <1$, and,
if $x \in \xi_{t}^{0}$ for some $t \in [0,2\alpha]$, the ancestor of $x$ at time $u_{t}$
indeed belongs to $\Xi_{N+1}$. \PA
 We will show that, under $P_{\alpha^{*}}$,
 $\Xi_{N+1}$ intersects $B\left(0,\frac{A}{2}\sqrt{\alpha}\right)$
 with a probability which decays exponentially with $A$.
\begin{lemma} \label{levels}
There exist positive constants $K_{3},K_{4}$ such that for any $A>0$,
for any $\alpha >1$,
\begin{eqnarray} \label{level}
 P_{\alpha}^{*} \left(\Xi_{N+1} \nsubseteq B(0,\frac{A}{2}
\sqrt{\alpha}) \right) \leq K_{3} \exp(-K_{4}A).
\end{eqnarray}
\end{lemma}
Then, we will argue that the probability under $P_{\alpha}^{*}$ for
$\bigcup_{t \leq 2 \alpha} \xi_{t}^{0}$ to escape the ball $B(0,A
\sqrt{\alpha})$ and simultaneously to have
 $\Xi_{N+1} \subset B(0,\frac{A}{2}\sqrt{\alpha})$
also decays exponentially with $A$.
This is seen below as a consequence of the following result.
\begin{lemma} \label{quickescape}
There exist positive constants $K_{5},K_{6}$ such that for any $A>0$,
$$ P\left( \exists t \in [0,1] \ \exists x \in \mathbb{Z}^{d}\setminus
B(0,A) \ \exists y \in B(0,\frac{A}{2})  :  (t,x) \searrow
 (0,y) \right) \leq K_{5} \exp(-K_{6}A). $$
\end{lemma}
Let us postpone the proofs of Lemmas \ref{levels} and \ref{quickescape} and show how Lemma \ref{speed} is deduced from these two results.
Introduce the event
$$\mathcal{A} := \bigg\{\exists
t \leq 2 \alpha  \ \exists x \in \xi_{t}^{0} :
 |x| > A \sqrt{\alpha}\bigg\}   \bigcap  \bigg\{
\Xi_{N+1} \subset B(0,\frac{A}{2} \sqrt{\alpha})  \bigg\}.$$
Clearly, $\mathcal{A} = \bigcup_{n=1}^{2^{N+2}} \mathcal{A}_{n}$, where
$$ \mathcal{A}_{n} := \bigg\{ \exists
t \in B_{n} \ \exists x \in \xi_{t}^{0} \ : \
 |x| > A \sqrt{\alpha} \bigg\}  \bigcap \bigg\{\Xi_{N+1} \subset
B(0,\frac{A}{2} \sqrt{\alpha}) \bigg\}.$$
As we noticed earlier, when $x \in \xi_{t}^{0}$, the ancestor of $x$ at time $u_{t}$ belongs to $\Xi_{N+1}$.
Hence, using the Markov property at time $u_{t}$, we get, for
every $n \in \big[\!|1,2^{N+2} |\!\big]$,
$$ P(\mathcal{A}_{n}) \leq P\!\!\left(
\exists s \!\in \!\left[0,\frac{\alpha}{2^{N}}\right] \
\exists x \!\in \!\mathbb{Z}^{d}\!
\setminus \!
B(0,A\sqrt{\alpha}) \ \exists y \!\in\! B(0,\frac{A}{2}\sqrt{\alpha})
\ : \ (s,x) \searrow
 (0,y) \right), $$
where we used that $t-u_{t} \leq \alpha 2^{-N}$. Using the fact that
$\alpha 2^{-N} <1$, it then follows from Lemma \ref{quickescape}
that
$$ P \left(\mathcal{A} \right) \leq 2^{N+2} K_{5} \exp(-K_{6}A \sqrt{\alpha}).
$$
Since $2^{N+2} \leq 8\alpha$ from the definition of $N$, it follows
from the above that \\$P_{\alpha}^{*} \left(\mathcal{A} \right) \leq
8 K_{5} \alpha^{2} \exp(-K_{6}A \sqrt{\alpha})$. Hence, there exists
positive constants $K_{1}', K_{2}'$ such that
$P_{\alpha}^{*}(\mathcal{A}) \leq K_{1}' \exp(-K_{2}' A).$ This fact
and Lemma \ref{levels} imply Lemma \ref{speed}. $\quad \Box$
\vskip0.2cm \PA
It now remains to prove Lemmas \ref{levels} and \ref{quickescape}. We first establish Lemma \ref{quickescape}. \\
\underline{Proof of Lemma \ref{quickescape}} :
   Fix $x \in \mathbb{Z}^{d} \setminus B(0,A)$. There is a Poisson number $n_{x}$
with parameter $1$ of arrows
pointing towards $x$ during the time interval $[0,1]$. Denote by
$1 \geq T_{1} > T_{2}>  ... > T_{n_{x}}\geq 0$ the times at which these
arrows occur and by $z_{1},z_{2},...,z_{n_{x}}$
the respective origins of these arrows. We also set $T_{i}=0$ when $i>n_{x}$.
For $t \in [0,1]$ and
$y \in B(0,A/2)$, a path
up $(0,y) \nearrow (x,t)$ has to ``follow''
one of the $n_{x}$ arrows pointing towards $x$ in the time interval $[0,1]$,
say the $i$th one at time $T_{i}$, in this case we then have
$(z_{i},T_{i}) \searrow (y,0)$.
 \PA For $t \in [0,1]$, let us define
$\mathcal{G}_{t}$ the $\sigma$-field which is generated by the random sets
 $\left( \Lambda(x,y)\cap[1-t,1]\right)$ for all $x,y \in \mathbb{Z}^{d}$.  \PA
 The times
$1-T_{i}, i \in \mathbb{N}$ are stopping times for the filtration
$(\mathcal{G}_{t})_{t \in [0,1]}$, and conditionally on $\{n_{x}=k
\}$, the points $z_{i},1 \leq i \leq k$ are located independently
according to $p(x,.)$. In particular, using the exponential moments
assumption on $p$, there exist positive $\tilde{\kappa}_{3},
\tilde{\kappa}_{4}$ such that for any $1 \leq i \leq k$,
\begin{equation} \label{loczi}
P\left(z_{i} \in B\left(0,\frac{3|x|}{4}\right) \bigg| n_{x} =k
\right) \leq \tilde{\kappa}_{3}  \exp(-\tilde{\kappa}_{4}|x|).
\end{equation}
For $i \in \big[\!|1,k|\!\big]$, let us define
$(\tilde{Z}_{s}^{x,T_{i}})_{0 \leq s \leq T_{i}}$ as follows
\begin{itemize}
\item for $0 < s \leq T_{i}$, $\tilde{Z}_{s}^{x,T_{i}}=Z_{s}^{x, T_{i}}$
\item $\tilde{Z}_{0}^{x,T_{i}} = z_{i}$.
\end{itemize}
For $t>0$, conditionally on $\{T_{i}=t\}$,
$(\tilde{Z}_{s}^{x,T_{i}})_{0 \leq s \leq t}$ is a rate $1$
random walk with jump kernel $p$ started from $z_{i}$, and is
thus distributed as $(Z_{s})_{0\leq s \leq t}$ under $\mathbb{P}_{z_{i}}$.
Furthermore, using (\ref{exp}), we have for any
$z \in \mathbb{Z}^{d} \setminus B\left(0,\frac{3|x|}{4} \right)$
\begin{equation} \label{pzi}
\mathbb{P}_{z} \left(\exists s \in [0,1] : Z_{s} \in
B\left(0,\frac{|x|}{2} \right) \right) \leq \kappa_{3}
\exp\left(-\kappa_{4}\frac{|x|}{4}\right).
\end{equation}
Combining (\ref{loczi}) and (\ref{pzi}), we see that
there exist positive constants $K_{5}', K_{6}'$ such that for any
$x \in \mathbb{Z}^{d} \setminus B(0,A)$, for any $k \in \mathbb{N}$
\begin{equation*}
P\left(\exists i \in \big[\!| 1,k |\!\big] :  Z_{T_{i}}^{x,T_{i}}
\in B\left(0,\frac{A}{2} \right) \ \bigg| \ n_{x}=k \right) \leq
K_{5}' k \exp \left(-K_{6}'|x|\right).
\end{equation*}
We thus get
\begin{eqnarray*}
&& P\left( \exists t \in [0,1] \ \exists x \in \mathbb{Z}^{d}\setminus
B(0,A) \ \exists y \in B\left(0,\frac{A}{2}\right) \ : \ (t,x) \searrow
 (0,y) \right) \\ &&
\leq \sum_{x \in \mathbb{Z}^{d} \setminus B(0,A)}
\sum_{k=0}^{\infty} P\left(n_{x}=k\ , \  \exists i \in
\big[\!|1,k|\!\big] \ : \ Z_{T_{i}}^{x,T_{i}} \in B(0,\frac{A}{2})
\right)
\\ && \leq \sum_{x \in \mathbb{Z}^{d} \setminus B(0,A)}
\frac{1}{e} \sum_{k=0}^{\infty} \frac{1}{(k-1)!} K_{5}'
\exp\left( -K_{6}' |x| \right).
\end{eqnarray*}
Lemma \ref{quickescape} follows. $\qquad \Box$ \vskip0.2cm
To prove Lemma \ref{levels}, we need the following key result.
 \begin{lemma} \label{AA'}
 Let $t \geq 0,s>r>0$ and $A > A' \geq 0$.
  \begin{eqnarray*}
 P\left( \!\bigcup_{u \in [0,t]} \!\!\Omega_{u}^{t+s} \subset B(0,A'),
 \!\!\!\bigcup_{u \in [t,t+r]} \!\!\!\!
 \Omega_{u}^{t+s} \nsubseteq B(0,A) \!\right) \! \leq
 p_{s-r} \kappa_{3} \exp \left( -\kappa_{4} \frac{A-A'}{\sqrt{r}}\right).
 \end{eqnarray*}
 \end{lemma}
  \underline{Proof of Lemma \ref{AA'}}:
The event $\left\{ \Omega_{t}^{t+s} \subset B(0,A'),
 \bigcup_{u \in [t,t+r]} \!
 \Omega_{u}^{t+s} \nsubseteq B(0,A)  \right\}$
 considered in Lemma \ref{AA'} is contained in the event that there
exists a point $z \in \xi_{t+r}^{0}$ having descendants at time $t+s$,
such that the ancestor of $z$ at time $t$ belongs to $B(0,A')$, and moreover,
$z$ has an ancestor in $B(0,A)^{c}$ at a time belonging to $[t,t+r]$.
More precisely, using duality over the time
interval $[0,t+r]$, and then decomposing over all possible values of the point $z$,
\begin{eqnarray} \label{fir}
&& P\left( \Omega_{t}^{t+s} \subset B(0,A'),
\bigcup_{u \in [t,t+r]} \Omega_{u}^{t+s} \nsubseteq B(0,A) \right)\nonumber \\
&& \leq  P \left( \exists z \in \xi_{t+r}^{0} : \hat{\xi}_{s-r}^{z,t+r} \neq
\emptyset, \sup_{u \in [0,r]} |Z_{u}^{z,t+r}| > A,
|Z_{r}^{z,t+r}| \leq A', Z_{t+r}^{z,t+r}=0 \right) \nonumber \\
&& \leq \sum_{z \in \mathbb{Z}^{d}}
 P \left( \hat{\xi}_{s-r}^{z,t+r} \neq
\emptyset, \sup_{u \in [0,r]} |Z_{u}^{z,t+r}| > A,
|Z_{r}^{z,t+r}| \leq A', Z_{t+r}^{z,t+r}=0 \right).
\end{eqnarray}
Using the Markov property at
time $(t+r)$, we obtain that the quantity in the right-hand side of (\ref{fir})
is equal to
\begin{eqnarray} \label{sec}
&& p_{s-r} \sum_{z \in \mathbb{Z}^{d}} \mathbb{P}_{z}\left(
\sup_{u \in [0,r]} |Z_{u}| > A,
|Z_{r}| \leq A', Z_{t+r}=0 \right) \nonumber \\
&& = p_{s-r} \sum_{z \in \mathbb{Z}^{d}} \mathbb{P}_{0}\left(
 |Z_{r}| \leq A', \sup_{u \in [t,t+r]} |Z_{u}| > A,
Z_{t+r}=z \right) \nonumber \\
  && \leq p_{s-r}\mathbb{P}_{0}\left(
 |Z_{r}| \leq A', \sup_{u \in [t,t+r]} |Z_{u}| > A \right),
\end{eqnarray}
where, at the second line above, we used a time-reversal argument together with the
symmetry assumption we made on the jump kernel $p$. From (\ref{fir}),
(\ref{sec}) and the Markov property for the random walk at time
 $t$, we now obtain
\begin{eqnarray*}
&& P\left( \Omega_{t}^{t+s} \subset B(0,A'),
\bigcup_{u \in [t,t+r]} \Omega_{u}^{t+s} \nsubseteq B(0,A) \right) \\
&& \leq p_{s-r}\mathbb{P}_{0} (\sup_{u\in [0,r]} |Z_{u}| > A-A' ),
\end{eqnarray*}
and we conclude using (\ref{exp}). $ \qquad \Box $ \vskip0.2cm
\underline{Proof of Lemma \ref{levels}}:
Let us first note that we only need to establish the existence
of positive $K_{3}', K_{4}'$ such that (\ref{level}) holds for any
$A \geq 1$ and $\alpha >1$. Indeed, Lemma \ref{levels} will follow from
 taking  $K_{3}=K_{3}' \vee \exp(K_{4}'), K_{4}:=K_{4}'$.
 For $p \in \big[\!| 0, N+1|\!\big]$, let us introduce the sets
$$ \Xi_{p} := \bigcup_{\stackrel{n=1}{n \ \mathrm{ odd}}}^{2^{p+1}} \ \ \bigcup_{u \in [(n-1) 2^{-p}\alpha, n2^{-p}
\alpha ]} \Omega_{u}^{(n+1)2^{-p}\alpha} .$$ For convenience, we
also set $\Xi_{-1} := \emptyset$. In the case $p=N+1$, this is of
course consistent with our definition of $\Xi_{N+1}$. For $p \geq
0$, let $A_{p} := \frac{A}{14} \sum_{i=0}^{p} 2^{-i/4}$ and set
$A_{-1} =0$. so that for any $k \in \mathbb{N}$, $A_{k} \leq
\frac{A}{2}$.
 For  $p \in \big[\!| 0,N+1 |\!\big] $, let
$$ \mathcal{E}_{p} : = \left\{ \Xi_{p-1} \subset B(0,A_{p-1}\sqrt{\alpha}),
 \Xi_{p} \nsubseteq B(0,A_{p} \sqrt{\alpha})\right\}.$$
Note that $\mathcal{E}_{p}$ is a subset of
\begin{eqnarray*}
&& \mathcal{F}_{p} := \bigg\{ \exists n \in \big[\! | 1, 2^{p+1}|\!\big], n \mbox{ odd } :
 \Omega_{(n-1)2^{-p}\alpha}^{(n+1)2^{-p}\alpha}
\subset B(0, A_{p-1}\sqrt{\alpha}), \\
 && \qquad \qquad \qquad \qquad  \bigcup_{u \in [(n-1)2^{-p}\alpha,
 n2^{-p}\alpha ]}
\Omega_{u}^{(n+1)2^{-p}\alpha} \nsubseteq B(0, A_{p}\sqrt{\alpha})\bigg\}.
\end{eqnarray*}
 Hence,
\begin{equation} \label{xipxip-1}
P_{\alpha}^{*} \left(\Xi_{N+1} \nsubseteq B(0,\frac{A}{2}
\sqrt{\alpha}) \right) \leq \sum_{p=0}^{N+1}
P_{\alpha}^{*} \left(\mathcal{E}_{p} \right) \leq \sum_{p=0}^{N+1}
P_{\alpha}^{*}\left(\mathcal{F}_{p}\right)
\end{equation}
From Lemma \ref{AA'}, we obtain
$$ P \left(\mathcal{F}_{p} \right) \leq \sum_{\stackrel{n=1}{n \ \mathrm{odd}}}^{2^{p+1}}
 p_{2^{-p}\alpha} \exp \left( -\kappa_{4} 2^{p/2} (A_{p}-A_{p-1})
\right).
$$
Hence, using our definition of the numbers $A_{p}, p \geq -1$, then
(\ref{asy}), we get

$$P_{\alpha}^{*}(\mathcal{F}_{p}) \leq p_{\alpha}^{-1} p_{2^{-p}\alpha} \exp \left( -\kappa_{4} 2^{p/4} \frac{A}{14}\right)
\leq \kappa_{0} 2^{2p} \exp \left( -\kappa_{4} 2^{p/4} \frac{A}{14}
\right). $$
 From (\ref{xipxip-1}) and the last inequality,
elementary arguments then give Lemma \ref{levels}. $\qquad \Box$
\vskip0.1cm This completes the proof of Lemma \ref{speed}. To finish
the one of Proposition \ref{uppd3}, it remains to establish Lemma
\ref{enough}.

 \subsubsection{Proof of Lemma \ref{enough}}
 Fix $T>0, \rho >0, \eta>0$. Recall $K_{1}, K_{2}$ are the constants appearing
in the statement of Lemma \ref{speed}. We can choose $\varepsilon_{0} \in (0,1)$ so that
for any $ \varepsilon \in (0, \varepsilon_{0}]$,
$$ K_{1} \exp\left( -K_{2} \frac{\eta}{
4 \sqrt{\varepsilon}} \right)  \leq \frac{\rho^{2}\varepsilon^{2}}{8 T^{2}}.$$
Let us now fix $ \varepsilon \in (0,\varepsilon_{0}]$.
We can then choose $\delta>0$ small enough so that
$$  \mathbb{N}_{0}^{(4\varepsilon)} \left[
\inf_{t \in [\varepsilon,3\varepsilon]} Y_{t}(1) \leq
\delta \right] \leq \frac{ \rho^{2} \varepsilon^{2}}{8 T^{2}}.$$
 The reasons for our choices of $\epsilon_{0}$ and $\delta$ will become clear
in the following. \PA
We first need to reduce the problem to a finite
time interval. Notice that
 $P_{c^{2}T}^{*}(\xi_{c^{2}\frac{4T}{\rho}} \neq \emptyset)
 \leq p_{c^{2}T}^{-1} p_{c^{2}\frac{4T}{\rho}}$ which, using (\ref{pt}), is
bounded by $\rho/2$ for $c$ large enough.
Thus, to establish Lemma
\ref{enough} we only need to prove that provided $c$ is sufficiently large,
\begin{equation} \label{intfini}
P_{c^{2}T}^{*}\left( \exists t \in
[4 \varepsilon c^{2}, \frac{4T}{\rho} c^{2}] \
\exists x \in \xi_{t}^{0}
 : \inf_{s \in [t-3\varepsilon c^{2},t-2\varepsilon c^{2}]} |\xi_{s}^{0}\cap
\overline{B}(x,\eta c) |< \delta m_{c^{2}} \right) \leq \frac{\rho}{2},
\end{equation}
 \PA
Set $M := \left\lceil \frac{4T}{\rho \varepsilon} \right\rceil$. Let
us discretize the time scale via introducing the levels $
\mathcal{L}_{k}:= k\varepsilon c^{2}, k\in \big[\!|0,M-4|\!\big]$.
\PA We are going to establish, using Lemma \ref{speed}, that
 with arbitrarily high probability, when $c$ is large enough,
each point holding opinion $1$ at such a level $\mathcal{L}_{k}$
 and having descendants at time $\mathcal{L}_{k}+ 4\varepsilon c^{2}$ is close
(at a distance less than $\eta c/2$) to all its descendants during
the time interval $[\mathcal{L}_{k}, \mathcal{L}_{k}+5\varepsilon
c^{2}]$. Then, using Theorem 2, we will prove that such a point has
more than
 $\delta m_{c^{2}}$ descendants
in the time interval $[\mathcal{L}_{k}+\varepsilon c^{2},
\mathcal{L}_{k}+3\varepsilon c^{2}]$.
 \PA Let us be more precise.
We shall prove that if $c$ is large enough,
\begin{eqnarray} \label{descenproches}
&& P_{c^{2}T}^{*} \bigg( \exists k \in \big[\!|0,M-4|\!\big] \ \ \exists y \in
\Omega_{k \varepsilon c^{2}}^{(k+4)\varepsilon c^{2}}  :
  \bigcup_{s \in [0, 5\varepsilon c^{2}]}
\hat{\xi}_{s}^{y,k \varepsilon c^{2}} \not\subset
\overline{B}(y,\eta c /2) \bigg)
\leq \frac{\rho}{4},
\end{eqnarray}
\begin{eqnarray} \label{assezdescen}
&& P_{c^{2}T}^{*} \bigg( \exists k \in \big[\!|0,M-4|\!\big] \ \ \exists y \in
\Omega_{k \varepsilon c^{2}}^{(k+4)\varepsilon c^{2}}  :
  \inf_{s \in [\varepsilon c^{2}, 3\varepsilon c^{2}]}
|\hat{\xi}_{s}^{y,k \varepsilon c^{2}}| < \delta m_{c^{2}} \bigg)
\leq \frac{\rho}{4}.
\end{eqnarray}
Let us postpone the proof of these two results and show how (\ref{intfini})
 follows from (\ref{descenproches}) and (\ref{assezdescen}).
Consider $ t \in [4 \varepsilon c^{2}, M\varepsilon c^{2}]$ and
$x \in \xi_{t}^{0}$. Introduce
$$ k_{t} := \sup\{ n \in \mathbb{N} :
n\varepsilon c^{2} \leq t - 4 \varepsilon c^{2}\}; \ \
z_{x}:= Z_{t-k_{t}\varepsilon c^{2}}^{x,t},$$
so that
$z_{x} \in \Omega_{k_{t}\varepsilon c^{2}}^{t} \subset
\Omega_{k_{t}\varepsilon c^{2}}^{(k_{t}+4)\varepsilon c^{2}}$ (indeed
$z_{x}$
is the ancestor of $x \in \xi_{t}^{0}$
at time $k_{t}\varepsilon c^{2}$). Thus,
\begin{itemize}
\item using (\ref{descenproches}), with probability at least
$1 - \rho/4$, for any
$x \in \xi_{t}^{0}$, $t \in [4\varepsilon c^{2},M\varepsilon c^{2}]$,
 all descendants of $z_{x}$
until time $(k_{t}+5)\varepsilon c^{2}$ belong to $\overline{B}(x,\eta c)$.
\item using (\ref{assezdescen}), with probability at least $1-\rho/4$, for any
$x \in \xi_{t}^{0}$, $t \in [4\varepsilon c^{2},M\varepsilon c^{2}] $,
 there are more than
 $\delta m_{c}^{2}$ descendants of $z_{x}$ at every time \\
$s \in [(k_{t}+1)\varepsilon c^{2}, (k_{t}+3) \varepsilon c^{2}]$.
\end{itemize}
Since $[t-3\varepsilon c^{2},t-2 \varepsilon c^{2}]
\subset [(k_{t}+1)\varepsilon c^{2}, (k_{t}+3) \varepsilon c^{2}]$,
we now deduce from the above that with
probability at least $1-\rho/2$, for any $x \in \xi_{t}^{0}$,
$z_{x}$ has at least $\delta m_{c}^{2}$ descendants
in $\overline{B}(x, \eta c)$ at every time
$s \in [t-3\varepsilon c^{2},t-2 \varepsilon c^{2}]$.
 Assertion (\ref{intfini}) follows. \PA
 Let us now prove (\ref{descenproches}).
Let us consider $k \in \big[\!|0,M-4|\!\big]$, and $c$ large enough so that
$4 \varepsilon c^{2} >1$. Using the Markov property at time
$k\varepsilon c^{2}$ and the fact that
$(\hat{\xi}_{s}^{y,k\varepsilon c^{2}})_{s \geq 0}$ has the same law as
$(y+ \xi_{s}^{0})_{s \geq 0}$ we get :
 \begin{eqnarray} \label{lem5preuve}
&& P\left( \exists y \in \Omega_{k \varepsilon c^{2}}^{(k+4)\varepsilon c^{2}}
 : \bigcup_{s \in [0, 5\varepsilon c^{2}]}
\hat{\xi}_{s}^{y,k \varepsilon c^{2}} \not\subset
\overline{B}(y,\eta c /2) \right)
\nonumber \\
&& \ \leq  \sum_{y \in \mathbb{Z}^{d}}
P(y \in \xi_{k \varepsilon c^{2}}^{0}) P
 \left( \xi_{4 \varepsilon c^{2}}^{0} \neq \emptyset,
\bigcup_{s \in [0, 5\varepsilon c^{2}]} \xi_{s}^{0} \not\subset
 \overline{B}(0,\eta c /2)  \right)  \nonumber \\
&& \ = \sum_{y \in \mathbb{Z}^{d}} q_{k\varepsilon c^{2}}(y)
p_{4\varepsilon c^{2}} P_{4\varepsilon c^{2}}^{*}
 \left( \sup_{s \leq 5\varepsilon c^{2}} \sup_{x \in \xi_{s}^{0}} |x| >
\frac{\eta c}{2}  \right). \nonumber \\
&& \ \leq p_{4\varepsilon c^{2}}  K_{1} \exp\left( -K_{2} \frac{\eta}{
4\sqrt{\varepsilon}} \right),
\end{eqnarray}
 where at the last line we used
Lemma \ref{speed} with $\alpha = 4 \varepsilon c^{2} > 1$ and
$A= \eta (4\sqrt{\varepsilon})^{-1} $, and the fact that
$\sum_{z \in \mathbb{Z}^{d}}q_{k \varepsilon c^{2}}(y)=1$ from the symmetry
assumption on $p$.
 Since $M-3 \leq 4T (\varepsilon \rho)^{-1}$, we deduce from
the above that
\begin{eqnarray} \label{descend}
&& P_{c^{2}T}^{*} \left( \exists k \in \big[\!|0,M-4|\!\big] \ \exists y \in
\Omega_{k \varepsilon c^{2}}^{(k+4)\varepsilon c^{2}}
 : \bigcup_{s \in [0, 5\varepsilon c^{2}]}
\hat{\xi}_{s}^{y,k \varepsilon c^{2}} \not\subset
\overline{B}(y,\eta c /2) \right)
\nonumber \\
&& \ \ \qquad \qquad
 \leq  p_{c^{2}T}^{-1}  \frac{4T}{\varepsilon \rho}
p_{4\varepsilon c^{2}} K_{1} \exp\left( -K_{2} \frac{\eta}{
4\sqrt{\varepsilon}} \right).
\end{eqnarray}
 Provided $c$ is sufficiently large, we then deduce
(\ref{descenproches}) from (\ref{descend}), (\ref{pt}), and our
choice of $\varepsilon_{0}$.
 \PA
Let us now prove (\ref{assezdescen}). Fix $k \in \big[\!|0,M-4|\!\big]$.
 Using the same arguments as in the proof of (\ref{descenproches}), we obtain
\begin{eqnarray} \label{decomp}
&& P\left(  \exists y \in \Omega_{k \varepsilon c^{2}}^{(k+4)\varepsilon c^{2}}
 : \inf_{s \in [\varepsilon c^{2}, 3 \varepsilon c^{2}]}
|\hat{\xi}_{s}^{y,k \varepsilon c^{2}}| < \delta m_{ c^{2}} \right)
\nonumber \\
&& \ \ \leq  \sum_{y \in \mathbb{Z}^{d}}
q_{k\varepsilon c^{2}}(y)
 P\left(\inf_{s \in [\varepsilon c^{2}, 3\varepsilon c^{2}]}
|\xi_{s}^{0}| < \delta m_{ c^{2}}, \xi_{4\varepsilon c^{2}}^{0} \neq
\emptyset  \right)
\nonumber \\
&& \ \ = p_{4\varepsilon c^{2}}
 P_{4 \varepsilon c^{2}}^{*} \left( \inf_{s \in [\varepsilon c^{2},
3 \varepsilon c^{2}]} |\xi_{s}^{0}| < \delta m_{ c^{2}} \right).
\end{eqnarray}
 Furthermore, by rescaling, when $c$ is large enough so that
$(m_{c^{2}})^{-1}m_{\varepsilon c^{2}} \leq 1$ (recall $\varepsilon <1$ from
our choice of $\varepsilon_{0}$), we get
\begin{eqnarray} \label{scaling4eps}
P_{4 \varepsilon c^{2}}^{*} \left( \inf_{s \in [\varepsilon c^{2},
3 \varepsilon c^{2}]} |\xi_{s}^{0}| < \delta m_{\varepsilon c^{2}}
\right) &= & P_{4\varepsilon c^{2}}^{*}
\left( \inf_{s \in [\varepsilon,3\varepsilon]}
X_{s}^{c^{2},0}(1) < \delta \frac{m_{\varepsilon c^{2}}}{m_{c^{2}}}
\right) \nonumber \\ & \leq & P_{4\varepsilon c^{2}}^{*}
\left( \inf_{s \in [\varepsilon,3\varepsilon]}
X_{s}^{c^{2},0}(1) \leq \delta \right) .
\end{eqnarray}
Since the set $\{ \omega \in D(\mathbb{R}_{+}, M_{F}(\mathbb{R}^{d})) :
\inf_{s \in [\varepsilon, 3\varepsilon]}
\omega_{s}(1) \leq  \delta \}$ is closed
for the Skorohod $J_{1}$ topology,
Theorem \ref{inv} implies
\begin{equation*}
 \limsup_{c \to \infty}
P_{4 \varepsilon c^{2}}^{*} \left( \inf_{s \in [ \varepsilon,
3  \varepsilon]} X_{s}^{c^{2},0}(1) \leq \delta \right)
\leq \mathbb{N}_{0}^{(4\varepsilon)} \left[
\inf_{t \in [\varepsilon,3\varepsilon]} Y_{t}(1) \leq
\delta \right]\leq \frac{\varepsilon^{2}\rho^{2}}{8T^{2}},
\end{equation*}
by our choice of $\delta$.
Assertions (\ref{decomp}), (\ref{scaling4eps}), and the above now imply
\begin{eqnarray*}
&&  \limsup_{c \to \infty} P_{c^{2}T}^{*}
\left(  \exists y \in \Omega_{k \varepsilon c^{2}}^{(k+4)\varepsilon c^{2}}
 : \inf_{s \in [\varepsilon c^{2}, 3 \varepsilon c^{2}]}
|\hat{\xi}_{s}^{y,k \varepsilon c^{2}}| < \delta m_{ c^{2}} \right) \\
&& \ \ \leq \limsup_{c \to \infty}
p_{c^{2}T}^{-1} p_{4\varepsilon c^{2}} \frac{\varepsilon^{2} \rho^{2}}{8 T^{2}}
= \frac{\varepsilon \rho^{2}}{32 T},
\end{eqnarray*}
where we used (\ref{pt}) at the last line.
Hence, using the fact that
$M-3 \leq 4T(\rho \varepsilon)^{-1}$, we get (\ref{assezdescen}),
 provided $c$ is sufficiently large. This ends the proof of Lemma
\ref{enough}. $\qquad \Box$

We have thus finished the proof of the asymptotic upper bound on
 $P_{c^{2}T}^{*}(T_{c[x]_{c}}< \infty)$ (Proposition \ref{uppd3}).
However, to complete the proof of the
asymptotic upper bound for $d=2$ or $3$ in Theorem 1, we need to establish a
corresponding result under the measure $P$. Let us briefly explain how Lemma
\ref{speed} allows us to do so.

\subsection{Back to non-conditioned results.}
First, we shall prove a result corresponding to Lemma \ref{speed} without
conditioning upon survival.
\begin{claim} \label{noncond}
There exists a positive
 $K_{0}$ such that for any $\alpha>1$, for any $A \geq 1$,
\begin{equation*}
P \left(\sup_{t \leq 2 \alpha} \sup_{x\in \xi_{t}^{0}} |x| >
A \sqrt{\alpha} \right) \leq K_{0}p_{\alpha}  \exp(-K_{2} A).
\end{equation*}
\end{claim}
\underline{Proof of Claim \ref{noncond}}:
For any $i \in \big[\!|0,N-1 | \!\big]$ we have
\begin{eqnarray*}
\!\!\!&&\!\!\!
P \left( \sup_{t \leq 2^{1-i}\alpha} \sup_{x\in \xi_{t}^{0}} |x| >
A \sqrt{\alpha} \right) \\
\!\!\!&&\!\!\! = p_{2^{-i} \alpha } P_{2^{-i}\alpha}^{*}
 \left( \sup_{t \leq 2^{1-i}\alpha}
\sup_{x\in \xi_{t}^{0}} |x| > A \sqrt{\alpha}\right)\! + \!
P \left(\sup_{t \leq 2^{-i}\alpha} \sup_{x\in \xi_{t}^{0}} |x| >
A \sqrt{\alpha}, \xi_{2^{-i}\alpha}^{0} \!= \!\emptyset \right) \\
\!\!\!&&\!\!\!
\leq p_{2^{-i} \alpha } K_{1} \exp \left( -K_{2} A 2^{i/2} \right) +
P \left(\sup_{t \leq 2^{-i}\alpha} \sup_{x\in \xi_{t}^{0}} |x| >
A \sqrt{\alpha}, \xi_{2^{-i}\alpha}^{0} = \emptyset \right),
\end{eqnarray*}
where we used Lemma 4 at the last line.
It easily follows that
\begin{eqnarray} \label{sumoveri}
&& P \left(\sup_{t \leq 2 \alpha} \sup_{x\in \xi_{t}^{0}} |x| >
A \sqrt{\alpha} \right)\nonumber \\ && \leq \sum_{i=0}^{N-1}
p_{2^{-i} \alpha } K_{1} \exp \left( -K_{2} A 2^{i/2} \right)
+ P \left(\sup_{t \leq 2^{1-N}\alpha} \sup_{x\in \xi_{t}^{0}} |x| >
A \sqrt{\alpha} \right).
\end{eqnarray}
  Furthermore, by an easy application of Lemma \ref{quickescape},
\begin{eqnarray*}
 && P \left(\sup_{t \leq 2^{1-N}\alpha}
 \sup_{x\in \xi_{t}^{0}} |x| >
A \sqrt{\alpha} \right) \leq 2K_{5} \exp\left( -K_{6}A \sqrt{\alpha}/2 \right).
\end{eqnarray*}
 Thus, from (\ref{sumoveri})
and (\ref{asy}), we obtain
$$ \!\!P \left(\sup_{t \leq 2 \alpha} \sup_{x\in \xi_{t}^{0}} |x| >
A \sqrt{\alpha} \right) \! \leq \! \kappa_{0} K_{1} \sum_{i=0}^{N-1}
 2^{i}p_{\alpha} \exp \left( -K_{2} A 2^{i/2} \right) \!+ \!
2K_{5}\exp(-K_{6}A \sqrt{\alpha}/2), $$
and Claim \ref{noncond} follows. $\qquad \Box$ \vskip0.1cm \PA
Let us now finish the proof of the upper bound in Theorem \ref{main} in
 dimensions $2$ and $3$.
 As before, $x \in \mathbb{R}^{d}\setminus 0$ is fixed.
Simply observe that, for every $T>0$,
\begin{eqnarray} \label{relation}
P(T_{c[x]_{c}}< \infty) =
 P(\xi_{c^{2}T}^{0}= \emptyset, T_{c[x]_{c}}<\infty)
+ p_{c^{2}T} P_{c^{2}T}^{*}(T_{c[x]_{c}}<\infty).
\end{eqnarray}
 On the one hand
\begin{eqnarray*}
  \phi_{d}(c) P(\xi_{c^{2}T}^{0}= \emptyset, T_{c[x]_{c}}<\infty)
&\leq &\phi_{d}(c) P \left(\sup_{t \leq c^{2}T} \sup_{y \in \xi_{t}^{0}} |y| \geq
c|[x]_{c}|  \right) \\
 &\leq & \phi_{d}(c) K_{0}p_{c^{2}T} \exp\left(-K_{2}\frac{|[x]_{c}|}{\sqrt{T}}
\right),
\end{eqnarray*}
where we used Claim \ref{noncond} at the last line. We can now use
(\ref{pt}) to obtain
$$ \limsup_{c \to \infty}
\phi_{d}(c) P(\xi_{c^{2}T}^{0}= \emptyset, T_{c[x]_{c}}<\infty)
\leq   \frac{K_{0}}{\beta_{d}T} \exp\left( -K_{2} \frac{|x|}{2 \sqrt{T}}
\right), $$
which goes to $0$ as $T \to 0$. \PA
On the other hand, using Proposition \ref{uppd3} and (\ref{pt}), we get, for
every $T>0$,
\begin{equation*}
\limsup_{c \to \infty}  \phi_{d}(c) p_{c^{2}T} P_{c^{2}T}^{*}(T_{c[x]_{c}} < \infty)
\leq \frac{1}{T \beta_{d}} \mathbb{N}_{0}^{(T)}(x \in \mathcal{R}),
\end{equation*}
and by (\ref{N0xT}), the right-hand side converges, as $T \to 0$, to
$$ \frac{2\sigma^{2}}{\beta_{d}} \left( 2-\frac{d}{2} \right) |x|^{-2}.$$
 From (\ref{relation}) and the preceeding observations we get
\begin{eqnarray} \label{uppbdd2d3}
 \limsup_{c \to \infty}\phi_{d}(c)P(T_{c[x]_{c}}<\infty)
\leq \frac{2 \sigma^{2}}{\beta_{d}}
\left( 2- \frac{d}{2} \right) |x|^{-2},
\end{eqnarray}
which completes the proof of the upper bound in Theorem 1, in the
 case $d=2$ or $3$. We already noticed that the case $d \geq 5$ follows from
Proposition \ref{uppd5}. Finally, note that in the case $d=4$,
Proposition \ref{uppd3} and a similar proof imply
$$ \limsup_{c \to \infty} c^{2}P(T_{c[x]_{c}}<\infty) =0, $$
as we already mentioned in the introduction.

\section{Lower bound.}
 In this section we finish the proof of Theorem 1 by establishing the required
asymptotic lower bounds on $P(T_{c[x]_{c}}<\infty)$.
 We also prove a similar result in dimension $4$. \PA
  \begin{proposition} \label{lowerbd}
Fix $x \in \mathbb{R}^{d}$, $x \neq 0$.
\begin{description}
\item[{\sc Rough lower bound }] Let $d \geq 4$. There exists
a positive constant $a_{d}$ depending on $|x|$ and $d$ such that
$$ \liminf_{c \to \infty} \phi_{d}(c) P(\exists t \geq 0 : c[x]_{c} \in
\xi_{t}^{0}) \geq  a_{d}, $$
where we recall that for $d \geq 5$, $\phi_{d}(c)=c^{d-2}$, and $\phi_{4}(c)=c^{2}\ln(c)$.
\item[{\sc Sharp lower bound }] Let $d=2$ or $3$. Recall $\phi_{2}(c)=c^{2}(\ln(c))^{-1}$ and $\phi_{3}(c)=c^{2}$. Then
$$ \liminf_{c \to \infty} \phi_{d}(c)  P(\exists t \geq 0 : c[x]_{c} \in
\xi_{t}^{0}) \geq  \frac{2 \sigma^{2}}{\beta_{d}}
\left( 2- \frac{d}{2} \right) |x|^{-2}. $$
\end{description}
\end{proposition}

\subsection{Proof of the rough lower bound, $d \geq 4$}
 For $T>0$ let us introduce the random variable
$$ U_{T} := \int_{0}^{T} \mathbf{1}_{\{c[x]_{c} \in \xi_{c^{2}s}^{0}\}} ds,$$
so that $c^{2}U_{T}$ is the occupation time of opinion $1$ for the voter
at $c[x]_{c}$ in the time interval $[0,c^{2}T]$.  \\
 We clearly have for any $T>0$,
$P(\exists t \geq 0 : c[x]_{c} \in
\xi_{t}^{0}) \geq P(U_{T}>0)$. Using the Cauchy-Schwarz inequality, we
thus obtain
\begin{equation} \label{UT}
P(\exists t \geq 0 : c[x]_{c} \in
\xi_{t}^{0}) \geq \frac{(E[U_{T}])^{2}}{E[(U_{T})^{2}]}.
\end{equation}
Hence, proving the lower bound reduces to establishing
the following two estimates
\begin{eqnarray} \label{EUT}
c^{d} E[U_{T}] \underset{c \to \infty}{\longrightarrow} \int_{0}^{T} ds
p_{s}(x),
\end{eqnarray}
\begin{eqnarray} \label{EUT2}
\limsup_{c \to \infty} c^{2d} \phi_{d}(c)^{-1}E[(U_{T})^{2}] \leq L,
\end{eqnarray}
where $L$ is a constant depending only on $|x|, d$ and $T$. \PA The
first moment of $U_{T}$ is
$$ E[U_{T}] =  \int_{0}^{T} ds P\left( Z_{c^{2}s}^{c[x]_{c},c^{2}s} =0
\right) = \int_{0}^{T} ds q_{c^{2}s}(c[x]_{c}),$$
so that (\ref{EUT}) is a consequence of (\ref{EUT1}).
Let us now estimate the second moment of $U_{T}$. We have
$$   \frac{1}{2} E[(U_{T})^{2}] = \int_{0}^{T} dt \int_{t}^{T} dr  P\left[
c[x]_{c} \in \xi_{c^{2}t}^{0}, c[x]_{c} \in \xi_{c^{2}r}^{0} \right]. $$
Let us fix $r$ and $t$ with $0<t<r \leq T$.
 Using duality over the time interval $[0,c^{2}r]$ and setting
$s : = r - t $, we see that
$P\left[ c[x]_{c} \in \xi_{c^{2}t}^{0}, c[x]_{c} \in \xi_{c^{2}r}^{0} \right]$
 is the probability
for two coalescing random walks starting at point $c[x]_{c}$
respectively at times $0$ and $c^{2}s$, to be both located at point
$0$ at time $c^{2}r$. Using the symmetry properties of $p$ and the
Markov property for the first walk at time $c^{2}s$, we get
\begin{equation} \label{EUT2'}
\frac{1}{2} E[(U_{T})^{2}] = \int_{0}^{T} dt \int_{0}^{T-t} ds
\sum_{y \in \mathbb{Z}^{d}} q_{c^{2}s}(y) \mathbb{P}_{0,y} \left[
T_{1} \leq c^{2}t, Z_{c^{2}t}^{1}= c[x]_{c} \right],
\end{equation}
recalling that the notation $\mathbb{P}_{0,y}$ was introduced in
Section 2.5. \PA With a slight abuse, in the remaining part of the
section we use $L$ to denote a positive constant that only depends
on $T,d$ and $|x|$ and may change from line to line. We suppose that
$c \geq 1 \vee |x|^{-2}$ in order to use Lemma \ref{boundcoal}. Let
us set
$$ I(y) := \int_{0}^{T} dt \left( \int_{0}^{T-t}\!\! ds\
q_{c^{2}s}(y)\right) \mathbb{P}_{0,y}
\left[ T_{1} \leq c^{2}t, Z_{c^{2}t}^{1}=c[x]_{c} \right].$$
 Note that $I(y)$ also depends on $d,T$ and $x$,
although this does not appear in our
notation. From (\ref{EUT2'}), we have
$E[(U_{T})^{2}]=2\sum_{y \in \mathbb{Z}^{d}} I(y)$.
Hence, we need to bound $I(y)$ over different regions of
$\mathbb{Z}^{d}$, in order to control
$c^{2d} (\phi_{d}(c))^{-1} E[(U_{T})^{2}]$.
Recall from Section 2.4 that $\psi_{d}(c)= c^{d-2}$ for $d \geq 4$.
Using (\ref{boundqc2s}) twice in the case $y=0$, and using (\ref{boundqc2s}) together with Lemma \ref{boundcoal} in the case $y \neq 0$,
we get
\begin{itemize}
\item for $y=0$,
\begin{equation*}
 I(0) \leq L c^{-2d} \psi_{d}(c),
\end{equation*}
\item for $y \in \mathbb{Z}^{d}, y \neq 0$,
 \begin{equation*}
 I(y) \leq L c^{-2d} \psi_{d}(c/|y|) \psi_{d}(|y|)^{-1},
\end{equation*}
\item for $y \in \mathbb{Z}^{d}, |y| > c^{2}$,
$$  I(y) \leq L c^{-d-2} \exp(-L_{2}' \sqrt{|y|}) \psi_{d}(|y|)^{-1}.$$
\end{itemize}
Since $\phi_{d}(c) \geq \psi_{d}(c)$ when $d \geq 4$, we then obtain
\begin{equation} \label{I0T}
 c^{2d} (\phi_{d}(c))^{-1} I(0)\leq  L.
\end{equation}
Furthermore, if $y \in \mathbb{Z}^{d} \setminus 0, d \geq 4$, we have
$\psi_{d}(c/|y|) = \psi_{d}(c)(\psi_{d}(|y|))^{-1}$. Hence,
\begin{eqnarray} \label{yleqc2}
 c^{2d} (\phi_{d}(c))^{-1} \!\!\!
\sum_{y \in \mathbb{Z}^{d},0< |y| \leq c^{2} } \!\! I(y) \!\!& \!\!
\leq \!\!&\!\! L (\phi_{d}(c))^{-1} \psi_{d}(c) \!\!\! \sum_{y \in
\mathbb{Z}^{d},0< |y| \leq c^{2} }
 \!\! |y|^{4-2d} \nonumber \\ \!\!&\!\! \leq \!\!&\!\!
 L (\phi_{d}(c))^{-1} \psi_{d}(c) \sum_{k=1}^{c^{2}} k^{d-1} k^{4-2d} \leq L,
\end{eqnarray}
Finally,
 \begin{equation} \label{y>c2}
c^{2d} (\phi_{d}(c))^{-1} \!\!\!\!\!\!\!\!
\sum_{y \in \mathbb{Z}^{d}, |y| > c^{2} }\!\!\!\!\!\!\! I(y)  \leq
 L c^{d-2} (\phi_{d}(c))^{-1} \!\!\!
\sum_{|y|>c^{2}}\!\! |y|^{-d+2} \exp(- L_{2}' \sqrt{|y|}) \!
\underset{c \to \infty}{\longrightarrow} \! 0.
\end{equation}
Since $E[(U_{T})^{2}]=2\sum_{y \in \mathbb{Z}^{d}} I(y)$,
the desired result (\ref{EUT2}) follows from
 (\ref{I0T}), (\ref{yleqc2}) and (\ref{y>c2}).
 This finishes the proof of the
 lower bound for $d \geq 4$. $\quad \Box$ \\
A similar proof in the case $d=2$ or $3$
would give us a rough lower bound, but we need to get sharper
estimates.

\subsection{Outline of the proof of the sharp lower bound, $d=2$ or $3$}
$\qquad\ \ \ \ \ \ $
Fix $x \in \mathbb{R}^{d} \setminus \{ 0\}$. For any $T>0$, $c>0$, we have
$$ \phi_{d}(c) P(T_{c[x]_{c}}< \infty) \geq \phi_{d}(c) p_{c^{2}T}
P_{c^{2}T}^{*}(T_{c[x]_{c}}< \infty).$$
We deduce from (\ref{pt}) that for any $T>0$,
$\lim_{c \to \infty} \phi_{d}(c) p_{c^{2}T} = (\beta_{d}T)^{-1}$ so that
\begin{equation} \label{liminfphi}
  \liminf_{c \to \infty} \phi_{d}(c) P(T_{c[x]_{c}}< \infty) \geq
  (\beta_{d}T)^{-1}  \liminf_{c \to \infty}
P_{c^{2}T}^{*}(T_{c[x]_{c}}< \infty).
\end{equation}
\begin{claim}  \label{lowercl}
For any $\rho>0$, if $T>0$ is sufficiently small,
\begin{equation*} \label{slbP*}
  \liminf_{c \to \infty} P_{c^{2}T}^{*}(T_{c[x]_{c}} < \infty) \geq
(1-\rho) \mathbb{N}_{0}^{(T)}(x \in \mathcal{R}).
\end{equation*}
\end{claim}
The desired lower bound follows from (\ref{liminfphi}), the above
claim and (\ref{N0xT}) by letting $T$ go to $0$. Let us now outline
the proof of Claim
 \ref{lowercl}. \PA For $z \in \mathbb{R}^{d}$,
$\varepsilon ' > \varepsilon >0$, let
$$\mathcal{C}(z,\varepsilon,\varepsilon ')=
\left\{ y \in \mathbb{R}^{d} : \varepsilon < |y-z| < \varepsilon '\right\},
\qquad h(\varepsilon) : = \varepsilon^{2} \ln(\ln(\varepsilon^{-1})). $$
We also set for $r \in (0,1)$
\begin{eqnarray*}
 g_{d}(r) = \begin{cases} 2^{-\left(\frac{\ln(r)}{\ln(2)}\right)^{4}} &
\mbox{ if } d=2 ,\\
 r^{16} & \mbox{ if } d=3. \end{cases}
\end{eqnarray*}
For $\alpha >0, T>0$ and
$\varepsilon_{0} \in (0,1)$ we consider the events
\begin{eqnarray*}
&& \mathcal{E}_{\varepsilon_{0}}^{(c)} = \left\{
 \exists s \geq 0 \ \exists \varepsilon \in
(g_{d}(\varepsilon_{0}), \varepsilon_{0})
: X_{s}^{0}
\left(\mathcal{C}\left(c[x]_{c},\frac{c \varepsilon}{4},2c\varepsilon\right)
\right) \geq \alpha h(\varepsilon) \phi_{d}(c) \right\}, \\
&& \mathcal{F}_{\varepsilon_{0}}^{(c)} = \left\{
 \exists s \geq 0 \ \exists \varepsilon \in
(g_{d}(\varepsilon_{0}), \varepsilon_{0})
: X_{s}^{0}
\left(\mathcal{C}\left(cx,\frac{c \varepsilon}{2},c\varepsilon\right)\right)
 > \alpha h(\varepsilon) \phi_{d}(c) \right\}.
\end{eqnarray*}
For $c$ large enough, we have $\mathcal{F}_{\varepsilon_{0}}^{(c)}
\subset \mathcal{E}_{\varepsilon_{0}}^{(c)}$, hence
\begin{equation} \label{conditio}
 P_{c^{2}T}^{*} (T_{c[x]_{c}}< \infty)
\geq P_{c^{2}T}^{*}
\left( \mathcal{F}_{\varepsilon_{0}}^{(c)}
  \right) \times
P_{c^{2}T}^{*} \left( T_{c[x]_{c}}<\infty  \bigg|
\mathcal{E}_{\varepsilon_{0}}^{(c)} \right).
\end{equation}
 The idea of the proof of Claim \ref{lowercl} is the following. Rescaling and
using Theorem \ref{inv}, we will show that for $c$ large, the
first term of the product in (\ref{conditio}), namely
$P_{c^{2}T}^{*}( \mathcal{F}_{\varepsilon_{0}}^{(c)} )$,
 is bounded below by a
corresponding rescaled quantity under $\mathbb{N}_{0}^{(T)}$. For
$\alpha$ small enough, this quantity will then be bounded from below
by a quantity arbitrarily close to $\mathbb{N}_{0}^{(T)}( x \in
\mathcal{R})$ (see assertions (\ref{limittheorem1}) and
(\ref{couronnecond}) below). To finish the proof of Claim
\ref{lowercl} we shall then
 establish that if we take
$T, \varepsilon_{0}$ small enough, the second term of the product
 in (\ref{conditio}), namely $P_{c^{2}T}^{*} \left( T_{c[x]_{c}}<\infty  \big|
\mathcal{E}_{\varepsilon_{0}}^{(c)} \right)$,
is, for $c$ large, arbitrarily close to $1$ (see Lemma
\ref{highprob} below). \PA Let us reformulate the preceeding discussion
in more precise terms. Using rescaling we have
\begin{equation} \label{EsousP*}
  P_{c^{2}T}^{*}
\left( \mathcal{F}_{\varepsilon_{0}}^{(c)}
  \right) =
 P_{c^{2}T}^{*} \left( \exists s \geq 0 \ \exists \varepsilon \in
(g_{d}(\varepsilon_{0}), \varepsilon_{0}) : X_{s}^{c^{2},0}
\left(\mathcal{C}(x,\frac{\varepsilon}{2},\varepsilon)\right)
 > \alpha h(\varepsilon)\right).
\end{equation}
 It is easy to see that the set
$$ \left\{ \omega \in
D(\mathbb{R}_{+}, M_{F}(\mathbb{R}^{d})) : \exists s \geq 0 \ \exists
\varepsilon \in (g_{d}(\varepsilon_{0}), \varepsilon_{0}), \  \omega_{s}
\left(\mathcal{C}(x,\frac{\varepsilon}{2},\varepsilon)
\right)  > \alpha h(\varepsilon)
 \right\}$$
is open for the Skorohod $J_{1}$ topology. Theorem \ref{inv} thus implies that
\begin{equation} \label{limittheorem1}
 \liminf_{c\to \infty} P_{c^{2}T}^{*}
\left( \mathcal{F}_{\varepsilon_{0}}^{(c)}
  \right)\! \geq  \mathbb{N}_{0}^{(T)}
( \exists s \geq 0 \ \exists \varepsilon \in (g_{d}(\varepsilon_{0}),
\varepsilon_{0}) : Y_{s}(C(x,\frac{\varepsilon}{2},\varepsilon)) > \alpha
h(\varepsilon) ).
\end{equation}
\begin{lemma} \label{couronne}
Let $d=2$ or $3$. We can choose
$\alpha > 0$ so that, for any $\delta>0$, there exists
$\varepsilon_{1}\in \left(0,1 \wedge \frac{|x|}{2}\right)$
such that for any $\varepsilon_{0} \in (0,\varepsilon_{1})$,
\begin{equation*}
 \mathbb{N}_{0} \left[
\exists s\geq 0 \ \exists \ \varepsilon \in
(g_{d}(\varepsilon_{0}),\varepsilon_{0}) : Y_{s}\left(
\mathcal{C}\left(x,\frac{\varepsilon}{2}, \varepsilon \right)\right) > \alpha
h(\varepsilon) \ \bigg|\ x \in \mathcal{R} \right] \geq 1-\delta.
\end{equation*}
\end{lemma}
In the following, we fix $\alpha$ as in Lemma \ref{couronne}.
The following lemma
 estimates the second term of the product in the right-hand side of
(\ref{conditio}).
 \begin{lemma} \label{highprob}
For any fixed $\gamma>0$,
there exists $\varepsilon_{2} \in (0,1)$ such that for any
$\varepsilon_{0} \in (0,\varepsilon_{2})$, we have
\begin{enumerate}
\item[(a)] \quad
$\underset{c \to \infty}{\liminf} \
 P\left(T_{c[x]_{c}} < \infty \  \big| \ \mathcal{E}_{
\varepsilon_{0}}^{(c)} \right) \geq 1- \gamma,$
\vskip0.2cm
\item[(b)] \quad
$ \underset{T \to 0}{\liminf} \ \left( \underset{c \to \infty}{\liminf} \
 P_{c^{2}T}^{*} \left( T_{c[x]_{c}}<\infty  \big|
\mathcal{E}_{\varepsilon_{0}}^{(c)}
 \right) \right) \geq 1-\gamma.$
\end{enumerate}
\end{lemma}
Let us now fix $\delta \in (0,1)$, $\gamma >0$
and let
$\varepsilon_{1}$ and $\varepsilon_{2}$ be as in Lemma \ref{couronne} and
 \ref{highprob} respectively.
Since
$$\mathbb{N}_{0} \left( \left\{ \sup \{ |y| : y \in \mathcal{R} \} >
\frac{|x|}{2}\right\} \cap
\bigg\{Y_{T}=0 \bigg\}  \right) \underset{T \to 0}{\longrightarrow} 0, $$
we deduce from Lemma \ref{couronne} that for any
$\varepsilon_{0} \in (0,\varepsilon_{1})$, for $T>0$ sufficiently small,
\begin{equation} \label{couronnecond}
 \mathbb{N}_{0}^{(T)}
\left[ \exists s \geq 0 \ \exists \varepsilon \in
(g_{d}(\varepsilon_{0}),\varepsilon_{0})
 : Y_{s}\left(C\left(x,\frac{\varepsilon}{2},\varepsilon\right)\right)
> \alpha h(\varepsilon) \right]  \geq (1-2\delta)\mathbb{N}_{0}^{(T)}(x \in \mathcal{R}).
\end{equation}
From (\ref{limittheorem1}) and (\ref{couronnecond}), we have for $T$ sufficiently small
\begin{equation}\label{deltaalpha}
\liminf_{c \to \infty} P_{c^{2}T}^{*}(\mathcal{F}_{\varepsilon_{0}}^{(c)})
\geq (1-2\delta) \mathbb{N}_{0}^{(T)}(x \in \mathcal{R}).
\end{equation}
Now use (\ref{conditio}) and Lemma \ref{highprob} (b) to get for $T$ small,
$$  \liminf_{c \to \infty}
 P_{c^{2}T}^{*} (T_{c[x]_{c}} <\infty)
\geq (1-2\delta) (1 - 2\gamma) \mathbb{N}_{0}^{(T)}(x \in \mathcal{R}),$$
which gives
 Claim \ref{lowercl}, hence Proposition \ref{lowerbd}.
\PA To complete our proof of the lower bound Proposition \ref{lowerbd},
we still need to establish Lemma \ref{couronne} and Lemma \ref{highprob}.
 Establishing that part (b) of Lemma \ref{highprob} follows from part (a)
requires a result which is a consequence of Lemma \ref{couronne}. However, we first
give the proof of Lemma \ref{highprob} (Section 4.3 below), because it is more
closely related to our results. We then provide a proof of Lemma \ref{couronne}
in Section 4.4. In these two sections, we will assume for simplicity that $\sigma =1$. Adapting
the proofs to a general $\sigma$ is easy.

 \subsection{Proof of Lemma \ref{highprob}}
We assume in this section that Lemma \ref{couronne} has been proved, and in particular that
(\ref{deltaalpha}) holds.
Let us first explain how to derive part (b) from part (a). We have
 \begin{eqnarray} \label{condtononcond}
&& P_{c^{2}T}^{*}\left(T_{c[x]_{c}} < \infty \ \bigg| \ \mathcal{E}_{
\varepsilon_{0}}^{(c)} \right) \nonumber \\&& = \frac{P\left(
\{T_{c[x]_{c}} < \infty \} \cap \mathcal{E}_{
\varepsilon_{0}}^{(c)}\right)}{P\left( \mathcal{E}_{
\varepsilon_{0}}^{(c)}\cap \{\xi_{c^{2}T}^{0} \neq \emptyset \}
\right)}
-\frac{P\left(
\{T_{c[x]_{c}} < \infty \} \cap \mathcal{E}_{\varepsilon_{0}}^{(c)}
\cap\{ \xi_{c^{2}T} = \emptyset\}\right)}
{P\left( \mathcal{E}_{\varepsilon_{0}}^{(c)}
\cap \{\xi_{c^{2}T}^{0} \neq \emptyset\}
\right)} \nonumber \\
&& \geq  P\left(T_{c[x]_{c}}<\infty \ \bigg| \ \mathcal{E}_{
\varepsilon_{0}}^{(c)} \right) - \frac{P\left(
\{T_{c[x]_{c}} < \infty \}
 \cap\{ \xi_{c^{2}T} = \emptyset\}\right)}{p_{c^{2}T}
P_{c^{2}T}^{*}(\mathcal{E}_{\varepsilon_{0}}^{(c)})}.
\end{eqnarray}
Take $\delta = 1/2$ in Lemma \ref{couronne}, and choose $\varepsilon_{1}$ so
that the conclusion of this lemma holds. From the fact that
$\mathcal{F}_{\varepsilon_{0}}^{(c)} \subset
\mathcal{E}_{\varepsilon_{0}}^{(c)}$, and then from (\ref{deltaalpha}), we get that
for $\varepsilon_{0} \in (0,\varepsilon_{1})$, for $T>0$ small,
\begin{equation} \label{bordeldelta}
 \liminf_{c \to \infty} T^{-1} P_{c^{2}T}^{*}
(\mathcal{E}_{\varepsilon_{0}}^{(c)})
 \geq  \frac{1}{2T} \mathbb{N}_{0}^{(T)}(x \in \mathcal{R}) \geq
\frac{1}{2 \beta_{d}} \left( 2- \frac{d}{2} \right) |x|^{-2},
\end{equation}
using (\ref{N0xT}).
Since $|[x]_{c}| > |x|/2$ for $c$ large enough, we have
$$ P\left( \{T_{c[x]_{c}}< \infty \}\cap \{ \xi_{c^{2}T}^{0} = \emptyset \}
\right) \leq P\left( \sup_{t \leq c^{2}T} \sup_{y \in \xi_{t}^{0}}
|y| > c|x|/2\right).$$ For $T$ small enough so that
$\frac{|x|}{\sqrt{2T}} \geq 1$, and $c$ large enough so that
$\frac{c^{2}T}{2} >1$, we can use Claim \ref{noncond} to deduce that
$$ P\left(
\{T_{c[x]_{c}} < \infty \}
 \cap \{ \xi_{c^{2}T} = \emptyset\} \right) \leq  p_{c^{2}T} K_{0}
\exp\left(-K_{2}\frac{|x|}{\sqrt{2T}} \right).$$
Combining (\ref{bordeldelta}) and this last inequality, we obtain that for
 any $\varepsilon_{0} \in (0,\varepsilon_{1})$
\begin{eqnarray*}
&& \limsup_{T \to 0} \limsup_{c \to \infty}
\frac{P\left(
\{T_{c[x]_{c}} < \infty \} \cap\{ \xi_{c^{2}T} = \emptyset\}\right)}
{ p_{c^{2}T}P_{c^{2}T}^{*}\left( \mathcal{E}_{\varepsilon_{0}}^{(c)}
\right)} = 0.
\end{eqnarray*}
It is now clear from (\ref{condtononcond}) and the above that part
(b) of Lemma \ref{highprob} follows from part (a).  $\qquad \Box$
\vskip0.2cm \PA Proving part (a) of Lemma \ref{highprob} requires
the following intermediate result. Recall that for $A \subset
\mathbb{Z}^{d}$ the voter model $\xi$ starts from $\xi_{0} =A$ under
$P_{A}$.
\begin{lemma} \label{variable:u}
For any $\gamma>0$, there exists $M>0$ and $U>0$ such that for every $u \geq U$
and every subset $A$ of $\mathbb{Z}^{d} \cap \mathcal{C}(0, \frac{u}{4},2u)$ with
$|A| \geq M \phi_{d}(u)$, one has
$$ P_{A}(\exists t \geq 0 : 0 \in \xi_{t}) \geq 1-\gamma. $$
\end{lemma}
\PA Let us postpone the proof of Lemma \ref{variable:u}
and proceed to the proof of Lemma \ref{highprob} (a). \PA
Let us fix $\gamma>0$. Let us choose $M>0$ and $U>0$
such that the conclusion of Lemma \ref{variable:u} holds. We can then choose
$\varepsilon_{2}>0$
small enough so that $$ \frac{\alpha}{2} \ln \left(\ln
\left( \frac{1}{\varepsilon_{2}} \right) \right) \geq M. $$
Let us fix $\varepsilon_{0} \in (0, \varepsilon_{2})$. Let $c>0$ be large enough so that
$$ c g_{d}(\varepsilon_{0}) \geq U, \ \ \ \ \mbox{ and } \ \ \ \  2 \ln(c g_{d}(\varepsilon_{0})) \geq \ln(c).$$
We then set
$$T_{\varepsilon_{0}}^{(c)}:= \inf \left\{ t \geq 0 :
X_{t}^{0}\left(\mathcal{C}\left(c[x]_{c},\frac{c\varepsilon}{4},{2c\varepsilon}
\right)\right) \! \geq \! \alpha
h(\varepsilon) \phi_{d}(c)   \mbox{ for some }
 \varepsilon \in (g_{d}(\varepsilon_{0}),\varepsilon_{0}) \right\}\!.$$
Clearly
$T_{\varepsilon_{0}}^{(c)}$ is a stopping time of the filtration
generated by the voter model, and
 $\mathcal{E}_{\varepsilon_{0}}^{(c)} =
\{T_{\varepsilon_{0}}^{(c)}< \infty\}$. \PA
From the definition of $T_{\varepsilon_{0}}^{(c)}$,
on the event $\mathcal{E}_{\varepsilon_{0}}^{(c)}$ we can
choose a random \\ $\varepsilon \in (g_{d}(\varepsilon_{0}),\varepsilon_{0})$
such that
 $$X_{T_{\varepsilon_{0}}^{(c)}}^{0}
\left(\mathcal{C}\left(
c[x]_{c},\frac{c\varepsilon}{4},2c\varepsilon \right)\right)
 \geq \alpha h(\varepsilon)\phi_{d}(c).$$
On the event $\mathcal{E}_{\varepsilon_{0}}^{(c)}$, we can consider the set
$$ \mathcal{A}_{\varepsilon}^{(c)} := \left(
\xi_{T_{\varepsilon_{0}}^{(c)}}^{0} \cap
\mathcal{C}\left(c[x]_{c},\frac{c\varepsilon}{4}, 2c\varepsilon \right)
\right) - c[x]_{c}, $$
where, for $A \subset \mathbb{Z}^{d}$, $z \in \mathbb{Z}^{d}$,
 $A-z = \{ y \in \mathbb{Z}^{d} : y+z \in A\}$.
The random set $\mathcal{A}_{\varepsilon}^{(c)}$
is a subset of
$\mathbb{Z}^{d} \cap \mathcal{C}(0,\frac{c \varepsilon}{4},2 c \varepsilon)$,
and has cardinality
$|\mathcal{A}_{\varepsilon}^{(c)}| =
\lceil \alpha h(\varepsilon) \phi_{d}(c)\rceil$. \PA
Let us argue on $\mathcal{E}_{\varepsilon_{0}}^{(c)}$ and set
$u(\varepsilon,c)= c \varepsilon$.
Note that $u(\varepsilon,c) \geq c g_{d}(\varepsilon_{0}) \geq U$.
\\ When $d=3$, we have
$$ | \mathcal{A}_{\varepsilon}^{(c)}| \geq \alpha h(\varepsilon) c^{2} \geq
\alpha \ln\left( \ln \left( \frac{1}{\varepsilon_{2}} \right)
\right) \varepsilon^{2} c^{2} \geq M \phi_{d}(u(\varepsilon,c)).$$
When $d=2$, noticing that $2\ln(c\varepsilon) \geq
2\ln(cg_{d}(\varepsilon_{0})) \geq \ln(c)$, we also have
$$   | \mathcal{A}_{\varepsilon}^{(c)}| \geq
\alpha h(\varepsilon) \frac{c^{2}}{\ln(c)} \geq \alpha
\ln\left( \ln \left( \frac{1}{\varepsilon_{2}} \right) \right)
\frac{c^{2}\varepsilon^{2}}{2 \ln(c \varepsilon )} \geq
 M \phi_{d}(u(\varepsilon,c)).$$
From Lemma \ref{variable:u}, we deduce that, on the event
$\mathcal{E}_{\varepsilon_{0}}^{(c)}$,
\begin{equation*}
  P_{\mathcal{A}_{\varepsilon}^{(c)}}
( \exists t \geq 0 : 0 \in \xi_{t}) \geq 1-\gamma.
\end{equation*}
\PA  Using the strong Markov property for $\xi^{0}$ at time
$T_{\varepsilon_{0}}^{(c)}$, then the fact that \\
$\mathcal{A}_{\varepsilon}^{(c)}+c[x]_{c}
\subset \xi_{T_{\varepsilon_{0}}^{(c)}}^{0}$, we obtain
\begin{eqnarray*}
 P \left(
 \mathcal{E}_{\varepsilon_{0}}^{(c)} \cap \{ T_{c[x]_{c}} < \infty \}\right)
& \geq & E \left(\mathbf{1}_{\{T_{\varepsilon_{0}}^{(c)}<\infty\}}
 P_{\xi_{T_{\varepsilon_{0}}^{(c)}}^{0}}
( \exists t \geq 0 : c[x]_{c} \in \xi_{t})  \right) \nonumber \\
& \geq & E \left( \mathbf{1}_{\{T_{\varepsilon_{0}}^{(c)}<\infty\}}
P_{\mathcal{A}_{\varepsilon}^{(c)}}
( \exists t \geq 0 : 0 \in \xi_{t})  \right) \\
& \geq & (1-\gamma) P(\mathcal{E}_{\varepsilon_{0}}^{(c)}),
\end{eqnarray*}
which gives part (a) of Lemma \ref{highprob}.
$\qquad \Box$
\vskip0.2cm \PA
Let us now fix $\gamma>0$ and establish Lemma \ref{variable:u}.
First, notice that the function $A \to P_{A}(\exists t \geq 0 : 0 \in \xi_{t})$
is increasing.
It thus suffices to find $M>0$ and $U>0$ such that for $u \geq U$,
\begin{equation} \label{card=Mphid}
 \inf \left\{
P_{A}(\exists t \geq 0 : 0 \in \xi_{t}) : A \subset \mathbb{Z}^{d} \cap
\mathcal{C}\left(0, \frac{u}{4}, 2u \right), |A| = \lceil M \phi_{d}(u) \rceil
\right\} \geq 1-\gamma.
\end{equation}
For $M>0, u>0$ let us introduce
$$ \mathfrak{A}_{u}^{(M)} : = \left\{
 A \subset \mathbb{Z}^{d} \cap
\mathcal{C}\left(0,\frac{u}{4},2u\right) \ : \ |A| = \lceil M
\phi_{d}(u)\rceil \right\}.$$
We then use a similar method as for
establishing the rough lower bound.
 Let us set $V_{T} = \int_{0}^{T}
\mathbf{1}_{\{0 \in \xi_{u^{2}t}\}}dt$.
As in Section 4.1, we use the Cauchy-Schwarz inequality to get for any $u>0$,
$A \in \mathfrak{A}_{u}^{(M)}$
\begin{equation} \label{moments12}
P_{A} ( \exists t \geq 0 : 0 \in \xi_{t}) \geq
\frac{(E_{A}[V_{T}])^{2}}{E_{A}[(V_{T})^{2}]}.
\end{equation}
We will verify that for any fixed $M>0$, there exists a constant $U(M)$ such that
if $u \geq U$, then for any $A \in \mathfrak{A}_{u}^{(M)}$, we have
\begin{equation} \label{1mom}
E_{A}[V_{T}] \geq \frac{M}{2}
(\psi_{d}(u))^{-1} \int_{T/2}^{T} \frac{1}{(2\pi s)^{d/2}}
\exp\left(- \frac{2}{ s} \right) ds,
\end{equation}
\begin{equation} \label{2mom}
E_{A}[(V_{T})^{2}] \leq \left( E_{A}[V_{T}] \right)^{2} +
L' M \left(\psi_{d}(u)\right)^{-2},
\end{equation}
where $L'$ is a constant depending only on $d$ and $T$.
Let us postpone the proof of these two assertions and finish the proof of
Lemma \ref{variable:u}.
 We can choose $M>0$ sufficiently large so that
$$ L' M \leq \gamma \frac{M^{2}}{4}
\left( \int_{T/2}^{T} ds \frac{1}{(2\pi s)^{d/2}} \exp \left( -
\frac{2}{s} \right) \right)^{2}. $$
From (\ref{1mom}) and (\ref{2mom}), we then deduce that for $u \geq U(M)$,
 for any $A \in \mathfrak{A}_{u}^{(M)}$, we have
$$ E_{A}[(V_{T})^{2}] \leq (1+\gamma) \left( E_{A}[V_{T}] \right)^{2}.
$$
 Lemma \ref{variable:u} now follows from (\ref{moments12}).
\vskip0.2cm
\underline{Proof of (\ref{1mom})}:
Let us now fix $M>0$ and
 establish that (\ref{1mom}) is valid
 for $u$ sufficiently large, and for any $A \in \mathfrak{A}_{u}^{(M)}$.
Note that for any $A \in \mathfrak{A}_{u}^{(M)}$, $u^{-1}A$ is a
subset of  $u^{-1}\mathbb{Z}^{d} \cap \mathcal{C}(0,1/4,2)$ and has
cardinality $\lceil M \phi_{d}(u) \rceil$.
 For any $u>0$, $A \in \mathfrak{A}_{u}^{(M)}$, duality gives
\begin{equation} \label{1momentprel}
E_{A}[V_{T}] =
\int_{0}^{T} \mathbb{P}_{0}(Z_{u^{2}t} \in
A )dt  =
\!\!\sum_{y \in u^{-1}A}
\int_{0}^{T} \! q_{u^{2}t}(u y) dt \geq \!\!\sum_{y \in u^{-1}A}
\int_{T/2}^{T} \! q_{u^{2}t}(u y) dt.
\end{equation}
It is easy to deduce from Theorem \ref{LLT} that uniformly in $t \in [T/2,T]$,
$$ \lim_{u \to \infty} \sup_{y \in u^{-1}\mathbb{Z}^{d}}
|u^{d} q_{u^{2}t}(u y) - p_{t}(y)| =0.$$
Thus, if $u$ is sufficiently large, for any $y \in u^{-1}\mathbb{Z}^{d}$ with
$|y| \leq 2$,
$$ u^{d} q_{u^{2}t}(u y) \geq \frac{1}{2} p_{t}(y) \geq \frac{1}{2}
(2\pi t)^{-d/2} \exp\left(-\frac{2}{t}\right).$$
We deduce from the above and (\ref{1momentprel}) that for $u$ large enough,
and for any $A \in \mathfrak{A}_{u}^{(M)}$, we have
$$E_{A}[V_{T}] \geq \frac{1}{2} u^{-d}|A|
\int_{T/2}^{T} (2\pi s)^{-d/2} \exp\left(-\frac{2}{s}\right) dt.$$
 (\ref{1mom}) now follows from the fact that $u^{-d}\phi_{d}(u)=(\psi_{d}(u))^{-1}$.
\vskip0.1cm
\underline{Proof of (\ref{2mom})}: Let us now
estimate the second moment of $V_{T}$ and prove (\ref{2mom}). Using the same
arguments as in the proof of the rough estimate, we obtain for any $u>0$,
$A \in \mathfrak{A}_{u}^{(M)}$,
\begin{equation*}
 E_{A}[(V_{T})^{2}] = 2 \int_{0}^{T}\!\!\! dt
\int_{0}^{T-t}\!\!\!\! ds \sum_{z' \in \mathbb{Z}^{d}}
q_{u^{2}s}(z') \mathbb{E}_{0,z'} [Z_{u^{2}t}^{1} \in A,
Z_{u^{2}t}^{2} \in A].
\end{equation*}
It follows that
\begin{equation} \label{2moment}
E_{A}[(V_{T})^{2}] = 2 \int_{0}^{T} dt
\int_{0}^{T-t} ds
\sum_{z' \in \mathbb{Z}^{d}} q_{u^{2}s}(z') (H_{1,u}(t,z')+
H_{2,u}(t,z')),
\end{equation}
where
\begin{eqnarray*}
&& H_{1,u}(t,z') := \sum_{y \in A  }
\mathbb{P}_{0,z'}[T_{1}\leq u^{2}t, Z_{u^{2}t}^{1}=y], \\
&& H_{2,u}(t,z') := \sum_{y \in A}
\sum_{y' \in A}
\mathbb{P}_{0,z'}[T_{1} > u^{2}t, Z_{u^{2}t}^{1}=y,Z_{u^{2}t}^{2}=y' ].
\end{eqnarray*}
Since two coalescing walks behave independently before they meet,
we can bound
$\mathbb{P}_{0,z'}[T_{1} > u^{2}t, Z_{u^{2}t}^{1}=y,Z_{u^{2}t}^{2}=y' ]$
by $q_{u^{2}t}(y) q_{u^{2}t}(z'-y') $ so we obtain
\begin{eqnarray} \label{HtildeH}
&& 2 \int_{0}^{T} dt
\int_{0}^{T-t} ds
\sum_{z' \in \mathbb{Z}^{d}} q_{u^{2}s}(z') H_{2,u}(t,z') \nonumber \\
& \leq & 2 \int_{0}^{T} dt \int_{0}^{T-t} ds
 \sum_{z' \in \mathbb{Z}^{d}} q_{u^{2}s}(z')
\sum_{y \in A} \sum_{y' \in A}
q_{u^{2}t}(y) q_{u^{2}t}(z'-y') \nonumber \\
& = & 2 \int_{0}^{T} dt \int_{0}^{T-t} ds \sum_{y \in A}
\sum_{y' \in A} q_{u^{2}t}(y) q_{u^{2}(t+s)}(y')
\nonumber \\
& = & \left(\int_{0}^{T} dt \sum_{y \in A} q_{u^{2}t}(y) \right)^{2} =
\left( E_{A}[V_{T}] \right)^{2}
\end{eqnarray}
With a slight abuse of notation, in the remaining part of the section we use
$L$ to denote a constant depending only on $d$ and $T$ and which may change
from line to line. \PA
Using (\ref{boundqc2s}), we obtain
\begin{eqnarray} \label{decoupen0,u}
&& \int_{0}^{T} dt \int_{0}^{T-t} ds
\sum_{z' \in \mathbb{Z}^{d}} q_{u^{2}s}(z') H_{1,u}(t,z') \nonumber
 \\&& \leq
L_{2} u^{-d} \bigg[ L_{2}\psi_{d}(u)u^{-d} \sum_{y \in A}
\psi_{d}(|y|^{-1}u) + \!\!\!
\sum_{z'\in \mathbb{Z}^{d}, 0<|z'|\leq u} \psi_{d}\left(\frac{u}{|z'|}\right) \int_{0}^{T} H_{1,u}(t,z')  \nonumber \\
&& \qquad \quad \qquad \qquad \qquad \qquad \qquad +\sum_{z'\in
\mathbb{Z}^{d}, |z'|>u} \exp\left(-L_{2}'\frac{|z'|}{u} \right)
\int_{0}^{T} H_{1,u}(t,z')\bigg].
\end{eqnarray}
Note that we used (\ref{boundqc2s}) a second time to bound
$\int_{0}^{T} dt H_{1,u}(t,0)$ and get the first term in the sum
above. Then, Lemma \ref{boundcoal} is exactly what we need to bound
$\int_{0}^{T} H_{1,u}(t,z')$ when $z' \neq 0$. However, the cases
$d=2$ and $d=3$ are slightly different. \PA When $d=3$, for any
$u>0$, for any $A \in \mathfrak{A}_{u}^{(M)}$, and any $z'\neq 0$,
we obtain from Lemma \ref{boundcoal} that
\begin{eqnarray} \label{express1pointd=3}
 \int_{0}^{T} dt  H_{1,u}(t,z')
& \leq & L_{4} u^{-3} \left(\psi_{3}(|z'|\vee 1)\right)^{-1}
 \sum_{y \in u^{-1}A} \psi_{3}\left( \frac{1}{|y|}\right)
\end{eqnarray}
From the fact that $\min_{y \in u^{-1}A} |y| >1/4$ and $|A|=\lceil M \phi_{3}(u) \rceil= \lceil M u^{2} \rceil$,
we obtain that
 \begin{eqnarray*}
\sum_{y \in u^{-1}A} \psi_{3}\left( \frac{1}{|y|}\right)   \leq  4  \lceil M u^{2} \rceil.
\end{eqnarray*}
Hence, from (\ref{decoupen0,u}) and (\ref{express1pointd=3}), we deduce that, when $d=3$, for $u$
sufficiently large, and for any $A \in \mathfrak{A}_{u}^{(M)}$,
\begin{eqnarray*}
&& \int_{0}^{T} dt \int_{0}^{T-t} ds
\sum_{z' \in \mathbb{Z}^{d}} q_{u^{2}s}(z') H_{1,u}(t,z')
\\& \leq & L u^{-6}  M u^{2}  \left[ M + \!\!\!\sum_{z'\in \mathbb{Z}^{3}, 0< |z'| \leq u}\!\!\!
|z'|^{-2} + \!\!\!\sum_{z' \in \mathbb{Z}^{3},|z'| \geq u}\!\!\!
|z'|^{-1} \exp\left( -L_{2}'\frac{|z'|}{u}\right) \right]
\\ & \leq & L  u^{-6}  M u^{4} = L M (\psi_{3}(u))^{-2},
\end{eqnarray*}
where we used that $\sum_{z' \in \mathbb{Z}^{3},|z'| \geq u}
|z'|^{-1} \exp\left( -L_{2}' \frac{|z'|}{u} \right) \leq L
\int_{u}^{\infty} \rho \exp(-L_{2}' \rho/u) d\rho$.
 \PA
Similarly, when $d=2$, for any $u>1$, for any $A \in \mathfrak{A}_{u}^{(M)}$, and any $z'\neq 0$, we get
from Lemma \ref{boundcoal} that
\begin{eqnarray} \label{express1pointd=2}
 \int_{0}^{T} dt  H_{1,u}(t,z')
& \leq & L_{4} u^{-2} \frac{\psi_{2}(u/|z'|)}{\psi_{2}(|z'|)}
\sum_{y \in u^{-1}A} \ln(|y|^{-1}) \nonumber
 \\ & \leq & \ln(4) L_{4}  u^{-2} \lceil M \frac{u^{2}}{\ln(u)} \rceil \frac{\psi_{2}(u/|z'|)}{\psi_{2}(|z'|)},
\end{eqnarray}
Furthermore, we have $\min_{y \in u^{-1}\!A} |y| >1/4$ and
$|A|=\lceil M \phi_{2}(u) \rceil= \lceil M u^{2}/ (2\ln(u)) \rceil$,
so that
\begin{eqnarray*}
\sum_{y \in u^{-1}A} \psi_{2}\left( \frac{1}{|y|}\right)   \leq  4  \lceil M u^{2}/2(\ln(u)) \rceil.
\end{eqnarray*}
Hence, from (\ref{decoupen0,u}) and (\ref{express1pointd=2}), we deduce that for any $u$ sufficiently large,
for any $A \in \mathfrak{A}_{u}^{(M)}$,
\begin{eqnarray*}
&& \int_{0}^{T} dt \int_{0}^{T-t} ds
\sum_{z' \in \mathbb{Z}^{d}} q_{u^{2}s}(z') H_{1,u}(t,z')
\\& \leq &
L u^{-4}  M \frac{u^{2}}{2\ln(u)}  \bigg[ \ln(u) + \!\!\!\sum_{z'\in
\mathbb{Z}^{3}, 0< |z'| \leq u}\!\!\! \ln\left(\frac{u}{|z'|}\vee
e\right)^{2} \times \frac{1}{\ln(|z'|\vee e)}
\\ && \qquad \qquad \qquad \qquad \qquad \qquad \qquad  +
\!\!\!\sum_{z' \in \mathbb{Z}^{3},|z'| \geq u}\!\!\!
\frac{1}{\ln(|z'|)} \exp\left( -L_{2}'\frac{|z'|}{u}\right) \bigg]
\\ & \leq & L  u^{-4}  M \frac{u^{2}}{2\ln(u)}  \left[ \ln(u) +
\int_{\sqrt{2}}^{u} \rho \frac{(\ln(u/\rho))^{2}}{\ln(\rho)} d\rho +
\int_{u}^{\infty} \frac{\rho}{\ln(\rho)} \exp(-L_{2}' \rho/u) d\rho
\right]
\\ & \leq & L u^{-2} (\ln(u))^{-1}  M \frac{u^{2}}{2\ln(u)} = 2 L M (\psi_{2}(u))^{-2}.
 \end{eqnarray*}
To get to the last line above, we have used elementary computations
to check that for $u$ large, one has
$$\int_{\sqrt{2}}^{u} \rho \frac{(\ln(u/\rho))^{2}}{\ln(\rho)} d\rho
\leq \frac{1}{2} \frac{u^{2}}{\ln(u)}, \quad \int_{u}^{\infty}
\frac{\rho}{\ln(\rho)} \exp(-L_{2}' \rho/u) d\rho \leq L
\frac{u^{2}}{\ln(u)}. $$

In both $d=2$ and $d=3$, we have thus obtained that for any $u$
sufficiently large, for any $A \in \mathfrak{A}_{u}^{(M)}$,
\begin{equation*}
\int_{0}^{T} dt \int_{0}^{T-t} ds
\sum_{z' \in \mathbb{Z}^{d}} q_{u^{2}s}(z') H_{1,u}(t,z')
 \leq L' M (\psi_{d}(u))^{-2},
\end{equation*}
where $L'$ is a constant depending on $d$ and $T$.
From (\ref{2moment}), (\ref{HtildeH}) and the above, we deduce (\ref{2mom}). As explained earlier,
this completes the proof of Lemma \ref{variable:u}. $\quad \Box$

\subsection{Proof of Lemma \ref{couronne}}
The proof of Lemma \ref{couronne} is somewhat lengthy.
It is inspired by the first part of \cite{LGP},
where an upper bound for the Hausdorff measure of the support of
two-dimensional super-Brownian motion is established.
In particular, we use the Brownian snake as a main tool. The Brownian snake
gives an alternative construction of super-Brownian motion
under its excursion measure. Moreover, this object introduces
time dynamics in the
analysis of super-Brownian motion which prove to be critical for our arguments to work. We briefly introduce
the Brownian snake and related notation in paragraph 4.4.1, then
discuss the link between Brownian snake and super-Brownian motion.  \PA
For convenience, we work in this section with super-Brownian motion
with branching rate $4$ and diffusion coefficient $1$
under its excursion measure.
Simple scaling arguments then give the general case. \PA
We only give a detailed proof of Lemma \ref{couronne} in the three-dimensional case (paragraph 4.4.4), after having
 summarized the basic idea (paragraph 4.4.2), and presented three intermediate lemmas (paragraph 4.4.3).
 Using the results of the first part of \cite{LGP}, the case $d=2$ easily adapts. In fact, we even establish a stronger
result in the plane (see Lemma \ref{couronned=2} below), which we discuss in paragraph 4.4.5.

\subsubsection{Brownian snake}
For a precise definition of the Brownian snake, we refer to \cite{LG}, Chapter IV.
 Let $\mathcal{W}$ be the set of continuous finite paths from
$\mathbb{R}_{+}$ into $\mathbb{R}^{d}$.
 For $w \in \mathcal{W}$, we denote by
$\zeta_{w}$ the lifetime of $w$, and by $\hat{w}$ the terminal point of the
path $w$, that is $w(\zeta_{w})$. The trivial path $\overline{y}$ is the
path with initial point $y \in \mathbb{R}^{d}$ and lifetime $0$.
The space $\mathcal{W}$ is Polish when equipped with the distance
$$ d(w,w') = |\zeta_{w}-\zeta_{w'}| + \sup_{r \geq 0} | w(r \wedge \zeta_{w})
-w'(r\wedge \zeta_{w'})|.$$
\PA We then consider $\Omega=\mathcal{C}(\mathbb{R}_{+}, \mathcal{W})$,
the space of continuous paths from $\mathbb{R}_{+}$ into $\mathcal{W}$
with the topology of uniform convergence on compact sets,
and $\mathcal{F}=\mathcal{B}(\Omega)$ the Borel $\sigma$-field on $\Omega$.
The canonical
 process on this space is denoted $(W_{s}, s \geq 0)$, and we define
for $s \geq 0$, $\zeta_{s}:= \zeta(W_{s})$, and
$\sigma(\zeta) := \inf\{s > 0 : \zeta_{s}=0 \}$. We also let
$ (\mathcal{F}_{t})_{t \geq 0}$ be the canonical filtration on $\Omega$.
\PA For $w \in \mathcal{W}$,
 we let $\Pi_{w}$ be the law on $(\Omega,\mathcal{F})$ of the Brownian snake
 starting from the path $w$.
Under $\Pi_{w}$, $(W_{s},s\geq 0)$ is a $\mathcal{W}$-valued diffusion and
$(\zeta_{s}, s \geq 0)$ is a one-dimensional reflecting
Brownian motion. Informally, when
$\zeta_{s}$ ``increases'', the path $W_{s}$ grows like a
$d$-dimensional Brownian motion,
whereas it is erased when $\zeta_{s}$ ``decreases'' (see \cite{LG}, Chapter IV for more precisions). \PA
For $y \in \mathbb{R}^{d}$, the
measure $\mathbb{N}_{y}$ is the excursion measure of $W$ away from the
trivial
path $\overline{y}$. We abuse the notation by using the same notation
$\mathbb{N}_{y}$ for the excursion measure of the Brownian snake away from
$\overline{y}$ and for the excursion measure of super-Brownian motion
(cf Section 2.1). This abuse will be justified below when we construct
the excursion measure of super-Brownian motion from the Brownian snake under
$\mathbb{N}_{y}$.
 Under $\mathbb{N}_{y}$,
 the law of $\zeta$ is the It\^o
measure of positive Brownian excursions and $\sigma(\zeta)$ is the length
of this excursion.  \PA
Denote by $\Pi_{w}^{*}$ the law under $\Pi_{w}$ of
$(W_{s\wedge \sigma_{\zeta_{s}}},s \geq 0)$, that is the law of the Brownian
snake stopped when its lifetime process hits $0$.
The
strong Markov property of $W$ under $\mathbb{N}_{y}$
can be expressed in the following way.
Let $\theta _{t}$ denotes the
usual shift operator on $\Omega$. If $T$ is a
$(\mathcal{F}_{t})_{t \geq 0}$-stopping time such that $T>0$
$\mathbb{N}_{y}$-a.e., then, for any nonnegative
$\mathcal{F}_{T}$-measurable $F$, for any nonnegative $\mathcal{F}$-measurable
$G$,
\begin{equation} \label{MarkovW}
\mathbb{N}_{y}\left(\mathbf{1}_{\{T<\infty\}}  F \times G \circ \theta_{T}
\right) = \mathbb{N}_{y} \left( \mathbf{1}_{\{T<\infty\}} F \times
\Pi_{W_{T}}^{*}(G) \right).
\end{equation}
 \PA
The link between Brownian snake and super-Brownian motion can be
expressed as follows.
Let $L_{s}^{t}$ denote
the local time of $\zeta$ at time $s$ and level $t$. Since the law of
$(\zeta_{s}, s \geq 0)$ under $\mathbb{N}_{y}$ is the It\^o measure of positive
Brownian excursions, $(L_{s}^{t}, s \geq 0)$ is, for any $t \geq 0$,
 well-defined, increasing and continuous, $\mathbb{N}_{y}$-a.s. We
denote by $d_{s}L_{s}^{t}$ the measure associated with the
function $u \to L_{u}^{t}$ and we let
$(Y_{t}(W), t \geq 0)$ be the measure-valued process defined by the formula
$$  Y_{t}(W)(.) = \int_{0}^{\sigma(\zeta)}
d_{s}L_{s}^{t} \mathbf{1}_{\{\hat{W}_{s}\in .\}}. $$
 Then, the law of $(Y_{t}(W), t\geq 0)$ under $\mathbb{N}_{y}$ is the
excursion measure of super-Brownian motion with branching rate $4$
and diffusion coefficient $1$\footnote{Moreover, if we let $\mu \in M_{F}(\mathbb{R}^{d})$ and
$\sum_{i \in I} \delta_{y_{i},W_{i}}$ be a Poisson measure
with intensity $\mu(dy)\mathbb{N}_{y}(dW)$, then a
super-Brownian motion $(\overline{Y}_{t}, t\geq 0)$ starting from $\mu$ can
be obtained by setting
$$ \overline{Y}_{t}= \sum_{i \in I} Y_{t}(W_{i}).$$}.

\subsubsection{Outline of the proof of Lemma \ref{couronne}}
 Using a symmetry argument,
we can interchange the roles of $0$ and $x$, and we will thus work under the
probability measure $\mathbb{N}_{x}(. \big| 0 \in \mathcal{R}) =
\mathbb{N}_{x}(.\big| T_{0}< \infty)$, where
 $T_{0} = \inf \{ t \geq 0 : \hat{W}_{t}=0\}$. It is possible to precise the law of $(|W_{T_{0}}|)_{t \leq \zeta_{T_{0}}}$
under $\mathbb{N}_{x}(.|0 \in \mathcal{R})$ (see Lemma \ref{loiWT0} below). \PA
For $j \in \mathbb{N}$, let us introduce $r_{j}= \exp(-j^{2})$. To $n_{1} \in \mathbb{N}$ we associate $\varepsilon_{1}:= r_{2^{n_{1}}}$,
and for $\varepsilon_{0}>0$, we set
$n_{0} := \min\{p \in \mathbb{N}^{*} : r_{2^{p}} \leq \varepsilon_{0}\}$. Note that we have
$$ r_{j} \in (g_{d}(\varepsilon_{0}), \varepsilon_{0}) \ \ \forall j \in \big[\!| 2^{n_{0}}, 2^{n_{0}+1}-1 |\!\big]. $$
\begin{claim} \label{zetaWT0}
One can choose $\alpha>0$ such that, for any $\delta>0$, there exists
$n_{1} \in \mathbb{N}$ such that for any $n_{0} \geq n_{1}$, one has
\begin{equation*}
 \mathbb{N}_{x}
\left( \forall j \in \big[\!| 2^{n_{0}}, 2^{n_{0}+1}-1 |\!\big] :
Y_{\zeta_{T_{0}}}(C(0,r_{j} /2,r_{j})) < \alpha h(r_{j}) \ \big| \
T_{0}<\infty \right) \leq \delta.
\end{equation*}
\end{claim}
From our preceeding remarks, Lemma \ref{couronne} follows from Claim
\ref{zetaWT0} (even if it means changing $\alpha$ to loosen the
inequality). \PA
 The idea of the proof of Claim \ref{zetaWT0} is
the following. For given $w \in \mathcal{W}$, $n_{0} \in \mathbb{N}$
and $j \in \big[\!| 2^{n_{0}}, 2^{n_{0}+1}-1 |\!\big]$, we will
express further the contribution $\mathfrak{Y}_{w}(r_{j})$ to
$Y_{\zeta_{w}}(\mathcal{C}(0,r_{j}/2,r_{j}))$ of particules which
split off the path $w$ in the time interval $[\zeta_{w}- r_{j}^{2}
\ln (1/r_{j}), \zeta_{w}-r_{j}^{2}]$ (see (\ref{Zwepsilon}) below).
We will observe that for large enough $n$, the contributions
$\mathfrak{Y}_{w}(r_{j}), j \in \big[\!|2^{n},2^{n+1} |\!\big]$ are
independent. Using estimates on these contributions (see Lemma
\ref{psiZI} below), this independence will lead us to a bound on the
probability that for any $j \in \big[\!| 2^{n}, 2^{n+1}-1 |\!\big]$,
$\mathfrak{Y}_{w}(r_{j})$ remains smaller than $\alpha h(r_{j})$
(see (\ref{couronne5}) below).  \PA For a well-choosen $\alpha>0$,
we will deduce from this bound and the knowledge of the law of the
path $|W_{T_{0}}|$ the existence of integers $N_{0}, N$, and of a
family of sets of ``good paths'' $(\mathbb{W}_{n}, n \geq N_{0})$
such that, with a probability arbitrarily close to $1$ when $N$ is
large enough,
\begin{itemize}
\item $(|W_{T_{0}}|)_{t \leq \zeta_{T_{0}}}$ belongs to $\mathbb{W}_{n}$ for any $n \geq N$.
\item for any $w \in \mathbb{W}_{n}$, $n \geq N$, there exists $j \in \big[\!| 2^{n}, 2^{n+1}-1 |\!\big]$ such that
$r_{j}^{2}Z_{w}(r_{j}) > \alpha h(r_{j})$.
\end{itemize}
The desired claim will follow (see assertions (\ref{couronne2}), (\ref{couronne3}) and (\ref{couronne4}) below).  \PA
We now present three intermediate lemmas.
\subsubsection{Preliminary results}
Let us start by investigating the law of the path $|W_{T_{0}}|$ under $\mathbb{N}_{x}[.|0 \in \mathcal{R}]$.
\begin{lemma} \label{loiWT0}
Under $\mathbb{N}_{x}[.|0 \in \mathcal{R}]$,
$|W_{T_{0}}(t)|_{t \in[0,\zeta_{T_{0}}]}$ has the law of a Bessel process with
index $-5/2$ started from $|x|$ and stopped when it first hits $0$.
\end{lemma}
 \PA
\underline{Proof of Lemma \ref{loiWT0}}.
 Introduce a $d$-dimensional Brownian
motion $(B_{t})_{t \geq 0}$ and the function
 $u(y):=\mathbb{N}_{y}(T_{0}<\infty) =
\left( 2-d/2 \right)|y|^{-2}$. Let us
denote by $Q$ the law of the solution of the
stochastic differential equation
$$ dY_{t} = dB_{t} + \frac{\nabla u}{u}(Y_{t}), Y_{0}=x, $$
stopped when it hits $0$. \PA
We know from \cite{De}, Proposition 1.4,
 that for any nonnegative continuous function $F$ on
$\mathcal{W}$,
\begin{equation} \label{lawWT0}
\mathbb{N}_{x} \left[\mathbf{1}_{\{T_{0} < \infty\}} F(W_{T_{0}}) \right] =
\int \pi_{0,x}(dW) F(W),
\end{equation}
where $\pi_{0,x} : = \left(2-\frac{d}{2}\right)|x|^{-2} Q.$ \PA
For $t \in [0,\zeta_{T_{0}}]$ let us set $R_{t} := |W_{T_{0}}(t)|.$
From  (\ref{lawWT0}) we deduce that under the probability measure
$\mathbb{N}_{x}[. | T_{0} < \infty]$,
$(R_{t})_{0 \leq t \leq \zeta_{T_{0}}}$
solves the stochastic differential equation
$$ dR_{t} = d\beta_{t} -\frac{2}{R_{t}} dt,$$
where $(\beta_{t})_{t \geq 0}$ is a linear Brownian motion. Thus,
$(|W_{T_{0}}|(t))_{0 \leq t \leq \zeta_{T_{0}}}$
has under $\mathbb{N}_{x}(.\big| \{T_{0}<\infty\})$
the law of a Bessel process with index $-5/2$ started
from $|x|$ and stopped when it first hits $0$. We have completed the proof
of Lemma \ref{loiWT0}. $\quad \Box$ \vskip0.2cm
For $x \in \mathbb{R}^{d}$, $p \in \mathbb{N}^{*}$, $\varepsilon>0$ and $t>0$, we will need a lower bound on
$$\psi(t,x,\varepsilon,p) :=
\mathbb{N}_{x}\left[ Y_{t}\left(\mathcal{C}(0,\varepsilon /2, \varepsilon)  \right)^{p} \right].$$
 Recall $P_{t}$ denotes the semigroup of $d$-dimensional Brownian motion.
We know (see \cite{LGP}, Proposition 3.2) that we have
$$\psi(t,x,\varepsilon,1) = P_{t}\mathbf{1}_{\mathcal{C}(0,\varepsilon/2,\varepsilon)}(x),$$
and the following recursion relation for
$p \geq 2$
\begin{eqnarray} \label{recursion}
\psi(t,x,\varepsilon,p) = 2 \sum_{j=1}^{p-1} {p \choose j} \int_{0}^{t} P_{t-s}\left(
 \psi(s,.,\varepsilon,j) \psi(s,.,\varepsilon,p-j)\right)(x)ds.
\end{eqnarray}
\PA  Fix $c_{1}>1$ and let $c_{2}=1-\frac{1}{2c_{1}}$.
Note that $1/2 < c_{2} <1$. First observe that there exists a positive $c_{3}$ such that
\begin{eqnarray} \label{psi1}
 \psi(t,x,\varepsilon,1) = P_{t}
\mathbf{1}_{\mathcal{C}(0,\varepsilon/2,\varepsilon)}(x)
\geq c_{3} \varepsilon^{3}
\exp \left( -\frac{|x|^{2}}{2t} \right)
t^{-3/2} \mathbf{1}_{\{t \geq \varepsilon^{2}\}}.
\end{eqnarray}
\begin{lemma} \label{psi}
For $d=3$, there exists a positive constant $c_{4}$ so that
 for any $p \in \mathbb{N}^{*}$, $t>0$, $x \in \mathbb{R}^{3}$ and
$\varepsilon \geq 0$,
\begin{equation*}
(\mathcal{H}_{p}) \ \ \ \ \ \ \
 \psi(t,x,\varepsilon,p) \geq c_{4}^{p} p!
\frac{\varepsilon^{2p+1}}{t^{3/2}}
 \exp\left(
-\frac{c_{1}|x|^{2}}{t} \right)
\mathbf{1}_{\{t \geq c_{2}^{-2} \varepsilon^{2}\}}.
\end{equation*}
\end{lemma}
Corollary 3.3 of \cite{LGP} is the corresponding result for the
two-dimensional case. \underline{Proof of Lemma \ref{psi}}. Note
that there exists a constant $c_{5} \geq 1$ such that for any $p
\geq 2$,
$$ \sum_{j=1}^{p-1} (j(p-j))^{3/2} \geq \frac{1}{c_{5}} p^{4}.$$
 Let us set $c_{6} := c_{3}c_{5}^{-1}(\pi c_{1})^{-3/2}$.
 We first verify that for any $p \in \mathbb{N}^{*}$,
\begin{equation*}
(\tilde{\mathcal{H}}_{p}) \ \quad \quad \psi(t,x,\varepsilon,p) \geq
\pi^{3/2} c_{5} c_{6}^{p} p! \varepsilon^{3p} p_{t/2pc_{1}}(x)
\left( \frac{1}{\varepsilon}-\frac{1}{\sqrt{t}}\right)^{p-1}
\mathbf{1}_{\{t \geq \varepsilon^{2} \}}.
\end{equation*}
 We use induction on $p$ to establish $(\tilde{\mathcal{H}}_{p})$. If $p=1$,
using (\ref{psi1}) and our definition of $c_{6}$, we obtain
$$  \psi(t,x,\varepsilon,1) \geq (\pi c_{1})^{3/2} c_{5} c_{6}
\varepsilon^{3} t^{-3/2} \exp \left( -\frac{|x|^{2}}{2t} \right)
\mathbf{1}_{\{t \geq \varepsilon^{2} \}}. $$
Since $c_{1} >1$, $(\tilde{\mathcal{H}}_{1})$ follows.
\PA  Let $p\geq 2$ and assume that the result holds for all $p'<p$.
Using (\ref{recursion}) and the induction assumption we get for
$t \geq \varepsilon^{2}$
\begin{eqnarray*}
\!\!\!\!&\!\!\!&\!\!\!\! \psi(t,x,\varepsilon,p)
\geq 2\pi^{3} c_{5}^{2} c_{6}^{p}  p! \varepsilon^{3p} \\
\!\!\!\!&\!\!\!&\!\!\!\! \qquad \qquad \qquad  \times
\sum_{j=1}^{p-1} \int_{\varepsilon^{2}}^{t} \!ds \!
\int_{\mathbb{R}^{3}} \!\!\!dy \ p_{t-s}(x-y)
p_{s/2jc_{1}}(y) p_{s/2(p-j)c_{1}}(y)
\left(\frac{1}{\varepsilon}-\frac{1}{\sqrt{s}}\right)^{p-2}.
\end{eqnarray*}
For any $j \in \big[\!|1,p-1|\!\big], s>0$ and $y \in \mathbb{R}^{3}$ we have
$$ p_{s/2jc_{1}}(y) p_{s/2(p-j)c_{1}}(y) =
\left(\frac{c_{1}j(p-j)}{\pi p s} \right)^{3/2} p_{s/2pc_{1}}(y) $$
From the last two displays, the choice of $c_{5}$ and the fact that
\\$p_{t-s} * p_{s/2pc_{1}} = p_{t-s+(s/2pc_{1})}$ we obtain
\begin{equation} \label{psi2eme}
 \psi(t,x,\varepsilon,p)
\geq 2 \pi^{3/2} c_{5} c_{6}^{p}  p! \varepsilon^{3p} p^{5/2}
c_{1}^{3/2} \int_{\varepsilon^{2}}^{t} \frac{ds}{s^{3/2}}
p_{t-s+s/2pc_{1}}(x)
\left(\frac{1}{\varepsilon}-\frac{1}{\sqrt{s}}\right)^{p-2}.
\end{equation}
Since $c_{1} >1$ the function $s \to t-s+s/2pc_{1}$ is decreasing,
so that for any $s \in [\varepsilon^{2},t]$, we have
 $t \geq t-s+\frac{s}{2pc_{1}} \geq \frac{t}{2pc_{1}}$. Thus,
$$ (2pc_{1})^{3/2}p_{t-s+\frac{s}{2pc_{1}}}(x) \geq p_{\frac{t}{2pc_{1}}}(x).$$
It follows that
     \begin{eqnarray*}
 \psi(t,x,\varepsilon,p)
\!\!&\!\! \geq \!\!&\!\!
 2^{-1/2} \pi^{3/2} c_{5} c_{6}^{p}  p.p!   \varepsilon^{3p}
  p_{\frac{t}{2pc_{1}}}(x)
\int_{\varepsilon^{2}}^{t} \frac{ds}{s^{3/2}}
\left(\frac{1}{\varepsilon}-\frac{1}{\sqrt{s}}\right)^{p-2} \\
\!\!&\!\! = \!\!&\!\! 2^{1/2} \pi^{3/2} c_{5} c_{6}^{p} p!
\varepsilon^{3p}
  p_{\frac{t}{2pc_{1}}}(x)
 \left(\frac{1}{\varepsilon}-\frac{1}{\sqrt{t}}\right)^{p-1},
\end{eqnarray*}
which finishes the proof of $(\tilde{\mathcal{H}}_{p})$. We have established
$(\tilde{\mathcal{H}}_{p})$ for any $p \in \mathbb{N}^{*}$.
Note in particular that (\ref{psi2eme}) holds for any $p \in \mathbb{N}^{*}$.
 \vskip0.2cm \PA Let us now complete the proof of Lemma \ref{psi}. We set
$c_{4}:=(1-c_{2}^{1/2})c_{6}$. The case $p=1$ follows from
(\ref{psi1}) since $c_{3} \geq c_{4}$. \PA Let us now suppose $p
\geq 2$ and $t \geq c_{2}^{-2} \varepsilon^{2}$. Since $c_{2}=
1-\frac{1}{2c_{1}} \leq \frac{2c_{1}-1}{2c_{1}-1/p}$, we get, for
any $s \leq c_{2}t$, $t-s+s/2pc_{1} \geq t/2c_{1}$, so that
$$ p_{t-s+\frac{s}{2pc_{1}}}(x) \geq (2\pi t)^{-3/2} \exp\left(
-\frac{c_{1}|x|^{2}}{t}\right). $$
Hence, it follows from (\ref{psi2eme}) that for any $p \geq 2$,
\begin{eqnarray*}
\!\! \psi(t,x,\varepsilon,p) \!\!\!& \!\!\! \geq \!\!\! & \!\!\!
2^{-1/2} c_{5} c_{6}^{p} p^{5/2} p! \varepsilon^{3p} t^{-3/2}
\exp\left( -\frac{c_{1}|x|^{2}}{t}\right)  \!\!
\int_{\varepsilon^{2}}^{c_{2}t}\!\! \frac{ds}{s^{3/2}}
\left(\frac{1}{\varepsilon}-\frac{1}{\sqrt{s}}\right)^{p-2}
\\ \!\!\! & \!\!\!\geq \!\!\! & \!\!\!
2^{1/2} c_{5} c_{6}^{p} p^{3/2} p! \varepsilon^{3p} t^{-3/2}
\exp\left( -\frac{c_{1}|x|^{2}}{t}\right)
\left(\frac{1}{\varepsilon} - \frac{1}{\sqrt{c_{2}t}}\right)^{p-1}.
\end{eqnarray*}
Since $t \geq c_{2}^{-2}\varepsilon^{2}$, we have
$\left(\frac{1}{\varepsilon} - \frac{1}{\sqrt{c_{2}t}}\right)^{p-1}
\geq \varepsilon^{-p+1}(1-c_{2})^{p-1}$. Moreover, $c_{5} \geq 1$,
hence,
 $(\mathcal{H}_{p})$ follows for $p \geq 2$, which completes the proof
of Lemma \ref{psi}. $\qquad \Box$ \vskip0.2cm
As explained briefly in paragraph 4.4.2, we will need to estimate, for a fixed $w \in \mathcal{W}$,
 the contribution under $\Pi_{w}$ to $Y_{\zeta_{w}}(\mathcal{C}(0,\varepsilon /2,\varepsilon))$
of particules which split off the path $w$ shortly before $\zeta_{w}$.
 For a given $t>0$, Lemma V.5 of \cite{LG} allows one to decompose $Y_{t}$ under $\Pi_{w}^{*}$ as
the sum of independent contributions corresponding to the decomposition of the path $\zeta$ into its excursions above its minimum-to-date.
Let us state this more precisely. \PA
 Fix $w \in \mathcal{W}$ with $\zeta_{w} >0$. Under $\Pi_{w}^{*}$
we can construct a Poisson point measure
 $\Lambda$ on $[0,\zeta_{w}]\times \Omega$
with intensity $2dt \mathbb{N}_{w(t)}(dW)$ such that
\begin{equation} \label{decompexcurs}
Y_{\zeta_{w}} = \int_{[0,\zeta_{w}]\times \Omega}
Y_{\zeta_{w}-t}(W) \Lambda(dt,dW), \qquad
\Pi_{w}^{*}-\mbox{a.s.}
\end{equation}
Hence, when $w \in \mathcal{W}$ satisfies $\zeta_{w} \geq
\varepsilon^{2}\ln(1/\varepsilon)$, the contribution to
$Y_{\zeta_{w}}(\mathcal{C}(0,\varepsilon /2, \varepsilon))$ of
particules which split off the path $w$ in the time interval
$[\zeta_{w}-\varepsilon^{2}\ln(1/\varepsilon), \zeta_{w} -
\varepsilon^{2}]$ can be written
\begin{equation} \label{Zwepsilon}
\mathfrak{Y}_{w}(\varepsilon) := \int_{\zeta_{w}-
\varepsilon^{2}\ln(1/\varepsilon)}^{\zeta_{w}-\varepsilon^{2}}
\int_{\Omega}
 Y_{\zeta_{w}-t}(\mathcal{C}(0,\varepsilon/2,\varepsilon)) \Lambda(dt,dW).
\end{equation}
For $\varepsilon>0, w \in \mathcal{W}$, we also introduce
$Z_{w}(\varepsilon) := \varepsilon^{-2}
\mathfrak{Y}_{w}(\varepsilon)$ and
 we then estimate
the moments of $Z_{w}(\varepsilon)$ under $\Pi_{w}^{*}$.
\begin{lemma} \label{psiZI}
For $ w \in \mathcal{W}$ and $\varepsilon >0$ such that $\zeta_{w} \geq \varepsilon \ln(1/\varepsilon)$, let us set
$$ I_{w}(\varepsilon) := \varepsilon
\int_{\varepsilon^{2}}^{\varepsilon^{2}\ln\left( 1/\varepsilon \right)}
\exp\left( -c_{1}
\frac{|w(\zeta_{w}-s)|^{2}}{s} \right)
s^{-3/2} \mathbf{1}_{\{s \geq \varepsilon^{2}
c_{2}^{-2}  \}} ds.$$
 There exist positive constants
$c_{7}, c_{8}, c_{9}$ such that for any $w \in \mathcal{W}$,
$\varepsilon_{2} \in (0,1/e)$ such that $\zeta_{w} \geq
\varepsilon_{2}^{2}\ln(1/\varepsilon_{2})$, the following holds.
\begin{description}
\item[(a)] For any $\varepsilon \in (0, \varepsilon_{2})$
and for any $p \in \mathbb{N}$
\begin{equation*}
c_{7}^{p} p^{p} \geq \Pi_{w}^{*}(Z_{w}(\varepsilon)^{p})  \geq
c_{8}^{p}p^{p}
 I_{w}(\varepsilon).
\end{equation*}
\item[(b)] For any $A>0$, let $p= \left\lceil 2A/c_{6}\right\rceil$. Then,
for any $\varepsilon \in (0,\varepsilon_{2})$
\begin{equation*}
 \Pi_{w}^{*}
(Z_{w}(\varepsilon) \geq A) \geq \exp(-c_{9}A) \left(
(2^{p}I_{w}(\varepsilon)-1)^{+}\right)^{2}.
\end{equation*}
\end{description}
\end{lemma}
The lower bounds on $\Pi_{w}^{*}(Z_{w}(\varepsilon)^{p}), p \in \mathbb{N}$
in Lemma \ref{psiZI} (a) are a direct consequence of Lemma \ref{psi}.
Furthermore, the proof of the upper bound in Lemma \ref{psiZI} (a) easily adapts from
the one of Lemma 3.4 in \cite{LGP}. Then, part (b) of Lemma \ref{psiZI}
is deduced from part (a) in the exact same manner as, in \cite{LGP},
Lemma 3.5 is deduced from Lemma 3.4. We leave details to the reader. $ \quad \Box$
\subsubsection{} Let us now complete the proof of Lemma \ref{couronne} by establishing Claim \ref{zetaWT0}.
 We let $\alpha=1/(4c_{9})$ and fix $\delta >0$. Note that $T_{0}$ is a stopping time of the filtration $(\mathcal{F}_{t})_{t \geq 0}$.
Using the strong Markov property (\ref{MarkovW}) at time
$T_{0}$, we have, for $n_{0}>0$,
\begin{eqnarray} \label{couronne2}
&&  \mathbb{N}_{x}\left[
 \forall j \in \big[\!| 2^{n_{0}}, 2^{n_{0}+1}-1 |\!\big]
: Y_{\zeta_{T_{0}}}(C(0,r_{j} /2,r_{j})) < \alpha
h(r_{j}), \ T_{0}<\infty \right]  \\
\!\!\!& \!\!\!\leq \!\!& \!\mathbb{N}_{x}\bigg[
\Pi_{W_{T_{0}}}^{*}\bigg( \forall j \in \big[\!| 2^{n_{0}},
2^{n_{0}+1}-1 |\!\big] : Y_{r}(C(0,r_{j} /2,r_{j})) < \alpha
h(r_{j}) \bigg)_{r=\zeta_{T_{0}}}  \!\!\!\!,   T_{0}<\infty \bigg].
\nonumber
\end{eqnarray}
Notice that (\ref{couronne2}) is an inequality and not an equality,
because we used that
$$ Y_{r} \left( \mathcal{C}(0,r_{j}/2,r_{j})\right) =
\int_{0}^{\sigma(\zeta)} d_{s}L_{s}^{r}
\mathbf{1}_{\mathcal{C}(0,r_{j}/2,r_{j})}(\hat{W}_{s}) \geq
\int_{T_{0}}^{\sigma(\zeta)} d_{s}L_{s}^{r}
\mathbf{1}_{\mathcal{C}(0,r_{j}/2,r_{j})}(\hat{W}_{s}), $$
on the event $\{ T_{0} < \infty \}$. \PA
Introduce the sequence $u_{j}:= r_{2^{j}}^{2} \ln(1/r_{2^{j}})$, and choose $n_{2}$ large enough so that
\begin{equation*}
 r_{2^{n_{2}}} \leq e^{-4}  \ \ \ \ \mbox{ and } \ \ \ \
\mathbb{N}_{x} \left[ \zeta_{T_{0}} \geq s_{n_{2}} \big| T_{0} < \infty \right] \leq \delta/4.
\end{equation*}
Let us set
 $$ \mathfrak{W} := \left\{ w \in \mathcal{W} \ : \
w(0)=x, \hat{w}=0, \zeta_{w} \geq s_{n_{2}} \right\}.$$
Our choice of $n_{2}$ ensures that
\begin{equation} \label{conditionn2}
\mathbb{N}_{x}[W_{T_{0}} \in \mathfrak{W}|T_{0}<\infty] \geq 1-\delta /4.
\end{equation}
Since $c_{2}>1/2$, it also
guarantees that for any $w \in \mathfrak{W}$, $j \geq 2^{n_{2}}$, $2r_{j}^{2}c_{2}^{-2} \leq \zeta_{w}$.
 We can then define, for $B>0$, $j \geq 2^{n_{2}}$,
$$ W_{B,j} := \left\{ w \in \mathfrak{W} \ : \
\sup_{s \leq 2r_{j}^{2}c_{2}^{-2}} |w(\zeta(w)-s)| > B
r_{j}c_{2}^{-1} \right\}. $$
For $w \in \mathfrak{W}$, $n \geq n_{2}$, we then introduce
$$ F_{n,B}(w) := 2^{-n}\sum_{j=2^{n}}^{2^{n+1}-1}
\mathbf{1}_{\left\{ w \in W_{B,j} \right\}},$$
 and for $n \geq n_{2}$, we finally let
$$ \mathbb{W}_{n} := \left\{ w \in \mathfrak{W} \ : \ \forall
p \geq n \ \ F_{p,B}(w) < 1/2 \right\}. $$
From (\ref{conditionn2}), it follows that
\begin{eqnarray*}
 \mathbb{N}_{x}\left[ W_{T_{0}} \notin \mathbb{W}_{n} \bigg| T_{0} <
\infty \right]  \leq  \frac{\delta}{4} + \!\sum_{p \geq n} \mathbb{N}_{x} \bigg[
W_{T_{0}}
\in \mathfrak{W},
F_{p,B}(W_{T_{0}}) \geq 1/2 \bigg| T_{0} < \infty \bigg].
\end{eqnarray*}
Using Lemma \ref{loiWT0} and following the arguments of the proof of Lemma 1 in \cite{LG2},
 one can easily establish that there exist constants $B>0, C>0$ such that
$$ \mathbb{N}_{x} \left[
W_{T_{0}} \in \mathfrak{W}, F_{p,B}(W_{T_{0}}) \geq 1/2 | T_{0} <
\infty \right] \leq C e^{-n}. $$ Hence, there exists $n_{3} \geq
n_{2}$ large enough so that for any $n \geq n_{3}$,
\begin{equation} \label{FnBWT0}
 \mathbb{N}_{x} \left[ W_{T_{0}} \notin \mathbb{W}_{n} \bigg| T_{0} < \infty \right] \leq \delta /2,
\end{equation}
which yields
\begin{eqnarray} \label{couronne3}
&& \!\!\! \mathbb{N}_{x}\left[ \Pi_{W_{T_{0}}}^{*}
\bigg( \forall j \in \big[\!| 2^{n}, 2^{n+1}-1 |\!\big]
: Y_{r}(C(0,r_{j} /2,r_{j})) < \alpha
h(r_{j}) \bigg)_{r=\zeta_{T_{0}}} \ \bigg|
\  T_{0}<\infty \right] \nonumber
\\ && \leq  \delta/2  +
\sup_{w \in \mathbb{W}_{n}}  \bigg\{ \Pi_{w}^{*} \big[ \forall j \in
\big[\!| 2^{n}, 2^{n+1}-1 |\!\big] : Y_{\zeta_{w}}(C(0,r_{j}
/2,r_{j})) < \alpha h(r_{j})  \big] \bigg\}.
\end{eqnarray}
Let $n \geq n_{3}$ and $w \in \mathbb{W}_{n}$. Since $r_{j}^{2} \geq r_{j+1}^{2} \ln (1/r_{j+1})$ for $2^{n} \leq j \leq 2^{n+1}-1$,
 the independence properties
of Poisson measures imply that for any $n \in \mathbb{N}$, the variables
$Z_{w}(r_{j}), 2^{n} \leq j \leq 2^{n+1}-1$
are independent under $\Pi_{w}^{*}$.
 Using (\ref{decompexcurs}) and the definition of $Z_{w}(\varepsilon)$, we then get
 \begin{eqnarray} \label{couronne4}
&& \Pi_{w}^{*} \left[\forall j \in \big[\!| 2^{n}, 2^{n+1}-1
|\!\big] : Y_{\zeta_{T_{0}}}(C(0,r_{j} /2,r_{j})) < \alpha h(r_{j})
 \right]  \\
 & \leq &  \Pi_{w}^{*} \bigg[
Z_{w}(r_{j}) < \alpha \theta(r_{j})  \ \forall j \in
\big[\!|2^{n},2^{n+1}-1
 |\!\big] \bigg] = \! \prod_{j=2^{n_{\mathrm{o}}}}^{2^{n_{0}+1}-1}
\Pi_{w}^{*} \bigg( Z_{w}(r_{j}) < \alpha \theta(r_{j})\bigg),
\nonumber
\end{eqnarray}
where we set, for $r>0$,
$\theta(r) := r^{-2} h(r)$.
Furthermore, Lemma \ref{psiZI} (b) leads to
\begin{eqnarray} \label{couronne5}
\!\!\!  \prod_{j=2^{n_{0}}}^{2^{n_{0}+1}-1} \Pi_{w}^{*} \bigg(
Z_{w}(r_{j}) < \alpha \theta(r_{j})\bigg) \!\!
 & \!\!\!\!\leq \!\!\!\!&  \!\!\prod_{j=2^{n_{0}}}^{2^{n_{0}+1}-1}
\left(1-e^{-c_{9}\alpha \theta(r_{j})} \left( \left(
2^{p_{j}}I_{w}(r_{j})-1 \right)^{+} \right)^{2}\right)
\nonumber  \\
&\!\!\leq \!\!&\exp \left\{ - \sum_{j=2^{n}}^{2^{n+1}-1}
j^{-2c_{9}\alpha} \left( \left( 2^{p_{j}}I_{w}(r_{j})-1 \right)^{+}
 \right)^{2} \right\}.
\end{eqnarray}
 We then note that, if $j\geq 2^{n_{2}}$ and
$w \in \mathfrak{W} \setminus W_{B,j}$,
 an easy computation provides
\begin{equation*}
I_{w}(r_{j}) \geq c_{2} \exp\left(-c_{1}B^{2}\right) =: K(B).
 \end{equation*}
Hence, from our choice of $\alpha$, there exists $n_{4} \geq n_{3}$ so that for any $w \in \mathbb{W}_{n}$, $n \geq n_{4}$, one has
$$ \sum_{j=2^{n}}^{2^{n+1}-1} j^{-2c_{9}\alpha}
 \left( \left( 2^{p_{j}}I_{w}(r_{j})-1 \right)^{+}
 \right)^{2} \geq 2^{n+1/2} .$$
 We finally choose
 $n_{0} > n_{1} \geq  n_{4} \geq n_{3} \geq n_{2}$ large enough so that $\exp(-2^{n+1/2}) \leq \delta /2$,
and combine (\ref{couronne2}), (\ref{couronne3}), (\ref{couronne4}) and (\ref{couronne5}) with the above inequality
 to obtain Claim \ref{zetaWT0}. As explained in paragraph 4.4.2, Lemma \ref{couronne} follows.$ \quad \Box$ \\

\subsubsection{The case $d=2$}
We know from \cite{Pe},
Section III.3 that for any $t>0$,
$h$ is for $d=3$ the correct Hausdorff measure function
of $\mathcal{R}_{t}$. On the other hand,
when $d=2$, the correct Hausdorff measure function
of $\mathcal{R}_{t}$, is, as it is proven in
\cite{LGP}, the function
$$ h_{2}(\varepsilon) = \varepsilon^{2} \ln(\varepsilon^{-1})
\ln(\ln(\ln(\varepsilon^{-1}))). $$
Not surprinsingly, when $d=2$, one can in fact establish a stronger result than Lemma
\ref{couronne}.
\begin{lemma} \label{couronned=2}
 We can choose
$\alpha > 0$ so that, for any $\delta>0$, there exists
$\varepsilon_{1}\in \left(0,1 \wedge \frac{|x|}{2}\right)$
such that for any $\varepsilon_{0} \in (0,\varepsilon_{1})$,
\begin{equation*}
 \mathbb{N}_{0} \left[
\exists s\geq 0 \ \exists \ \varepsilon \in
(g_{d}(\varepsilon_{0}),\varepsilon_{0}) : Y_{s}\left(
\mathcal{C}\left(x,\frac{\varepsilon}{2}, \varepsilon \right)\right) > \alpha
h_{2}(\varepsilon) \ \bigg|\ x \in \mathcal{R} \right] \geq 1-\delta.
\end{equation*}
\end{lemma}
Lemma \ref{couronned=2} clearly implies the two-dimensional case of Lemma \ref{couronne}.
The proof is similar to that of Lemma \ref{couronne} in the
three-dimensional case. Let us only point out the main differences,
and leave details to the reader. \PA
Obviously, one should work with $h_{2}$ instead of $h$, $g_{2}$ instead of $g$ and the function $\theta_{2}$ such that
$\theta_{2}(r) := \ln \ln \ln (1/r)$ instead of $\theta$.
Moreover, the sequence $r_{j}$ is to be replaced with $r_{j}^{(2)}=2^{-2^{j}}$, so that $r_{2^{n+1}-1}^{(2)} \geq g_{2}(r_{2^{n}}^{(2)})$.
We already noted that Lemma \ref{psi} for the three-dimensional case corresponds to Corollary 3.3 of \cite{LGP} in the plane.
In particular, note that $c_{1},c_{2},c_{4}$ should be replaced with
$c_{1}^{(2)}>1/2, c_{2}=(4c_{1}^{(2)}-2)/(4c_{1}^{(2)}-1),
c_{4}^{(2)}=c_{3}^{(2)}(2c_{1}^{(2)})^{-1}$.
We also already remarked that Lemma \ref{psiZI} (a) and (b) are to be respectively related with Lemma 3.4, respectively Lemma 3.5 of \cite{LGP}.
 Lemma \ref{psiZI} remains valid in the plane when one replaces $Z_{w}(\varepsilon)$ with
$$ Z_{w}^{(2)}(\varepsilon) :=  \big(\varepsilon^{2}\ln(\varepsilon^{-1}))^{-1}
\int_{\zeta_{w}-\varepsilon^{2}\ln(1/\varepsilon)}^{\zeta_{w}-\varepsilon^{2}}
 Y_{\zeta_{w}-t}(\mathcal{C}(0,\varepsilon/2,\varepsilon)) \Lambda(dt,dW),$$
 then $I_{w}(\varepsilon)$ with
\begin{eqnarray*}
 I_{w}^{(2)}(\varepsilon,p) \!\!\!\!\!&\!\!\!\!\!\!&\!\!\!\!\!:=
p \left(\log(\varepsilon^{-1})\right)^{-p} \\ &&\  \times
\int_{\varepsilon^{2}}^{\varepsilon^{2}\ln\left( 1/\varepsilon \right)}
\exp\left( -c_{1}^{(2)}
\frac{|w(\zeta_{w}-s)|^{2}}{s} \right) \left(\log^{+}\left(
s(c_{2}^{(2)})^{p}\varepsilon^{-2} \right)\right)^{p-1} s^{-1}  ds.
\end{eqnarray*}
Is is then straightforward to check that, once all these changes have been made, the exact same proof as in paragraph 4.4.4 leads to
 assertions similar to  (\ref{couronne2}),
 (\ref{conditionn2}), (\ref{FnBWT0}),
(\ref{couronne3}), (\ref{couronne4}), (\ref{couronne5}) and
 Claim \ref{zetaWT0}.

\section{Results on coalescing random walks}
In this section, we prove\footnote{Some of the ideas involved in the
results below are borrowed from \cite{LG'}, such as in particular,
the case $d \geq 3$ of Lemma
 \ref{leminter}. The proof of Lemma \ref{boundq}
 is also borrowed from this unpublished manuscript.
 Moreover, it is interesting to note that it would be
 possible to get precise asymptotics on
the quantities we are bounding below. In particular, for $\alpha
>0$ and $\phi \in C_{b}(\mathbb{R}^{d})$, it is possible to
get the exact asymptotics of quantities such as
$\mathbb{E}_{0,y}\left[T_{1} \leq |y|^{\alpha}t,
\phi(X_{T_{1}^{1}}/|y|\right]$, as $|y| \to \infty$. Such
asymptotics were computed in \cite{LG'} in the case $d \geq 3$,
$\alpha=2$, and it is possible to extend the results to a general
$\alpha>0$ and to the case $d=2$. For $d \geq 4$, these exact
estimates would allow one to get a more precise upper bound on
$E[(U_{T})^{2}]$ than the one obtained in Section 4.1, and therefore
improve the constants $a_{d}, d\geq 4$ appearing in the statement of
Proposition \ref{lowerbd} (however, this would not be enough to get
precise asymptotics on the hitting probability of a far point for $d
\geq 5$). We chose not to present
 these asymptotics here, as they only had this minor impact on our
main result.} Lemmas \ref{boundq} and \ref{boundcoal},
 which we used in the proof of Theorem \ref{main}.
First, let
us introduce further notation. We write $\mathbb{P}_{x}^{(2)}$ for a
probability measure under which $(Z_{t},t \geq 0)$ is a continuous
time random walk with rate $2$ (instead of rate $1$ for
$\mathbb{P}_{x}$) and jump kernel $p$, starting from $x$. Let us
also denote
\begin{eqnarray*}
 H_{t}(x) : = \mathbb{P}_{x}(Z \mbox{ hits } 0 \mbox{ before } t); \qquad
 G_{t}(x) := \mathbb{E}_{x}
\left[ \int_{0}^{t} \mathbf{1}_{\{ Z_{s}=0 \}}ds \right];
\end{eqnarray*}
and $H_{t}^{(2)}(x), G_{t}^{(2)}(x)$ the corresponding quantities under $\mathbb{P}_{x}^{(2)}$.
It is well-known that when $d=2$,
\begin{equation} \label{asgreend=2}
G_{c^{2}t}^{(2)}(0)
\underset{c \to \infty}{\sim} G_{c^{2}\varepsilon}^{(2)}(0)
\underset{c \to \infty}{\sim} \frac{1}{2\pi} \ln(c).
\end{equation}
When $d=3$, from the definition of $k_{d}$, we have
\begin{equation} \label{asgreend=3}
G_{c^{2}t}^{(2)}(0) \underset{c \to \infty}{\longrightarrow}
k_{d}^{-1}.
\end{equation}
In dimension $d=2$,
an easy adaptation of \cite{La}, Theorem 1.6.1, to the continuous time setting,
ensures the existence of
\begin{eqnarray*}
a(x) := \lim_{t\to \infty} [G_{t}(0)-G_{t}(x)], \ \ \ a^{(2)}(x) := \lim_{t\to \infty} [G_{t}^{(2)}(0)-G_{t}^{(2)}(x)].
\end{eqnarray*}
An easy consequence of the proof of Theorem 1.6.2 in \cite{La} is
the existence of a constant $\kappa_{5}$ depending only on $d$ such
that for any $t \geq 1$,
\begin{equation} \label{a(x)}
G_{t}(0)-G_{t}(x) \leq \kappa_{5}a(x), \ \ \
G_{t}^{(2)}(0)-G_{t}^{(2)}(x) \leq \kappa_{5}a^{(2)}(x).
\end{equation}
Furthermore, from
Theorem 1.6.2 of \cite{La}, both $a(x)$ and $a^{(2)}(x)$ are $O(\ln|x|)$.
\subsection{Proof of Lemma \ref{boundq}}
This proof was taken from \cite{LG'}. In the following, we denote by
$C_{i}, i \geq 0$  positive constants depending only on $d$.
 Let us denote by $(S_{n}, n \in \mathbb{N})$
a discrete time random walk with jump kernel $p$, starting from $x$
under the probability measure $\mathbb{Q}_{x}$. By combining the
well-known bound $\mathbb{Q}_{x}[S_{n}=y] \leq K n^{-d/2}$
 and the martingale inequality of Ledoux and Talagrand (\cite{LT}, Lemma 1.5),
we get for any $n \geq 1, y \in \mathbb{Z}^{d}$,
$$ \mathbb{Q}_{x}[S_{n}=y] \leq \frac{C_{1}}{n^{d/2}}
\exp\left( -\frac{C_{2}|y-x|^{2}}{n}\right).  $$
Let $(N_{t},t \geq 0)$ be a standard Poisson process. Then
\begin{eqnarray} \label{ledTal}
\mathbb{Q}_{0}[Z_{t}=y] &=& \sum_{n=0}^{\infty} P[N_{t}=n]
\mathbb{Q}_{0}[S_{n}=y] \nonumber \\
& \leq & \exp(-t) \mathbf{1}_{\{y=0\}} + \sum_{n=1}^{\infty} P[N_{t}=n]
\frac{C_{1}}{n^{d/2}}
\exp\left( -\frac{C_{2}|y|^{2}}{n}\right).
\end{eqnarray}
We also have for any $t \geq 0$,
$$ \sum_{n=1}^{\infty} \frac{1}{n^{d/2}} P[N_{t}=n] \leq C_{3}t^{-d/2}. $$
It follows from the above that
\begin{equation} \label{2t,2t}
 \sum_{n=1}^{\lfloor 2t \rfloor} P[N_{t}=n]
\frac{C_{1}}{n^{d/2}}
\exp\left( -\frac{C_{2}|y|^{2}}{n}\right) \leq
\frac{C_{1}C_{3}}{t^{d/2}}
\exp\left( -\frac{C_{2}|y|^{2}}{2t}\right) .
\end{equation}
For values of $n$ greater then $2t$, a simple large deviation estimate gives
for every $t>0, m \geq 1$,
$$P[N_{t} \geq 2^{m}t] \leq C_{4}\exp(-C_{5}2^{m}t). $$
    Hence,
\begin{eqnarray*}
&& \sum_{n>2t} P[N_{t}=n] n^{-d/2} \exp\left( -\frac{C_{2}|y|^{2}}{n}\right) \\
\qquad
& \leq & \sum_{m=1}^{\infty} C_{4} \exp(-C_{5}2^{m}t) (2^{m}t)^{-d/2} \exp
\left( -\frac{C_{2}|y|^{2}}{2^{m+1}t}\right) \\
\qquad
& = & C_{4} t^{-d/2}\sum_{m=1}^{\infty} 2^{-md/2} \exp\left(-C_{5}2^{m}t -
\frac{C_{2}|y|^{2}}{2^{m+1}t} \right).
\end{eqnarray*}
Setting $C_{6} = \sqrt{2C_{2}C_{5}}$ we now obtain from the above
$$\sum_{n>2t} P[N_{t}=n] n^{-d/2} \exp\left( -\frac{C_{2}|y|^{2}}{n}\right)
\leq  C_{7} t^{-d/2} \exp(-C_{6}|y|). $$
 Combining (\ref{ledTal}), (\ref{2t,2t}) and the above now gives
$$ q_{t}(y) \leq \exp(-t)\mathbf{1}_{\{y=0\}} + \frac{C_{1}C_{3}}{t^{d/2}}
\exp\left( -\frac{C_{2}|y|^{2}}{2t}\right) +
\frac{C_{7}}{t^{d/2}}\exp(-C_{6}|y|),$$
which clearly implies Lemma \ref{boundq}. $\qquad \Box$

\subsection{Proof of Lemma \ref{boundcoal}}
In the following, we use $K,K'$ to denote positive constants
depending only on $d,T$ and which may change from line to line.
 We consider only the case when $c$ and $x$ are such that $cx \in \mathbb{Z}^{d}$.
The general case immediately follows.
\PA We are first going to rule out small values of $t$. We deal with
 the integral over the interval $[0,c^{-2}]$. Note that
$$ \mathbb{P}_{0,y} \left[ T_{1} \leq c^{2}t, Z_{c^{2}t}^{1}=cx \right] \leq \min \left\{ q_{c^{2}t}(0,cx), q_{c^{2}t}(y,cx) \right\}. $$
 Considering separately the cases $|y|\leq 2c|x|$ and $|y| >2c|x|$ and
using (\ref{exp}), we easily obtain
\begin{eqnarray*}
&&  c^{d} \psi_{d}(|y|) \int_{0}^{c^{-2}} dt
\mathbb{P}_{0,y} \left[ T_{1}
 \leq c^{2}t, Z_{c^{2}t}^{1}=cx \right]
\\ & \leq & \kappa_{3} c^{d-2} \exp(-\kappa_{4}c|x|/2) \times \psi_{d}(|y|)
\exp(-\kappa_{4}(|y|/4)) \\
& \leq & K |x|^{2-d} \leq K \psi_{d}(|x|^{-1}).
\end{eqnarray*}
\PA Let us now deal with small values of $T_{1}$. In a similar way
as in the previous computation, one gets
\begin{eqnarray*}
&&  c^{d} \psi_{d}(|y|) \int_{0}^{c^{-2}} dt
\mathbb{P}_{0,y} \left[ T_{1}
 \leq 1, Z_{c^{2}t}^{1}=cx \right]
\\ & \leq & \kappa_{3} c^{d} \exp(-\kappa_{4}c|x|/2) \times \psi_{d}(|y|)
\exp(-\kappa_{4}(|y|/4\sqrt{2})) \\
& \leq & K|x|^{2-d} \leq K \psi_{d}(|x|^{-1}).
\end{eqnarray*}
Note that to obtain the last line above, we used the assumption $c \geq |x|^{-2}$. \PA
Let us then deal with large values of $T_{1}$. For $t \geq c^{-2}$, we have
\begin{eqnarray*}
 && c^{d} \psi_{d}(|y|) \mathbb{P}_{0,y} \left[
c^{2}t-1  \leq T_{1} \leq c^{2}t, Z_{c^{2}t}^{1}=cx \right] \\
& \leq & e c^{d} \psi_{d}(|y|) \mathbb{P}_{0,y}\left[c^{2}t-1  \leq T_{1} \leq c^{2}t,
Z_{c^{2}t}^{1}=Z_{T_{1}}^{1}=cx \right]
\\ & \leq & e c^{d}\psi_{d}(|y|)
\mathbb{P}_{0,y} \left[ Z_{c^{2}t}^{1}= Z_{c^{2}t}^{2}=cx \right]
\\ & \leq & e c^{-d}\psi_{d}(|y|) f_{t}(x) f_{t}(x-y/c),
\end{eqnarray*}
where we used Lemma \ref{boundq} at the last line above. Hence, from (\ref{integf}), we get that
\begin{eqnarray*}
 && c^{d} \psi_{d}(|y|) \int_{c^{-2}}^{T} dt \mathbb{P}_{0,y} \left[
c^{2}t-1 \leq T_{1}  \leq c^{2}t, Z_{c^{2}t}^{1}=cx \right]\\
 & \leq &  K c^{-d}
\psi_{d}(|y|)
\left( |x|^{2-2d} \wedge \left(\frac{|y|}{c}\right)^{2-2d} \right)
\leq  K |x|^{2-d} \leq K \psi_{d}(|x|^{-1}).
\end{eqnarray*}
Note that, at the last line above, we used the assumption $c \geq |x|^{-2}$.
\PA
We can now suppose
 $1 \leq c^{2}t -1$, that is $t \geq 2c^{-2}$, and restrict our attention to estimating
$$ \int_{2c^{-2}}^{T} dt \mathbb{P}_{0,y} \left[
1 \leq T_{1}  \leq c^{2}t-1 , Z_{c^{2}t}^{1}=cx \right].$$
Using the Markov property at time $T_{1}$, then Lemma \ref{boundq}, we obtain
 \begin{eqnarray} \label{start}
\!\!\!\!&&\!\!\!\! \ \ \ \ \ \
 \int_{2c^{-2}}^{T} dt \ \mathbb{P}_{0,y} \left[
1 \leq T_{1}  \leq c^{2}t-1 , Z_{c^{2}t}^{1}=cx \right]  \nonumber \\ && =
  \int_{2c^{-2}}^{T} dt \ \mathbb{E}_{0,y} \bigg[
\mathbf{1}_{\{1 \leq T_{1}  \leq c^{2}t-1 \}} \sum_{z\in \mathbb{Z}^{d}}
\mathbf{1}_{\{ Z_{T_{1}}^{1}=z \}}q_{c^{2}t-T_{1}}(z-cx) \bigg] \nonumber \\
&& \leq   Kc^{-d}
 \int_{2c^{-2}}^{T} dt \sum_{z \in \mathbb{Z}^{d}}
 \mathbb{E}_{0,y} \left[ \mathbf{1}_{\{1 \leq T_{1}  \leq c^{2}t-1,
Z_{T_{1}}^{1}=z \}}
f_{t-\frac{T_{1}}{c^{2}}}\left(\frac{z}{c}-x \right) \right].
\end{eqnarray}
In order to bound the above quantity, we need the following intermediate result.
\begin{lemma} \label{leminter}
  Let $d \geq 2$.
There exists a positive constant $L_{5}$ depending only on $d$ such
that for any $c \geq 1$, $z \in \mathbb{Z}^{d}$, $t \geq 2c^{-2}$,
$y \in \mathbb{Z}^{d} \setminus 0$ and every measurable function
$\phi : \mathbb{R}_{+}\times \mathbb{Z}^{d} \to \mathbb{R}_{+}$,
\begin{eqnarray*}
&& |y|^{d} \psi_{d}(|y|) \mathbb{E}_{0,y}
\left[ \mathbf{1}_{\{1 \leq T_{1} \leq c^{2}t-1,Z_{T_{1}}^{1}=z \}}
\phi\left(\frac{T_{1}}{|y|^{2}}, Z_{T_{1}}^{1} \right) \right] \nonumber \\
&& \leq L_{5}
\int_{|y|^{-2}}^{\frac{tc^{2}}{|y|^{2}}-\frac{1}{2|y|^{2}}}
 du \ \Phi \left( u,z \right)
 \tilde{f}_{u}\left(\frac{z}{|y|}\right) \tilde{f}_{u}\left(\frac{z-y}{|y|}
\right),
\end{eqnarray*}
where $\Phi(u,z) = \sup_{\left(u-\frac{1}{2|y|^{2}}\right)^{+} \leq r \leq u} \phi(r,z)$,
and $\tilde{f}_{u}$ was defined in Section 2.4.
\end{lemma}
\underline{Proof of lemma \ref{leminter}}.
In this proof, we use $L$ to denote a constant depending only on $d$ and
which may change from line to line. Obviously,
$$  \phi( r,z)
 \leq 2 |y|^{2} \int_{r}^{r+\frac{1}{2|y|^{2}}} \Phi(u,z)du.$$
It follows that
 \begin{eqnarray} \label{interme}
 && |y|^{d} \psi_{d}(|y|)\mathbb{E}_{0,y}
\left[ \mathbf{1}_{\{1 \leq T_{1} \leq c^{2}t-1,Z_{T_{1}}^{1}=z \}}
\phi \left(\frac{T_{1}}{|y|^{2}},Z_{T_{1}}^{1}\right) \right] \nonumber \\
&& \leq 2 |y|^{d+2} \psi_{d}(|y|) \mathbb{E}_{0,y}
\bigg[ \mathbf{1}_{\{ 1 \leq T_{1} \leq c^{2}t-1, Z_{T_{1}}^{1}=z \}}
\nonumber \\
&& \qquad \qquad \qquad \qquad \qquad \quad
\times \int_{|y|^{-2}}^{\frac{tc^{2}}{|y|^{2}}-\frac{1}{2|y|^{2}}}\!\! du
\Phi(u,z) \mathbf{1}_{\{ \frac{T_{1}}{|y|^{2}} \leq u \leq
\frac{T_{1}}{|y|^{2}} + \frac{1}{2|y|^{2}}\}}
\bigg] \nonumber \\
&& = 2|y|^{d+2} \psi_{d}(|y|)
\int_{|y|^{-2}}^{\frac{tc^{2}}{|y|^{2}}-\frac{1}{2|y|^{2}}} \!\!\!\!du
 \Phi(u,z) \\
&& \qquad \qquad \qquad \qquad \qquad \times \mathbb{P}_{0,y}
\left[  1 \leq T_{1} \leq u|y|^{2} \leq T_{1} +\frac{1}{2} \leq c^{2}t
-\frac{1}{2}, Z_{T_{1}}^{1}=z \right] \nonumber .
\end{eqnarray}
where we use the Fubini theorem at the last line. Hence, proving
Lemma \ref{leminter} reduces to establishing the following claim
\begin{claim} \label{eqinter}
If $u>|y|^{-2}$,
\begin{eqnarray*}
&&  2|y|^{d+2} \psi_{d}(|y|)\mathbb{P}_{0,y} \left[  1 \leq T_{1}
\leq u|y|^{2} \leq T_{1} +\frac{1}{2} \leq c^{2}t -\frac{1}{2},
Z_{T_{1}}^{1}=z \right] \nonumber \\ && \leq L_{5}
\tilde{f}_{u}\left(\frac{z}{|y|}\right)
\tilde{f}_{u}\left(\frac{z-y}{|y|}\right).
\end{eqnarray*}
\end{claim}
Let us first rule out the easy cases of Claim \ref{eqinter}. \PA
First note that the case $d \geq 3$ is simple, because
$$\mathbb{P}_{0,y} \left[
 T_{1} \leq |y|^{2}u \leq T_{1}+\frac{1}{2} ,
Z_{T_{1}}^{1}=z \right] \leq e q_{|y|^{2}u}(z)q_{|y|^{2}u}(y-z),$$
and we can use Lemma 1 to conclude. \PA In the case $d=2$, using the
same argument as in the case $d \geq 3$ only gives
\begin{eqnarray*}
&&    |y|^{4} \ln(|y|\vee e) \mathbb{P}_{0,y} \left[
1 \leq T_{1} \leq |y|^{2}u \leq T_{1}+\frac{1}{2} ,
Z_{T_{1}}^{1}=z   \right] \\ &&  \leq e \ln(|y|\vee e)
 f_{u}\left(\frac{z}{|y|}\right)
f_{u}\left(\frac{z-y}{|y|} \right).
\end{eqnarray*}
However, in the particular cases when $|y| \leq A$ for some fixed
constant $A \geq 1$, or when $ |y|^{-2} \leq u \leq |y|^{-1}$, we have
$$\ln(|y| \vee e)
 \exp \left( -\frac{\kappa_{2}|z|}{2 |y|\sqrt{u}}\right)
\exp \left( -\frac{\kappa_{2}|z-y|}{2 |y|\sqrt{u}}\right) \leq L. $$
This easily leads to the desired claim in these particular cases.
\vskip0.2cm We now suppose $d=2$, $|y| \geq A := 6^{8}$, and $u \geq
|y|^{-1}$ and outline of the proof of Claim \ref{eqinter}. We have
\begin{eqnarray} \label{sep}
&& \mathbb{P}_{0,y} \left[
1 \leq T_{1} \leq |y|^{2}u \leq T_{1}+\frac{1}{2} ,
Z_{T_{1}}^{1}=z   \right] \nonumber \\ &&\leq  e \mathbb{P}_{0,y} \left[
1 \leq T_{1} \leq |y|^{2}u \leq T_{1}+\frac{1}{2} ,
Z_{|y|^{2}u}^{1}=Z_{|y|^{2}u}^{2}=z   \right] \nonumber \\
&& \leq  e \mathbb{P}_{z,z} \left[
Z_{s}^{1} \neq Z_{s}^{2} \ \forall s \in [1,|y|^{2}u] ,
Z_{|y|^{2}u}^{1}=0, Z_{|y|^{2}u}^{2}=y   \right],
\end{eqnarray}
 where we used a time-reversal argument in the last line. \PA
We are going to use (\ref{sep}) and argue under $\mathbb{P}_{z,z}$.
On the one hand, with high probability, both $Z_{|y|^{3/2}(u\wedge
1)}^{1}$ and $Z_{|y|^{3/2}(u\wedge 1)}^{2}$ should remain close to
$z$. More precisely, if we set $$B_{z,y} := B
\left(z,|y|^{7/8}\left(\frac{|z|}{3|y|}\vee 1 \right)\right), \qquad
t(u,y) := |y|^{2}u -|y|^{3/2}(u\wedge 1), $$ we will establish that
\begin{eqnarray} \label{bouley}
&& \ln(|y|)|y|^{4}\mathbb{P}_{z,z} \!\bigg[ Z_{|y|^{3/2}(u\wedge
1)}^{1} \notin B_{z,y} \mbox{ or } Z_{|y|^{3/2}(u\wedge 1)}^{2}
\notin B_{z,y},  Z_{|y|^{2}u}^{1}=0, Z_{|y|^{2}u}^{2}=y  \bigg]
\nonumber \\ && \quad \leq L \tilde{f}_{u}\left(\frac{z}{|y|}\right)
\tilde{f}_{u}\left(\frac{z-y}{|y|}\right).
\end{eqnarray}
 On the other hand, when both $Z_{|y|^{3/2}(u\wedge 1)}^{1}, Z_{|y|^{3/2}(u\wedge 1)}^{2}$
 are close to $z$, we obtain from the Markov property for the walks $Z^{1}, Z^{2}$ at time
$|y|^{3/2}(u \wedge 1)$ that
\begin{eqnarray}  \label{markovy3/2}
\!\!\!\!&&\!\!\!\! \mathbb{P}_{z,z} \bigg[ Z_{s}^{1} \neq Z_{s}^{2}
\ \forall s \in [1,|y|^{3/2}(u \wedge 1)] , \
Z_{|y|^{2}u}^{1}=0, \ Z_{|y|^{2}u}^{2}=y, \nonumber \\
&& \qquad \qquad \qquad \qquad \qquad
 Z_{|y|^{3/2}(u\wedge 1)}^{1} \in B_{z,y}, \ Z_{|y|^{3/2}(u\wedge 1)}^{2} \in
B_{z,y} \bigg]  \nonumber \\ && \leq  \mathbb{P}_{z,z} \bigg[
Z_{s}^{1} \neq Z_{s}^{2} \ \forall s \in [1,|y|^{3/2}(u \wedge 1)]
\bigg] \sup_{(x_{1},x_{2}) \in B_{z,y}^{2}} q_{t(u,y)}(x_{1})
q_{t(u,y)}(x_{2}-y).
\end{eqnarray}
We will then establish, using Lemma \ref{boundq}, that
\begin{equation} \label{supf}
 \sup_{(x_{1},x_{2}) \in B_{z,y}^{2}}
|y|^{4}q_{t(u,y)}(x_{1}) q_{t(u,y)}(x_{2}-y)   \leq
\tilde{f}_{u}\left(\frac{z}{|y|}\right)
\tilde{f}_{u}\left(\frac{z-y}{|y|}\right).
\end{equation}
Moreover, we will finally prove that the probability for
$(Z^{1},Z^{2})$ to avoid each other in the time interval $[1,
|y|^{3/2}(u\wedge 1)]$ is of order $\ln(|y| \vee e)^{-1}$ :
\begin{equation} \label{evite}
\ln(|y| \vee e) \mathbb{P}_{z,z} \bigg[ Z_{s}^{1} \neq Z_{s}^{2} \
\forall s \in [1,|y|^{3/2}(u \wedge 1)] \bigg] \leq L.
\end{equation}
Combining (\ref{bouley}) (\ref{markovy3/2}), (\ref{supf}) and
(\ref{evite}), we obtain \begin{eqnarray*}
 && \ln(|y|) |y|^{4}
\mathbb{P}_{z,z} \left[ Z_{s}^{1} \neq Z_{s}^{2} \ \forall s \in
[1,|y|^{2}u] , Z_{|y|^{2}u}^{1}=0, Z_{|y|^{2}u}^{2}=y   \right]\\
&& \leq L \tilde{f}_{u}\left(\frac{z}{|y|}\right)
\tilde{f}_{u}\left(\frac{z-y}{|y|}\right).
\end{eqnarray*}
 Claim \ref{eqinter} then follows from (\ref{sep}) and the above.
To complete the proof of Claim \ref{eqinter}, hence the one of Lemma
\ref{leminter}, it remains to establish (\ref{bouley}), (\ref{supf})
and (\ref{evite}). \PA {\em Proof of (\ref{bouley})} :
 As a consequence of (\ref{exp}),
\begin{eqnarray} \label{bouley1}
&& \ln(|y|)\mathbb{P}_{z,z} \!\bigg[ Z_{|y|^{3/2}(u\wedge 1)}^{1} \notin
B_{z,y}
 \mbox{ or } Z_{|y|^{3/2}(u\wedge 1)}^{2} \notin
B_{z,y} \bigg]
\nonumber\\
&& \leq 2 \kappa_{3} \ln(|y|) \exp\left(-\kappa_{4}
\frac{|y|^{1/8}}{\sqrt{u\wedge 1}}\left(\frac{|z|}{3|y|}
\vee 1 \right)\right) \nonumber \\
&& \leq L  \exp\left(-\kappa_{4} \frac{|y|^{1/8}}{2\sqrt{u\wedge
1}}\left(\frac{|z|}{3|y|} \vee 1 \right)\right),
\end{eqnarray}
where at the last line, we used the bound $\ln(|y|)
\exp\left(-\kappa_{4} |y|^{1/8}/2 \right) \leq L$. By studying
separately the cases $|z| \geq 3|y|$, $|z| \leq 3|y|$, and using the
fact that $\kappa_{4}/4 \geq \kappa_{2}$, we get
\begin{equation*}
 \exp\left(-\kappa_{4}
\frac{|y|^{1/8}}{2\sqrt{u\wedge 1}}\left(\frac{|z|}{3|y|} \vee 1
\right)\right)  \leq
\exp\left(-\kappa_{2}\frac{|z|}{4|y|\sqrt{u}}\right)
  \exp\left(-\kappa_{2}\frac{|z-y|}{4|y|\sqrt{u}}\right),
\end{equation*}
Furthermore, since  $t(u,y) \geq |y|^{2}u/2$,
we get from Lemma \ref{boundq} that
$$\sup \{ |y|^{2}
q_{t(u,y)} (x,x'), x \in \mathbb{Z}^{d}, x'\in \mathbb{Z}^{d} \}
\leq L u^{-1} $$ Hence, (\ref{bouley}) follows from using the Markov
property for the walks $Z^{1},Z^{2}$ at time $|y|^{3/2}(u\wedge 1)$
and combining the above remarks. \PA {\em Proof of (\ref{supf})} :
Lemma \ref{boundq} implies that
\begin{equation} \label{supq}
\sup_{(x_{1},x_{2}) \in B_{z,y}}
q_{t(u,y)}(x_{1}) q_{t(u,y)}(x_{2}-y) \leq 4|y|^{-4} \!\!
\sup_{(x_{1},x_{2}) \in B_{z,y}^{2}}
f_{u}\left(\frac{x_{1}}{|y|}\right)
f_{u}\left(\frac{x_{2}-y}{|y|}\right).
\end{equation}
The bound (\ref{supf}) in the case $|z| \geq 3|y|$ easily follows
from (\ref{supq}) and the fact that for $x_{1},x_{2}$ both in
$B_{z,y}$, we have
 $|x_{1}|\geq 2|z|/3$, $|x_{2}-y| \geq |z|/3 \geq |z-y|/4$. \PA
 Let us now suppose $|z| \leq 3|y|$, and recall that
we assumed $|y| \geq 6^8$, so that the balls $B(0,3|y|^{7/8})$ and
$B(y,3|y|^{7/8})$ are disjoint. Now, if $z \in B(0,2|y|^{7/8})$ we
easily see that for any $x_{2} \in B(z,|y|^{7/8})$, we have
$|x_{2}-y| \geq |z-y|/2 \geq |z|/2.$ It follows that
$$ \exp \left(-\frac{\kappa_{2}|x_{2}-y|}{|y|\sqrt{u}}\right) \leq
\exp \left(-\frac{\kappa_{2}|z-y|}{4|y|\sqrt{u}}\right) \exp
\left(-\frac{\kappa_{2}|z|}{4|y|\sqrt{u}}\right).$$ Assertion
(\ref{supq}) and the above imply
 (\ref{supf}) in the case $z \in B(0,2|y|^{7/8})$.
 We can use a similar argument to conclude in the case
$z \in B(y, 2|y|^{7/8})$.  At last, if $z$ satisfies $|z|\leq 3|y|$
but is not in any of the two aforementioned balls, then, for any
$x_{1}, x_{2} \in B(z,|y|^{7/8})$ we have $|x_{1}| \geq |z|/2$,
$|x_{2}-y| \geq |z-y|/2$, and (\ref{supf}) easily follows from
(\ref{supq}). \PA {\em Proof of (\ref{evite})} :
First note that under $\mathbb{P}_{z,z}$, $Z^{1}-Z^{2}$ has law
 $\mathbb{P}_{0}^{(2)}$.
Recall that the
notation $a^{(2)}(x)$, $G_{x}^{(2)}$, $H_{t}^{(2)}(x)$ have been introduced at the
beginning of the section.
From the simple bound $H_{t}^{(2)}(x) \geq G_{t}^{(2)}(0)^{-1}G_{t}^{(2)}(x)$, then
(\ref{a(x)}), we get
\begin{equation} \label{green}
 1-H^{(2)}_{|y|}(x) \leq
\frac{G^{(2)}_{|y|}(0)-G^{(2)}_{|y|}(x)}{G^{(2)}_{|y|}(0)} \leq
\frac{ \kappa_{5} a^{(2)}(x)}{\ln(|y|)}.
\end{equation}
From the Markov property for
$Z^{1}-Z^{2}$ at time $1$, we have
\begin{eqnarray*}
\!\!\!\!\!&&\!\!\!\!\!
\mathbb{P}_{z,z} \left[
Z_{s}^{1} \neq Z_{s}^{2} \ \forall s \in [1,|y|^{3/2}(u \wedge 1)] \right] \\
&=&   \mathbb{E}_{z,z} \left[
\sum_{x \in \mathbb{Z}^{d}} \mathbf{1}_{\{Z_{1}^{1}-Z_{1}^{2} = x\}}
\mathbf{1}_{ \{ Z_{s}^{1} - Z_{s}^{2} \neq 0 \
\forall s \in [1,|y|^{3/2}(u \wedge 1)]\} } \right] \\
&=& \sum_{x \in \mathbb{Z}^{d}} q_{2}(x)
\left( 1-\mathbb{P}_{x}^{(2)} \left( Z \mbox{ hits } 0
\mbox{ before time } |y|^{3/2}(u \wedge 1) \right) \right) \\
& \leq & \sum_{x \neq 0} q_{2}(x) \frac{\kappa_{5}
a^{(2)}(x)}{\ln\left(|y|^{3/2}(u \wedge 1)\right)}
\end{eqnarray*}
where we used (\ref{green}) at the last line. Since $u \geq |y|^{-1}$, we have
 $\ln\left(|y|^{3/2}(u \wedge 1)\right) \geq \ln(|y|)/2$. Hence, using Lemma
\ref{boundq} and the fact that $a^{(2)}(x)=O(\ln(|x|))$, we get
(\ref{evite}). \PA  As explained earlier on, this completes the
proof of Claim \ref{eqinter}, hence the one of Lemma \ref{leminter}.
$ \qquad  \Box $

\vskip0.2cm\PA Let us now complete the proof of Lemma
\ref{boundcoal}. We will apply Lemma \ref{leminter} to bound the
right-hand side of (\ref{start}). Fix $t\in [2c^{-2},T]$. Let us
consider the nonnegative functions $\Phi_{c}(u,z) =
f_{t-\frac{u|y|^{2}}{c^{2}}} \left(\frac{z}{c}-x\right)$. For $u \in
[\frac{1}{|y|^{2}}, \frac{2tc^{2}-1}{2|y|^{2}}]$ we have
$$\Phi_{c}(u,z) := \sup_{u-(2|y|^{2})^{-1} \leq r \leq u} \phi_{c}(r,z)
\leq \tilde{f}_{t-\frac{u|y|^{2}}{c^{2}}}\left( \frac{z}{c}-x \right).$$
Thus, from (\ref{start}) and Lemma \ref{leminter}, it follows that
 \begin{eqnarray} \label{decoupentc}
\!\!\!\!\!\!\!\!&\!\!\!& \ \ \ \ \ \
c^{d}\psi_{d}(|y|) \int_{2c^{-2}}^{T} dt \mathbb{P}_{0,y} \left[
1 \leq T_{1}  \leq c^{2}t-1 , Z_{c^{2}t}^{1}=cx \right]  \\
\!\!\!\!\!&\!\!\!&\!\!\!\! \leq   K|y|^{-d} \int_{2c^{-2}}^{T} \! dt
\sum_{z \in \mathbb{Z}^{d}}
 \int_{|y|^{-2}}^{\frac{tc^{2}}{|y|^{2}}-\frac{1}{2|y|^{2}}}\!\!\!\! du
\tilde{f}_{t-\frac{u|y|^{2}}{c^{2}}}\left( \frac{z}{c}-x \right)
\tilde{f}_{u}\left(\frac{z}{|y|}\right) \tilde{f}_{u}\left(\frac{z-y}{|y|}
\right)   \nonumber \\
\!\!\!\!\!&\!\!\!&\!\!\!\! \leq  K|y|^{-d} \int_{2c^{-2}}^{T} \! dt \!
\sum_{z \in \mathbb{Z}^{d}}
 \int_{|y|^{-2}}^{\frac{tc^{2}}{2|y|^{2}}}  \! du
\hat{f}_{t}\left( \frac{z}{c}-x \right)
  \tilde{f}_{u} \left(\frac{z}{|y|}\right)
\tilde{f}_{u} \left(\frac{z-y}{|y|} \right)
          \nonumber \\ \!\!\!\!\!&\!\!\!&\!\!\!\! \ +
K|y|^{-d} \!\!\int_{2c^{-2}}^{T}\! dt \sum_{z \in \mathbb{Z}^{d}}
 \int_{\frac{tc^{2}}{2|y|^{2}}}^{\frac{tc^{2}}{|y|^{2}}-\frac{1}{2|y|^{2}}} \!\!\! du
\tilde{f}_{t-\frac{u|y|^{2}}{c^{2}}}\left( \frac{z}{c}-x \right)
  \hat{f}_{\frac{tc^{2}}{|y|^{2}}}\left(\frac{z}{|y|}\right)
\hat{f}_{\frac{tc^{2}}{|y|^{2}}} \left(\frac{z-y}{|y|} \right)\!.
\nonumber
\end{eqnarray}
For convenience, let us define, for $z \in \mathbb{Z}^{d}$, $y \in \mathbb{Z}^{d} \setminus 0$,
$x \in c^{-1}\mathbb{Z}^{d}\setminus 0$ and $c \geq |x|^{-1} \vee 1$,
\begin{eqnarray*}
\!\!\!\!&\!\!\!&\!\!\!\! F_{1}(c,x,y,z):  =   \int_{2c^{-2}}^{T}\! dt
 \int_{|y|^{-2}}^{\frac{tc^{2}}{2|y|^{2}}} \! du
\hat{f}_{t}\left( \frac{z}{c}-x \right)
  \tilde{f}_{u} \left(\frac{z}{|y|}\right)
\tilde{f}_{u} \left(\frac{z-y}{|y|} \right)  \\
\!\!\!\!&\!\!\!&\!\!\!\! F_{2}(c,x,y,z):  = \int_{2|y|^{-2}}^{\frac{Tc^{2}}{|y|^{2}}}\! dt'
\int_{\frac{1}{2|c|^{2}}}^{\frac{t}{2}} \!\!\! du'
\tilde{f}_{u'}\left( \frac{z}{c}-x \right)
  \hat{f}_{t'}\left(\frac{z}{|y|}\right)
\hat{f}_{t'} \left(\frac{z-y}{|y|} \right).
\end{eqnarray*}
so that (\ref{decoupentc}) can be rewritten
\begin{eqnarray} \label{rhs}
\!\!\!\!&\!\!\!& \!\!\!\!
c^{d}\psi_{d}(|y|) \int_{2c^{-2}}^{T} dt \mathbb{P}_{0,y} \left[
1 \leq T_{1}  \leq c^{2}t-1 , Z_{c^{2}t}^{1}=cx \right]  \\
\!\!\!\!&\!\!\!&\!\!\!\! \leq
  K |y|^{-d} \sum_{z \in \mathbb{Z}^{d}}
F_{1}(c,x,y,z) + K|y|^{-d}
 \sum_{z \in \mathbb{Z}^{d}}  F_{2}(c,x,y,z). \nonumber
\end{eqnarray}
Thus,
completing the proof of Lemma \ref{boundcoal}, in the case $d=2$,
reduces to verify the bounds
\begin{equation} \label{F1d=2}
  |y|^{-2} \ln\left(\frac{c}{|y|}\vee e\right)^{-1} \sum_{z \in \mathbb{Z}^{2}}
F_{1}(c,x,y,z) \leq K \ln(|x|^{-1} \vee e),
\end{equation}
\begin{equation} \label{F2d=2}
 |y|^{-2} \ln\left(\frac{c}{|y|}\vee e\right)^{-1} \sum_{z \in\mathbb{Z}^{2}}
F_{2}(c,x,y,z) \leq K \ln(|x|^{-1}\vee e).
\end{equation}
Similarly, in the case $d=3$, in order to complete the proof
of Lemma \ref{boundcoal}, we need to establish that
\begin{equation} \label{F1d=3}
  |y|^{-3}  \sum_{z \in \mathbb{Z}^{3}}
F_{1}(c,x,y,z) \leq K (|x|^{-1} \vee 1),
\end{equation}
\begin{equation} \label{F2d=3}
 |y|^{-3} \sum_{z \in \mathbb{Z}^{3}}
F_{2}(c,x,y,z) \leq K (|x|^{-1}\vee 1).
\end{equation}
We first deal with the first term of the sum in the right-hand side of (\ref{rhs}). \\
\underline{Proof of (\ref{F1d=2}), (\ref{F1d=3})}:
From (\ref{integfoverR}), we obtain
\begin{eqnarray} \label{integF1:u}
 \int_{|y|^{-2}}^{\frac{tc^{2}}{2|y|^{2}}} du \tilde{f}_{u} \left(\frac{z}{|y|}\right)
\tilde{f}_{u} \left(\frac{z-y}{|y|} \right) \leq K \left( \frac{|z|}{|y|}\vee 1\right)^{2-2d}.
\end{eqnarray}
Then, from (\ref{integf}) and (\ref{integf2}), we easily get
\begin{eqnarray} \label{integF1:t}
\int_{2c^{-2}}^{T} dt \hat{f}_{t}(\frac{z}{c}-x) \leq \begin{cases}
K \psi_{d}(c) & \mbox{ if } z=cx, \\
K \psi_{d}\left(\left| \frac{z}{c}-x\right|^{-1} \right) & \mbox{ if } z \neq cx, \\
K \exp\left(-K'\left| \frac{z}{c}-x \right|\right) & \mbox{ if } \left| \frac{z}{c}-x\right| \geq \sqrt{T}.
\end{cases}
\end{eqnarray}
We then split $\mathbb{Z}^{d}$ into the following subsets
\begin{eqnarray*}
D_{0} \!\!&\!\! := \!\!& \!\! \{cx\},  \ D_{1} :=  \left(\mathbb{Z}^{d} \cap B(cx,c|x|/2)\right)\setminus D_{0},   \
D_{2} := \left(\mathbb{Z}^{d}\cap B(0,|y| \wedge 2c^{2}T) \right)\setminus D_{1}, \\
D_{3} \!\!&\!\! := \!\!&\!\!\left(\mathbb{Z}^{d}\cap B(0,|y| \vee 2c^{2}T) \right) \setminus (D_{2}\cup D_{1}),  \
D_{4} := \mathbb{Z}^{d} \setminus (D_{1}\cup D_{3}).
\end{eqnarray*}
We now combine the displays (\ref{integF1:u}), (\ref{integF1:t}), in order to
obtain bounds on $F_{1}(c,x,y,z)$ over the regions $D_{i}$, $0 \leq i \leq 4$.
We also use that, for $z \in D_{1}$, $|z| \geq Kc|x|$, while,
for $z \notin D_{1}$,
$\psi_{d}\left(\left|\frac{z}{c}-x\right|^{-1}\right) \leq K\psi_{d}(|x|^{-1})$. We have
\begin{equation} \label{boundF1}
 F_{1}(c,x,y,z) \leq K \times \begin{cases}
\psi_{d}(c)\left(\frac{c|x|}{|y|} \vee 1\right)^{2-2d} & \mbox{ if } z=cx, \\
\psi_{d}\left(\left| \frac{z}{c}-x\right|^{-1}\right)
 \left(\frac{c|x|}{|y|} \vee 1\right)^{2-2d} & \mbox{ if } z \in D_{1}, \\
  \psi_{d}(|x|^{-1}) & \mbox{ if } z \in D_{2}, \\
  \psi_{d}(|x|^{-1}) \left( \frac{|z|}{|y|} \vee 1 \right)^{2-2d}
& \mbox{ if } z \in D_{3}\cup D_{4}, \\
\exp \left(-K' \frac{|z|}{c} \right) |z|^{2-2d} |y|^{2d-2} & \mbox{ if } z \in D_{4}.
\end{cases}
\end{equation}
Then, observe that
\begin{eqnarray} \label{sumoverz}
&& \sum_{z \in D_{1}} \psi_{d}\left(\left| \frac{z}{c}-x\right|^{-1}\right) \leq
K c^{d-2} (c|x|)^{2}, \nonumber \\
&& |D_{2}| \leq K (|y|\wedge c^{2})^{d}, \nonumber\\
&& |D_{3}| \leq  K (|y|\vee c^{2})^{d}, \nonumber \\
&& \mbox { if } d=3, \ \  \sum_{z \in D_{3} \cup D_{4}} |z|^{-4} \leq K(|y| \wedge c^{2})^{-1}   \\
&& \mbox{ if } d=2, \ \ \begin{cases} & \sum_{z \in D_{3}} |z|^{-2}  \leq K \ln\left( \frac{c}{|y|} \vee e \right),  \\
 & \sum_{z \in D_{4}} |z|^{-2} \exp\left(-K'\frac{|z|}{c}\right) \leq K. \end{cases} \nonumber
\end{eqnarray}
Combining the bounds (\ref{boundF1}) and (\ref{sumoverz}), and doing some elementary
computations then leads to (\ref{F1d=2}), (\ref{F1d=3}).  \vskip0.2cm
\underline{Proof of (\ref{F2d=2}), (\ref{F2d=3})}:
From (\ref{integf}) and (\ref{integf2}), we obtain
\begin{eqnarray} \label{integF2:u}
 \int_{\frac{1}{2c^{2}}}^{\frac{t}{2}}
du' \tilde{f}_{u'}\left(\frac{z}{c}-x \right)  \leq \begin{cases}
K \psi_{d}(c) & \mbox{ if } z=cx, \\
K \psi_{d}\left(\left| \frac{z}{c}-x\right|^{-1} \right) & \mbox{ if } z \neq cx, \\
K \exp\left(-K'\left| \frac{z}{c}-x \right|\right) & \mbox{ if } \left| \frac{z}{c}-x\right| \geq \sqrt{T}.
\end{cases}
\end{eqnarray}
Also, from (\ref{integfoverR}),
\begin{eqnarray} \label{integF2:t}
 \int_{2|y|^{-2}}^{\frac{Tc^{2}}{|y|^{2}}} dt' \hat{f}_{t'} \left(\frac{z}{|y|}\right)
\hat{f}_{t'} \left(\frac{z-y}{|y|} \right) \leq K \left( \frac{|z|}{|y|}\vee 1\right)^{2-2d}.
\end{eqnarray}
Thus, the bounds in (\ref{boundF1}) remain true when replacing $F_{1}$ with $F_{2}$,
 and (\ref{F2d=2}), (\ref{F2d=3}) follow.
This completes the proof of Lemma \ref{boundcoal}. $\quad \Box$

\medskip
{\bf Acknowledgement} : This work is part of my PhD thesis, which I
did at Ecole Normale Sup\'erieure under the supervision of
Jean-Fran\c cois Le Gall. I wish to thank him heartily for the
patient guidance he provided me during this research. This article
owes much to his numerous advices and constant help.

\end{document}